\title{$SU(3)$-Goodman-de la Harpe-Jones subfactors and the realisation of $SU(3)$ modular invariants}
\author{
        David E. Evans and Mathew Pugh \\ \\
        School of Mathematics, \\
        Cardiff University, \\
        Senghennydd Road, \\
        Cardiff, CF24 4AG, \\
        Wales, U.K.
}
\date{\today}

\documentclass[12pt]{article}
\usepackage{amssymb}
\usepackage{graphicx}

\newtheorem{Def}{Definition}[section]
\newtheorem{Prop}[Def]{Proposition}
\newtheorem{Lemma}[Def]{Lemma}

\begin{document}
\maketitle

\begin{abstract}
We complete the realisation by braided subfactors, announced by Ocneanu, of all $SU(3)$-modular
invariant partition functions previously classified by Gannon.
\end{abstract}

\section{Introduction}

In \cite{goodman/de_la_harpe/jones:1989} Goodman, de la Harpe and Jones constructed a subfactor $B \subset C$ given by the embedding of the Temperley-Lieb algebra in the AF-algebra for an $SU(2)$ $ADE$ Dynkin diagram.
We will present an $SU(3)$ analogue of this construction, where we embed the $SU(3)$-Temperley-Lieb or Hecke algebra in an AF path algebra of the $SU(3)$ $\mathcal{ADE}$ graphs. Using this construction, we are able to realize all the $SU(3)$ modular invariants by subfactors.

The algebraic structure behind the integrable statistical mechanical $SU(N)$-models are the Hecke algebras $H_n(q)$ of type $A_{n-1}$, for $q \in \mathbb{C}$, since the Boltzmann weights lie in $(\bigotimes_{\mathbb{N}} M_N)^{SU(N)}$ or $(\bigotimes_{\mathbb{N}} M_N)^{SU(N)_q}$. The Hecke algebra $H_n(q)$ is the algebra generated by unitary operators $g_j$, $j=1,2,\ldots,n-1$, satisfying the relations
\begin{eqnarray}
(q^{-1} - g_j)(q + g_j) & = & 0, \\
g_i g_j & = & g_j g_i, \quad |i-j|>1,\\
g_i g_{i+1} g_i & = & g_{i+1} g_i g_{i+1}. \label{eqn:braiding_relation}
\end{eqnarray}
When $q=1$, the first relation becomes $g_j^2 = 1$, so that $H_n(1)$ reduces to the group ring of the symmetric, or permutation, group $S_n$, where $g_j$ represents a transposition $(j,j+1)$.
Writing $g_j = q^{-1} - U_j$ where $|q|=1$, and setting $\delta = q+q^{-1}$, these relations are equivalent to the self-adjoint operators $\mathbf{1}, U_1, U_2, \ldots, U_{n-1}$ satisfying the relations
\begin{center}
\begin{minipage}[b]{15cm}
 \begin{minipage}[b]{2cm}
  \begin{eqnarray*}
  \textrm{H1:}\\
  \textrm{H2:}\\
  \textrm{H3:}
  \end{eqnarray*}
 \end{minipage}
 \hspace{2cm}
 \begin{minipage}[b]{7cm}
  \begin{eqnarray*}
  U_i^2 & = & \delta U_i,\\
  U_i U_j & = & U_j U_i, \quad |i-j|>1,\\
  U_i U_{i+1} U_i - U_i & = & U_{i+1} U_i U_{i+1} - U_{i+1}.
  \end{eqnarray*}
 \end{minipage}
\end{minipage}
\end{center}

To any $\sigma$ in the permutation group $S_n$, decomposed into transpositions of nearest neighbours $\sigma = \prod_{i \in I_{\sigma}} \tau_{i,i+1}$, we associate the operator
$$g_{\sigma} = \prod_{i \in I_{\sigma}} g_i,$$
which is well defined because of the braiding relation (\ref{eqn:braiding_relation}).
Then the commutant of the quantum group $SU(N)_q$ is obtained from the Hecke algebra by imposing an extra condition, which is the vanishing of the $q$-antisymmetrizer
\begin{equation}\label{SU(N)q condition}
\sum_{\sigma \in S_{N+1}} (-q)^{|I_{\sigma}|} g_{\sigma} = 0.
\end{equation}
For $SU(2)$ it reduces to the Temperley-Lieb condition
\begin{equation}
U_i U_{i \pm 1} U_i - U_i = 0, \label{eqn:SU(2)q_condition}
\end{equation}
and for $SU(3)$ it is
\begin{equation} \label{eqn:SU(3)q_condition}
\left( U_i - U_{i+2} U_{i+1} U_i + U_{i+1} \right) \left( U_{i+1} U_{i+2} U_{i+1} - U_{i+1} \right) = 0.
\end{equation}

We will say that a family of operators $\{ U_m \}$ satisfy the $SU(3)$-Temperley-Lieb relations if they satisfy the Hecke relations H1-H3 and the extra condition (\ref{eqn:SU(3)q_condition}).
The Temperley-Lieb algebra has diagrammatic representations due to Kauffman \cite{kauffman:1987}.
There are similar diagrammatic representations for the $SU(3)$-Temperley-Lieb based on the  spider relations of Kuperberg, which we will exploit in a later sequel \cite{evans/pugh:2009iii, evans/pugh:2009iv} going into $SU(3)$-planar algebras. However, for our purposes here to construct $SU(3)$-Goodman-de la Harpe-Jones subfactors, it is enough to work algebraically.
We will embed the  $SU(3)$-Temperley-Lieb algebra in the path algebra  of the candidate nimrep graphs for the $SU(3)$ modular invariants, using the Boltzmann weights we constructed in \cite{evans/pugh:2009i}. This is with the exception of the graph $\mathcal{E}_4^{(12)}$ for which we did not derive the Ocneanu cells which permitted the derivation of the Boltzmann weights. However this is still enough to realise all $SU(3)$-modular invariants, and compute their nimrep graphs with the exception of $\mathcal{E}_4^{(12)}$ which we will do in this paper, after first outlining the theory of modular invariants from $\alpha$-induction in the next section.

\section{$\mathcal{ADE}$ Graphs}

We start with the $SU(3)$ modular invariants. The list below of all $SU(3)$ modular invariants was shown to be complete by Gannon \cite{gannon:1994}. Let $\mathcal{P}^{(n)} = \{ \mu=(\mu_1,\mu_2) \in \mathbb{Z}^2 | \mu_1,\mu_2 \geq 0; \mu_1 + \mu_2 \leq n-3 \}$. These $\mu$ are the admissible representations of the Ka\v{c}-Moody algebra $su(3)^{\wedge}$ at level $k = n-3$. We define the automorphism $A$ of order 3 on the weights $\mu \in \mathcal{P}^{(n)}$ by $A(\mu_1,\mu_2) = (n-3-\mu_1-\mu_2,\mu_1)$.

There are four infinite series of $SU(3)$ modular invariants: the identity (or diagonal) invariant at level $n-3$ is
\begin{equation}
Z_{\mathcal{A}^{(n)}} = \sum_{\mu \in P^{(n)}_+} |\chi_{\mu}|^2, \qquad n \geq 4, \label{Z(A)}
\end{equation}
and its orbifold invariant is given by
\begin{eqnarray}
Z_{\mathcal{D}^{(3k)}} & = & \frac{1}{3} \sum_{\stackrel{\mu \in P^{(3k)}_+}{\scriptscriptstyle{\mu_1 - \mu_2 \equiv 0 \mathrm{mod}3}}} |\chi_{\mu} + \chi_{A \mu} + \chi_{A^2 \mu}|^2, \qquad k \geq 2, \label{Z(D3k)} \\
Z_{\mathcal{D}^{(n)}} & = & \sum_{\mu \in P^{(n)}_+} \chi_{\mu} \chi_{A^{(n-3)(\mu_1-\mu_2)} \mu}^{\ast}, \qquad n \geq 5, \; n \not \equiv 0 \mathrm{mod} 3. \label{Z(Dk)}
\end{eqnarray}
Two other infinite series are given by their conjugate invariants. The conjugate invariant $Z_{\mathcal{A}^{(n)\ast}} = C$ and the conjugate orbifold invariants $Z_{\mathcal{D}^{(n)\ast}} = Z_{\mathcal{D}^{(n)}} C$ are
\begin{eqnarray}
Z_{\mathcal{A}^{(n)\ast}} & = & \sum_{\mu \in P^{(n)}_+} \chi_{\mu} \chi_{\overline{\mu}}^{\ast}, \qquad n \geq 4, \label{Z(Astar)} \\
Z_{\mathcal{D}^{(3k)\ast}} & = & \frac{1}{3} \sum_{\stackrel{\mu \in P^{(3k)}_+}{\scriptscriptstyle{\mu_1 - \mu_2 \equiv 0 \mathrm{mod}3}}} (\chi_{\mu} + \chi_{A \mu} + \chi_{A^2 \mu}) (\chi_{\overline{\mu}}^{\ast} + \chi_{\overline{A \mu}}^{\ast} + \chi_{\overline{A^2 \mu}}^{\ast}), \qquad k \geq 2, \label{Z(D3kstar)} \\
Z_{\mathcal{D}^{(n)\ast}} & = & \sum_{\mu \in P^{(n)}_+} \chi_{\mu} \chi_{\overline{A^{(n-3)(\mu_1-\mu_2)} \mu}}^{\ast}, \qquad n \geq 5, \; n \not \equiv 0 \mathrm{mod} 3. \label{Z(Dkstar)}
\end{eqnarray}
There are also exceptional invariants, i.e. invariants which are not diagonal, orbifold, or their conjugates:
\begin{eqnarray}
Z_{\mathcal{E}^{(8)}} & = & |\chi_{(0,0)}+\chi_{(2,2)}|^2 + |\chi_{(0,2)}+\chi_{(3,2)}|^2 + |\chi_{(2,0)}+ \chi_{(2,3)}|^2 + |\chi_{(2,1)}+\chi_{(0,5)}|^2 \nonumber \\
& & \;\; + |\chi_{(3,0)}+\chi_{(0,3)}|^2 + |\chi_{(1,2)}+\chi_{(5,0)}|^2, \label{Z(E8)} \\
Z_{\mathcal{E}^{(8)\ast}} & = & |\chi_{(0,0)}+\chi_{(2,2)}|^2 + (\chi_{(0,2)}+\chi_{(3,2)})(\chi_{(2,0)}^{\ast}+\chi_{(2,3)}^{\ast}) \nonumber \\
& & \;\; + (\chi_{(2,0)}+\chi_{(2,3)})(\chi_{(0,2)}^{\ast}+\chi_{(3,2)}^{\ast}) + (\chi_{(2,1)}+\chi_{(0,5)})(\chi_{(1,2)}^{\ast}+\chi_{(5,0)}^{\ast}) \nonumber \\
& & \;\; + |\chi_{(3,0)}+\chi_{(0,3)}|^2 + (\chi_{(1,2)}+\chi_{(5,0)})(\chi_{(2,1)}^{\ast}+\chi_{(0,5)}^{\ast}), \label{Z(E8star))} \\
Z_{\mathcal{E}^{(12)}} & = & |\chi_{(0,0)}+\chi_{(0,9)}+\chi_{(9,0)}+\chi_{(4,4)}+\chi_{(4,1)}+\chi_{(1,4)}|^2 \nonumber \\
& & \;\; + 2 |\chi_{(2,2)}+\chi_{(2,5)}+\chi_{(5,2)}|^2, \label{Z(E1_12)} \\
Z_{\mathcal{E}_{MS}^{(12)}} & = & |\chi_{(0,0)}+\chi_{(0,9)}+\chi_{(9,0)}|^2 + |\chi_{(2,2)}+\chi_{(2,5)}+\chi_{(5,2)}|^2 + 2 |\chi_{(3,3)}|^2 \nonumber \\
& & \;\; + |\chi_{(0,3)}+\chi_{(6,0)}+ \chi_{(3,6)}|^2 + |\chi_{(3,0)}+\chi_{(0,6)}+\chi_{(6,3)}|^2 + |\chi_{(4,4)}+\chi_{(4,1)}+\chi_{(1,4)}|^2 \nonumber \\
& & \;\; + (\chi_{(1,1)}+\chi_{(1,7)}+\chi_{(7,1)})\chi_{(3,3)}^{\ast} + \chi_{(3,3)}(\chi_{(1,1)}^{\ast}+\chi_{(1,7)}^{\ast}+\chi_{(7,1)}^{\ast}), \label{Z(E5_12)} \\
Z_{\mathcal{E}_{MS}^{(12)\ast}} & = & |\chi_{(0,0)}+\chi_{(0,9)}+\chi_{(9,0)}|^2 + |\chi_{(2,2)}+\chi_{(2,5)}+\chi_{(5,2)}|^2 + 2 |\chi_{(3,3)}|^2 \nonumber \\
& & \;\; + (\chi_{(0,3)}+\chi_{(6,0)}+ \chi_{(3,6)})(\chi_{(3,0)}^{\ast}+\chi_{(0,6)}^{\ast}+\chi_{(6,3)}^{\ast}) \nonumber \\
& & \;\; + (\chi_{(3,0)}+\chi_{(0,6)}+\chi_{(6,3)})(\chi_{(0,3)}^{\ast}+\chi_{(6,0)}^{\ast}+ \chi_{(3,6)}^{\ast}) + |\chi_{(4,4)}+\chi_{(4,1)}+\chi_{(1,4)}|^2 \nonumber \\
& & \;\; + (\chi_{(1,1)}+\chi_{(1,7)}+\chi_{(7,1)})\chi_{(3,3)}^{\ast} + \chi_{(3,3)}(\chi_{(1,1)}^{\ast}+\chi_{(1,7)}^{\ast}+\chi_{(7,1)}^{\ast}), \label{Z(E4_12)} \\
Z_{\mathcal{E}^{(24)}} & = & |\chi_{(0,0)}+\chi_{(4,4)}+\chi_{(6,6)}+\chi_{(10,10)}+\chi_{(21,0)}+\chi_{(0,21)}+\chi_{(13,4)}+\chi_{(4,13)} \nonumber \\
& & \;\; +\chi_{(10,1)}+\chi_{(1,10)}+\chi_{(9,6)}+\chi_{(6,9)}|^2 \nonumber \\
& & \;\; + |\chi_{(15,6)}+\chi_{(6,15)}+\chi_{(15,0)}+\chi_{(0,15)}+\chi_{(10,7)}+\chi_{(7,10)}+\chi_{(10,4)} \nonumber \\
& & \;\; \quad +\chi_{(4,10)}+\chi_{(7,4)}+\chi_{(4,7)}+\chi_{(6,0)}+\chi_{(0,6)}|^2, \label{Z(E24)}
\end{eqnarray}
where $Z_{\mathcal{E}^{(12)}}$ and $Z_{\mathcal{E}^{(24)}}$ are self-conjugate, $Z_{\mathcal{E}^{(8)\ast}} = Z_{\mathcal{E}^{(8)}} C$ and $Z_{\mathcal{E}_{MS}^{(12)\ast}} = Z_{\mathcal{E}_{MS}^{(12)}} C$.

The modular invariants arising from $SU(3)_k$ conformal embeddings are (see \cite{evans:2003}):
\begin{itemize}
\item $\mathcal{D}^{(6)}$: $SU(3)_3 \subset SO(8)_1$, also realised as an orbifold $SU(3)_3 / \mathbb{Z}_3$,
\item $\mathcal{E}^{(8)}$: $SU(3)_5 \subset SU(6)_1$, plus its conjugate,
\item $\mathcal{E}^{(12)}$: $SU(3)_9 \subset (\mathrm{E}_6)_1$,
\item $\mathcal{E}^{(24)}$: $SU(3)_{21} \subset (\mathrm{E}_7)_1$.
\end{itemize}

The Moore-Seiberg invariant $\mathcal{E}_{MS}^{(12)}$ \cite{moore/seiberg:1989}, an automorphism of the orbifold invariant $\mathcal{D}^{(12)} = SU(3)_9 / \mathbb{Z}_3$, is the $SU(3)$ analogue of the $E_7$ invariant for $SU(2)$, which is an automorphism of the orbifold invariant $D_{10} = SU(2)_{16}/\mathbb{Z}_2$ (see Section 5.3 of \cite{bockenhauer/evans/kawahigashi:2000} for a realisation by a braided subfactor).

In the statistical mechanical models underlying this theory, the vertices and edges of the underlying graph are used to describe bonds on a two dimensional lattice, together with some Hamiltonian or family of Boltzmann weights. In the conformal field theory, or subfactor theory, the vertices of the graph appear as primary fields or endomorphisms of a type III factor.

\begin{figure}[tb]
\begin{center}
\includegraphics[width=70mm]{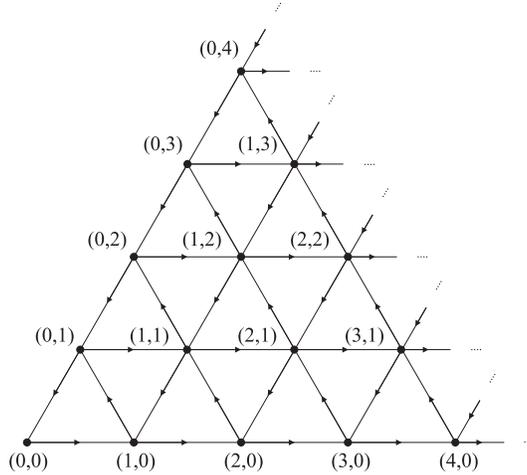}\\
 \caption{The infinite graph $\mathcal{A}^{(\infty)}$} \label{Fig:SU(3)-A(infty)}
\end{center}
\end{figure}

The simplest case of the diagonal invariant only involves the Verlinde algebra, whose fusion rules are determined by the graph $\mathcal{A}^{(n)}$.
The infinite graph $\mathcal{A}^{(\infty)}$ is illustrated in Figure \ref{Fig:SU(3)-A(infty)}, whilst for finite $n$, the graphs $\mathcal{A}^{(n)}$ are the subgraphs of $\mathcal{A}^{(\infty)}$, given by all the vertices $(\lambda_1, \lambda_2)$ such that $\lambda_1 + \lambda_2 \leq n-3$, and all the edges in $\mathcal{A}^{(\infty)}$ which connect these vertices.

The Verlinde algebra of $SU(3)$ at level $k=n-3$ will be represented by a finite system ${}_N \mathcal{X}_N$ of irreducible inequivalent endomorphisms of a type III factor $N$ \cite{wassermann:1998} which possesses a non-degenerate braiding, with unitary operator $\varepsilon(\lambda,\mu)$ intertwining $\lambda \mu$  and $\mu \lambda$ , called a braiding operator, which satisfy the Braiding Fusion Equations \cite[Def. 2.2]{bockenhauer/evans/kawahigashi:1999}.
For every braiding $\varepsilon^+ \equiv \varepsilon$ there is an opposite braiding $\varepsilon^-$ obtained by reversing the crossings.
If we have an inclusion $\iota:N \hookrightarrow M$ of type III factors together with a non-degenerately braided finite system ${}_N \mathcal{X}_N$ such that the dual canonical endomorphism $\theta = \overline{\iota} \iota $ decomposes as a sum of elements of ${}_N \mathcal{X}_N$ then we call $N \subset M$ a braided subfactor. The $\alpha$-induced morphisms $\alpha_{\lambda}^{\pm} \in \mathrm{End}(M)$, which extend $\lambda \in {}_N \mathcal{X}_N$, are defined by the Longo-Rehren formula \cite{longo/rehren:1995} $\alpha_{\lambda}^{\pm} = \overline{\iota}^{-1} \circ \mathrm{Ad}(\varepsilon^{\pm}(\lambda,\theta)) \circ \lambda \circ \overline{\iota}$,
A coupling matrix $Z$ can be defined \cite{bockenhauer/evans/kawahigashi:1999} by $Z_{\lambda,\mu} = \langle \alpha_{\lambda}^+, \alpha_{\mu}^- \rangle$, where $\lambda, \mu \in {}_N \mathcal{X}_N$, normalized so that $Z_{0,0} = 1$. By \cite{bockenhauer/evans:2000, evans:2003} this matrix $Z$ commutes with the modular $S$- and $T$-matrices, and therefore $Z$ is a modular invariant.
The right action of the $N$-$N$ system ${}_N \mathcal{X}_N$ on the $M$-$N$ system ${}_M \mathcal{X}_N$ yields a representation of the Verlinde algebra or a nimrep $G_{\lambda}$, of the original $N$-$N$ fusion rules,
i.e. a matrix representation where all the matrix entries are non-negative integers.
These nimreps give multiplicity graphs associated to the modular invariants (or at least associated to the inclusion, as a modular invariant may be represented by wildly differing inclusions).
The matrix $G_{\nu}$ has spectrum $S_{\lambda,\nu}/S_{\lambda,0}$ with multiplicity ${Z_{\lambda,\lambda}}$.
In particular the spectrum of the nimrep is determined by the diagonal part of the modular invariant and provides an automatic connection between the modular invariant and fusion graphs, which in the $SU(2)$ case reduces to the classification by Capelli-Itzykson-Zuber \cite{cappelli/itzykson/zuber:1987ii} of modular invariants by $ADE$ graphs.
As $M$-$N$ sectors cannot be multiplied among themselves there is no associated fusion rule algebra to decompose.
Nevertheless, when chiral locality does hold \cite{bockenhauer/evans:1999ii, bockenhauer/evans:2000} the nimrep graph ${}_M \mathcal{X}_N$ can be canonically identified with both chiral graphs  ${}_M \mathcal{X}_M^{\pm}$, the systems induced by the images of $\alpha$-induction, by $\beta \mapsto \beta \circ \iota$, $\beta \in {}_M \mathcal{X}_M^{\pm}$.

The question then arises whether or not every $SU(3)$ modular invariant can be realised by a subfactor.
This was claimed and announced by Ocneanu \cite{ocneanu:2002} in his bimodule setting. Most of these invariants  are understood in the literature.
Feng Xu \cite{xu:1998} (see also \cite{bockenhauer/evans:1998, bockenhauer/evans:1999i, bockenhauer/evans:1999ii}) looked at the conformal embedding invariants in the loop group setting of \cite{wassermann:1998}, taking $\alpha$-induction as the principal tool. These conformal inclusions are local or type I.
In particular, the chiral graphs for the $\mathcal{D}^{(6)}$, $\mathcal{E}^{(8)}$, $\mathcal{E}^{(12)}$ and $\mathcal{E}^{(24)}$ $SU(3)$ invariants were computed. Since these inclusions are type I, the chiral graphs coincide with their nimreps with corresponding graphs $\mathcal{D}^{(6)}$, $\mathcal{E}^{(8)}$, $\mathcal{E}_1^{(12)}$ and $\mathcal{E}^{(24)}$ respectively. These graphs are illustrated in Figures 10, 13, 14 and 16 of \cite{evans/pugh:2009i} respectively. Note that by the spectral theory of nimreps developed in \cite{bockenhauer/evans/kawahigashi:1999, bockenhauer/evans/kawahigashi:2000} and described above, these graphs and the other candidate graphs of di Francesco and Zuber will now automatically have spectra described by the diagonal part of the modular invariant.

The infinite series of orbifold invariants $\mathcal{D}^{(3k)}$ were considered by B\"ockenhauer and Evans in \cite{bockenhauer/evans:1999i}, yielding nimreps which produce the graphs $\mathcal{D}^{(3k)}$, which are the $\mathbb{Z}_3$-orbifolds of the graphs $\mathcal{A}^{(3k)}$.
B\"ockenhauer and Evans \cite{bockenhauer/evans:1999i} produced a method for analysing conjugates of conformal embedding invariants by taking an orbifold of the extended system of the level one theory of the ambient group.
In \cite{bockenhauer/evans:2001}, B\"ockenhauer and Evans realised all modular invariants for cyclic $\mathbb{Z}_n$ theories, in particular charge conjugation.
The conformal embedding modular invariant $\mathcal{E}^{(8)}$: $SU(3)_5 \subset SU(6)_1$ produces the $\mathcal{E}^{(8)}$ invariant and the nimrep graph $\mathcal{E}^{(8)}$. Then
taking the extension $SU(6)_1 \subset SU(6)_1 \rtimes  \mathbb{Z}_3$ describes charge conjugation on
the cyclic $\mathbb{Z}_6$ system for $SU(6)_1$. Then the inclusion $SU(3)_5 \subset SU(6)_1 \rtimes \mathbb{Z}_3$ produces its orbifold $\mathcal{E}^{(8)} / \mathbb{Z}_3$ for the conjugate modular invariant (see Figure \ref{fig:fig-GHJ-14}).
This procedure could be used to understand and realise $SU(3)_9 \subset (\mathit{E}_6)_1$, with two nimreps.
One was $\mathcal{E}_1^{(12)}$ through of course the $SU(3)_9 \subset (\mathit{E}_6)_1$ standard conformal embedding,
and another the orbifold $\mathcal{E}_2^{(12)} = {\mathcal{E}_1^{(12)}}/\mathbb{Z}_3$ obtained from the subfactor $SU(3)_9 \subset (\mathit{E}_6)_1 \rtimes \mathbb{Z}_3$.
The extension $(\mathit{E}_6)_1 \subset (\mathit{E}_6)_1 \rtimes \mathbb{Z}_3$ describes charge conjugation on the cyclic $\mathbb{Z}_6$ system for $(\mathit{E}_6)_1$. The conformal embedding inclusion is always local and so type I but its orbifold here is not local, so this particular modular invariant $\mathcal{E}^{(12)}$ is type I
for one subfactor realisation and type II for another, $\mathcal{E}_1^{(12)}$ and its $\mathbb{Z}_3$-orbifold $\mathcal{E}_2^{(12)}$ (see Figure \ref{fig:fig-GHJ-9}).

We now realise the remaining $SU(3)$ modular invariants $\mathcal{A}^{\ast}$, $\mathcal{D}^{\ast}$ and $\mathcal{E}_{MS}^{(12)}$ by subfactors, using an $SU(3)$ analogue of the Goodman-de la Harpe-Jones construction of a subfactor, where we embed the $SU(3)$-Temperley-Lieb or Hecke algebra in an AF path algebra of the $SU(3)$ $\mathcal{ADE}$ graphs.
These subfactors yield nimreps which produce the graphs $\mathcal{A}^{(n)\ast}$, $\mathcal{D}^{(n)\ast}$, $\mathcal{E}_{5}^{(12)}$ respectively (see Figures \ref{fig:fig-GHJ-12}, \ref{fig:fig-GHJ-13}, \ref{fig:fig-GHJ-11} respectively). We can also realize the conjugate invariant of the Moore-Seiberg invariant $\mathcal{E}_{MS}^{(12)}$ by a subfactor, since this is now a product of two modular invariants (the Moore-Seiberg and conjugate) which can both be realised by subfactors, and so by \cite[Theorem 3.6]{evans/pinto:2003} their product is also realised by an inclusion.
However, we have not yet been able to compute its nimrep as we have been unable to determine the cells for the graph $\mathcal{E}_4^{(12)}$ which would enable a direct computation of the desired nimrep graph using the $SU(3)$-Goodman-de la Harpe-Jones subfactor, or alternatively, compute the nimrep in the alternative inclusion given by the braided product of the Moore-Seiberg inclusion  and the conjugate inclusion.

Almost all the $\mathcal{ADE}$ graphs mentioned above were proposed by di Francesco and Zuber \cite{di_francesco/zuber:1990} by looking for graphs whose spectrum reproduced the diagonal part of the modular invariant, aided to some degree by first listing the graphs and spectra of fusion graphs of the finite subgroups of $SU(3)$. At that time, they proposed looking for 3-colourable graphs.
They succeeded, for $SU(3)$, in finding graphs and nimreps for the orbifold invariants, and the exceptional invariants (with three candidates for the conformal embedding $SU(3)_9 \subset (\mathit{E}_6)_1$ invariant). All these graphs were three-colourable, and they conjectured this to be the case for all $SU(3)$ modular invariants.
B\"ockenhauer and Evans \cite{bockenhauer:1999} understood that nimrep graphs for the conjugate $SU(3)$ modular invariants were not three colourable.
This was also realised simultaneously by Behrend, Pearce, Petkova and Zuber \cite{behrend/pearce/petkova/zuber:2000} and Ocneanu \cite{ocneanu:2002}.
Indeed Ocneanu announced in Bariloche \cite{ocneanu:2002} that all $SU(3)$ modular invariants can be realised by subfactors, and the classification of their associated nimreps. He ruled out the third candidate $\mathcal{E}_{3}^{(12)}$ for the $\mathcal{E}^{(12)}$ modular invariant by asserting that it did not support a valid cell system.
This graph was ruled out as a natural candidate in Section 5.2 of \cite{evans:2002}.

We now list the $\mathcal{ADE}$ graphs: four infinite series of graphs $\mathcal{A}^{(n)}$, $\mathcal{D}^{(n)}$, $\mathcal{A}^{(n)\ast}$ and $\mathcal{D}^{(n)\ast}$, $n \leq \infty$, and seven exceptional graphs $\mathcal{E}^{(8)}$, $\mathcal{E}^{(8)\ast}$, $\mathcal{E}_1^{(12)}$, $\mathcal{E}_2^{(12)}$, $\mathcal{E}_4^{(12)}$, $\mathcal{E}_5^{(12)}$ and $\mathcal{E}^{(24)}$.
We note that all the graphs are three-colourable, except for the graphs $\mathcal{D}^{(n)}$, $n \not \equiv 0 \textrm{ mod } 3$, $\mathcal{A}^{(n)\ast}$, $n \leq \infty$, and $\mathcal{E}^{(8)\ast}$. For the $\mathcal{A}$ graphs, the vertices are labelled by Dynkin labels $(\lambda_1,\lambda_2)$, $\lambda_1, \lambda_2 \geq 0$. We define the colour of a vertex $(\lambda_1,\lambda_2)$ of $\mathcal{A}^{(n)}$, $n < \infty$, to be $\lambda_1 - \lambda_2 \textrm{ mod } 3$. There is a natural conjugation on the graph defined by $\overline{(\lambda_1,\lambda_2)} = (\lambda_2,\lambda_1)$ for all $\lambda_1, \lambda_2 \geq 0$. This conjugation interchanges the vertices of colour 1 with those of colour 2, but leaves the set of all vertices of colour 0 invariant.
For all the other three-colourable graphs there is also a conjugation. The vertices of these graphs are coloured such that the conjugation again leaves the set of all vertices of colour 0 invariant. We use the convention that the edges on the graph are always from a vertex of colour $j$ to a vertex of colour $j+1$ (mod 3).
For the non-three-colourable graphs, we will not distinguish between the colour of vertices, so that all the vertices have colour $j$ for any $j \in \{ 1,2,3 \}$.
In this paper we will consider the finite graphs, i.e. $\mathcal{A}^{(n)}$, $\mathcal{D}^{(n)}$, $\mathcal{A}^{(n)\ast}$ and $\mathcal{D}^{(n)\ast}$, $n < \infty$, and the exceptional $\mathcal{E}$ graphs.

The figures for the complete list of the $\mathcal{ADE}$ graphs are given in \cite{behrend/pearce/petkova/zuber:2000, evans/pugh:2009i}.

\section{Ocneanu cells}

We will construct a representation of a Hecke algebra in the path algebra of an $\mathcal{ADE}$ graph. For more details on path algebras see \cite{evans/kawahigashi:1998}. This construction is not as straightforward as for $SU(2)$ where one only needs the Perron-Frobenius eigenvector for the $ADE$ Dynkin diagram.

The McKay graph $\mathcal{G}$ of $SU(3)$ is made of triangles, which are paths of length 3 on the graph such that the start and end vertices are the same. This corresponds to the fact that the fundamental representation $\rho$, which along with its conjugate representation $\overline{\rho}$ generates the irreducible representations of $SU(3)$, satisfies $\rho \otimes \rho \otimes \rho \ni \mathbf{1}$. To every triangle on $\mathcal{G}$ one can assign a complex number, called an Ocneanu cell.
More details are given in \cite{evans/pugh:2009i}.

These cells are axiomatized in the context of an arbitrary graph $\mathcal{G}$ whose adjacency matrix has Perron-Frobenius eigenvalue $[3] = [3]_q$, although in practice it will be any one of the $\mathcal{ADE}$ graphs. Here the quantum number $[m]_q$ is defined by $[m]_q = (q^m - q^{-m})/(q - q^{-1})$. We will frequently denote the quantum number $[m]_q$ simply by $[m]$, for $m \in \mathbb{N}$. Now $[3]_q = q^2 + 1 + q^{-2}$, so that $q$ is easily determined from the eigenvalue of $\mathcal{G}$. The quantum number $[2] = [2]_q$ is then simply $q + q^{-1}$.
If $\mathcal{G}$ is an $\mathcal{ADE}$ graph, the Coxeter number $n$ of $\mathcal{G}$ is the number in parentheses in the notation for the graph $\mathcal{G}$, e.g. the exceptional graph $\mathcal{E}^{(8)}$ has Coxeter number 8, and $q = e^{\pi i/n}$.

We define a type I frame in an arbitrary $\mathcal{G}$ to be a pair of edges $\alpha$, $\alpha'$ which have the same start and endpoint. A type II frame will be given by four edges $\alpha_i$, $i=1,2,3,4$, such that $s(\alpha_1) = s(\alpha_4)$, $s(\alpha_2) = s(\alpha_3)$, $r(\alpha_1) = r(\alpha_2)$ and $r(\alpha_3) = r(\alpha_4)$.

\begin{Def}[\cite{ocneanu:2002}]\label{cell_system}
Let $\mathcal{G}$ be an arbitrary graph with Perron-Frobenius eigenvalue $[3]$ and Perron-Frobenius eigenvector $(\phi_i)$. A \textbf{cell system} $W$ on $\mathcal{G}$ is a map that associates to each oriented triangle $\triangle_{ijk}^{(\alpha \beta \gamma)}$ in $\mathcal{G}$ a complex number $W \left( \triangle_{ijk}^{(\alpha \beta \gamma)} \right)$ with the following properties:

$(i)$ for any type I frame \includegraphics[width=16mm]{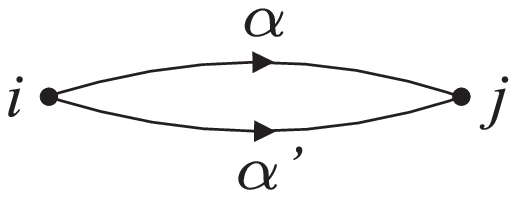} in $\mathcal{G}$ we have
\begin{flushright}
\begin{minipage}[b]{14cm}
 \begin{minipage}[t]{9cm}
  \mbox{} \\
 \includegraphics[width=80mm]{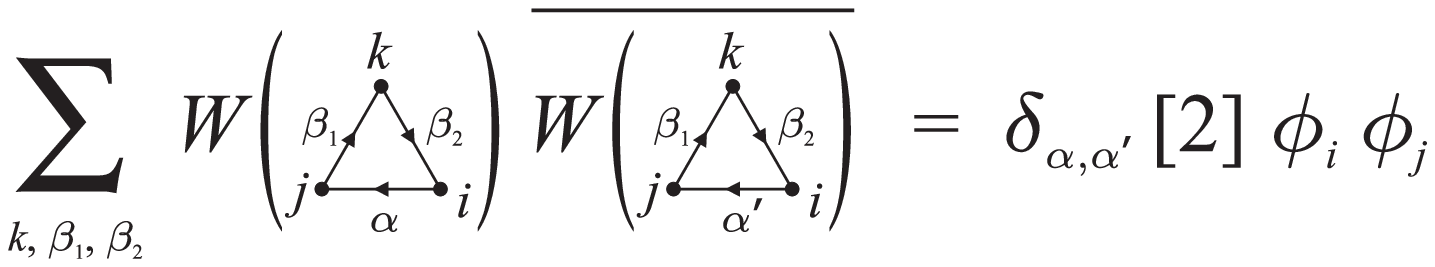}
 \end{minipage}
 \hfill
 \begin{minipage}[t]{1.5cm}
  \mbox{} \\
  \hfill
  \parbox[t]{7mm}{\begin{eqnarray}\label{eqn:typeI_frame}\end{eqnarray}}
 \end{minipage}
\end{minipage}
\end{flushright}

$(ii)$ for any type II frame \includegraphics[width=30mm]{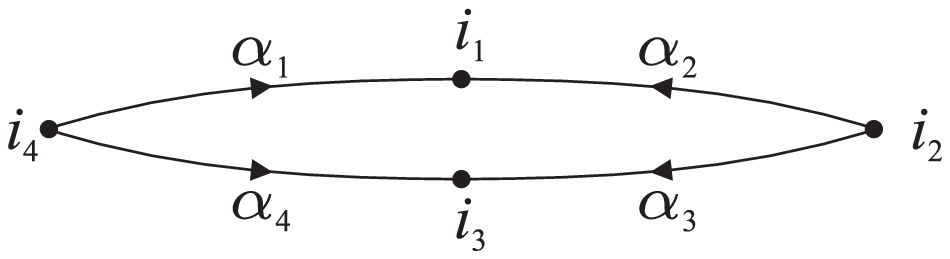} in $\mathcal{G}$ we have
\begin{flushright}
\begin{minipage}[b]{14cm}
 \begin{minipage}[b]{9cm}
  \mbox{} \\
 \includegraphics[width=90mm]{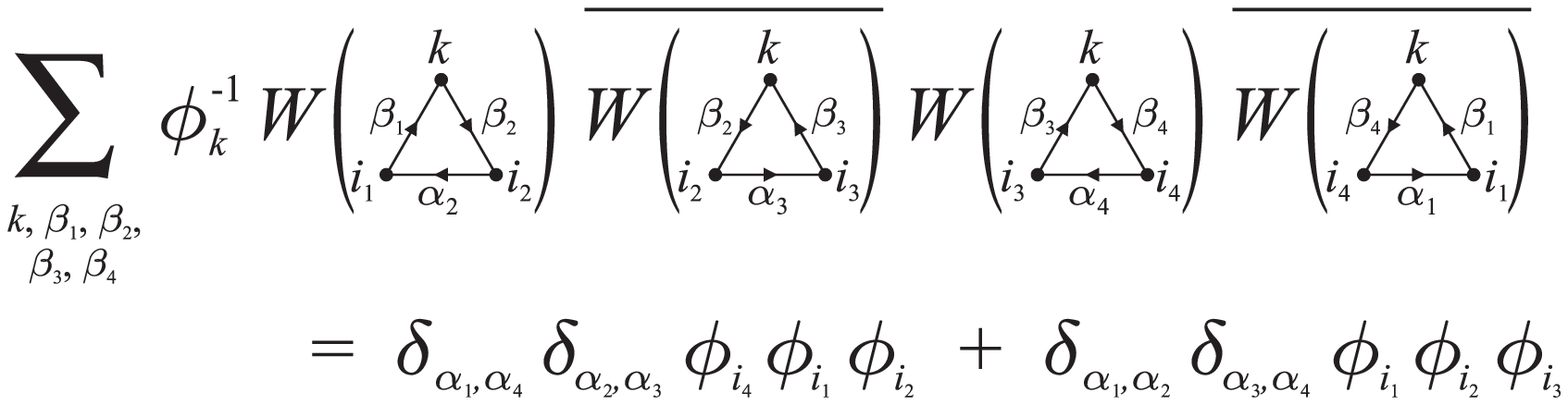}
 \end{minipage}
 \hfill
 \begin{minipage}[b]{1.5cm}
  \mbox{} \\
  \hfill
  \parbox[b]{7mm}{\begin{eqnarray}\label{eqn:typeII_frame}\end{eqnarray}}
 \end{minipage}
\end{minipage}
\end{flushright}
\end{Def}

Ocneanu cells for the $\mathcal{ADE}$ graphs were constructed in \cite{evans/pugh:2009i}, with the exception of the graph $\mathcal{E}_4^{(12)}$.
Using these cells we define the connection
$$X^{\rho_1,\rho_2}_{\rho_3,\rho_4} =
\begin{array}{c}
l \; \stackrel{\rho_1}{\longrightarrow} \; i \\
\scriptstyle \rho_3 \textstyle \big\downarrow \qquad \big\downarrow \scriptstyle \rho_2 \\
k \; \stackrel{\textstyle\longrightarrow}{\scriptstyle{\rho_4}} \; j
\end{array}$$
for the $\mathcal{ADE}$ graph $\mathcal{G}$ by
\begin{equation} \label{Connection_using_weights_W}
X^{\rho_1,\rho_2}_{\rho_3,\rho_4} = q^{\frac{2}{3}} \delta_{\rho_1,\rho_3} \delta_{\rho_2, \rho_4} - q^{-\frac{1}{3}} \mathcal{U}^{\rho_1,\rho_2}_{\rho_3,\rho_4},
\end{equation}
where $\mathcal{U}^{\rho_1,\rho_2}_{\rho_3,\rho_4}$ is given by the representation of the Hecke algebra, and is defined by
\begin{equation} \label{eqn:HeckeRep}
\mathcal{U}^{\rho_1,\rho_2}_{\rho_3,\rho_4} = \sum_{\lambda} \phi_{s(\rho_1)}^{-1} \phi_{r(\rho_2)}^{-1} W(\triangle_{j,l,k}^{(\lambda, \rho_3, \rho_4)}) \overline{W(\triangle_{j,l,i}^{(\lambda, \rho_1, \rho_2)})}.
\end{equation}
A representation $\mathcal{U}$ of the Hecke algebra corresponds to a picture
\begin{center}
\includegraphics[width=7mm]{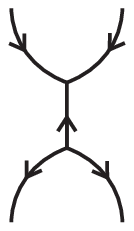}
\end{center}
in the $A_2$ web space. It will be proved in \cite{evans/pugh:2009iii} that a diagrammatic algebra generated by these pictures indeed gives a representation of the Hecke algebra. More details on the relation between the $A_2$ web space of Kuperberg and the Ocneanu cells are given in \cite{evans/pugh:2009i}.

The above connection corresponds to the natural braid generator $g_i$, up to a choice of phase.
It was claimed in \cite{ocneanu:2000ii} and proven in \cite{evans/pugh:2009i} that the connection satisfies the unitarity property of connections
\begin{equation} \label{eqn:unitarity_property_of_connections}
\sum_{\rho_3,\rho_4} X^{\rho_1,\rho_2}_{\rho_3,\rho_4} \; \overline{X^{\rho_1',\rho_2'}_{\rho_3,\rho_4}} = \delta_{\rho_1, \rho_1'} \delta_{\rho_2, \rho_2'},
\end{equation}
and the Yang-Baxter equation
\begin{equation} \label{eqn:YBE}
\sum_{\sigma_1, \sigma_2, \sigma_3} X^{\sigma_1,\sigma_2}_{\rho_1,\rho_2} \; X^{\rho_3,\rho_4}_{\sigma_1,\sigma_3} \; X^{\sigma_3,\rho_5}_{\sigma_2,\rho_6} = \sum_{\sigma_1, \sigma_2, \sigma_3} X^{\rho_3,\sigma_2}_{\rho_1,\sigma_1} \; X^{\sigma_1,\sigma_3}_{\rho_2,\rho_6} \; X^{\rho_4,\rho_5}_{\sigma_2,\sigma_3},
\end{equation}
provided that the cells $W(\triangle)$ satisfy (\ref{eqn:typeI_frame}), (\ref{eqn:typeII_frame}).

\section{General construction}

In this section we will construct the $SU(3)$-Goodman-de la Harpe-Jones subfactors. We first present some results that will be needed for this construction.

Let $U_1, U_2, \ldots U_{m-1}$ be operators which satisfy H1-H3 with parameter $\delta$. We let
\begin{equation} \label{Def:F_i}
F_i := U_i U_{i+1} U_i - U_i = U_{i+1} U_i U_{i+1} - U_{i+1},
\end{equation}
for $i=1,2,\ldots,m-2$. These operators $F_i$ correspond to the picture
\begin{center}
\includegraphics[width=9mm]{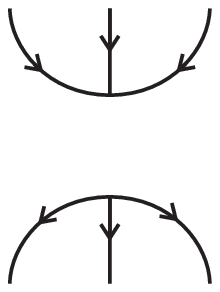}
\end{center}
in the $A_2$ web space.

\begin{Lemma} \label{Lemma:anitsymmetrizer_condition}
With $F_i$ defined as above, $F_i F_{i+1} F_i = \delta^2 F_i$ if and only if the $U_i$ satisfy the extra $SU(3)$ relation (\ref{eqn:SU(3)q_condition}).
\end{Lemma}
\emph{Proof}

The condition (\ref{eqn:SU(3)q_condition}) can be written as
\begin{eqnarray}
& & U_{i+2} U_{i+1} U_i U_{i+1} U_{i+2} U_{i+1} - U_i U_{i+1} - U_i U_{i+1} U_{i+2} U_{i+1} - U_{i+2} U_{i+1} U_i U_{i+1} \nonumber \\
& & \qquad = \;\; \delta (U_{i+1} U_{i+2} U_{i+1} - U_{i+1}) \label{eqn:quant_gp_antisymmetrizer2}.
\end{eqnarray}

We have
\begin{eqnarray*}
F_i F_{i+1} F_i & = & (U_{i+1} U_{i} U_{i+1} - U_{i+1})(U_{i+1} U_{i+2} U_{i+1} - U_{i+1})(U_{i+1} U_{i} U_{i+1} - U_{i+1}) \\
& = & (U_{i+1} U_{i} - \mathbf{1})(U_{i+1}^2 U_{i+2} U_{i+1}^2 - U_{i+1}^3)(U_{i} U_{i+1} - \mathbf{1}) \\
& = & (U_{i+1} U_{i} U_{i+1} - U_{i+1})(\delta^2 U_{i+2} - \delta \mathbf{1})(U_{i+1} U_{i} U_{i+1} - U_{i+1}) \\
& = & \delta (U_i U_{i+1} U_i - U_i)(\delta U_{i+2} - \mathbf{1})(U_i U_{i+1} U_i - U_i) \\
& = & \delta (\delta U_i U_{i+1} U_i U_{i+2} U_i U_{i+1} U_i - \delta U_i U_{i+1} U_i U_{i+2} U_i - \delta U_i U_{i+2} U_i U_{i+1} U_i \\
& & \quad + \delta U_i U_{i+2} U_i - U_i U_{i+1} U_i^2 U_{i+1} U_i + U_i U_{i+1} U_i^2 + U_i^2 U_{i+1} U_i - U_i^2).
\end{eqnarray*}
In the following we use relation H3 to transform each expression, and we indicate which terms have been replaced at each stage by enclosing them within square brackets $[ \; ]$. Since $U_i$, $U_{i+2}$ commute by H1, we have
\begin{eqnarray*}
& & \delta^2 (\delta U_i U_{i+1} U_{i+2} [U_i U_{i+1} U_i] - \delta U_i U_{i+1} U_{i+2} U_i - \delta U_{i+2} U_i U_{i+1} U_i + \delta U_i U_{i+2} \\
& & \quad - U_i [U_{i+1} U_i U_{i+1}] U_i + 2 U_i U_{i+1} U_i - U_i) \\
& = & \delta^2 (\delta U_i [U_{i+1} U_{i+2} U_{i+1}] U_i U_{i+1} - \delta U_i U_{i+1} U_{i+2} U_{i+1} + \delta U_i U_{i+1} U_{i+2} U_i - \delta U_i U_{i+1} U_{i+2} U_i \\
& & \quad - \delta U_{i+2} U_i U_{i+1} U_i + \delta U_i U_{i+2} - U_i^2 U_{i+1} U_i^2 - U_i U_{i+1} U_i + U_i^3 + 2 U_i U_{i+1} U_i - U_i) \\
& = & \delta^2 (\delta U_{i+2} [U_i U_{i+1} U_i] U_{i+2} U_{i+1} - \delta U_i U_{i+2} U_i U_{i+1} + \delta U_i [U_{i+1} U_i U_{i+1}] - \delta U_i U_{i+1} U_{i+2} U_i \\
& & \quad - \delta U_{i+2} [U_i U_{i+1} U_i] + \delta U_i U_{i+2} - (\delta^2 - 1)(U_i U_{i+1} U_i - U_i)) \\
& = & \delta^2 (\delta U_{i+2} U_{i+1} U_i U_{i+1} U_{i+2} U_{i+1} - \delta [U_{i+2} U_{i+1} U_{i+2}] U_{i+1} + \delta U_{i+2} U_i U_{i+2} U_{i+1} \\
& & \quad - \delta^2 U_{i+2} U_i U_{i+1} + \delta U_i^2 U_{i+1} U_i - \delta U_i^2 + \delta U_i U_{i+1} - \delta U_i U_{i+1} U_{i+2} U_{i+1} \\
& & \quad - \delta U_{i+2} U_{i+1} U_i U_{i+1} + \delta U_{i+2} U_{i+1} - \delta U_{i+2} U_i + \delta U_i U_{i+2} - (\delta^2 - 1)(U_i U_{i+1} U_i - U_i)) \\
& = & \delta^2 (\delta (U_{i+2} U_{i+1} U_i U_{i+1} U_{i+2} U_{i+1} + U_i U_{i+1} - U_i U_{i+1} U_{i+2} U_{i+1} - U_{i+2} U_{i+1} U_i U_{i+1}) \\
& & \quad - \delta U_{i+1} U_{i+2} U_{i+1}^2 + \delta U_{i+1}^2 - \delta U_{i+2} U_{i+1} + \delta U_{i+2} U_{i+1} + U_i U_{i+1} U_i - U_i) \\
& = & \delta^2 (\delta^2 (U_{i+1} U_{i+2} U_{i+1} - U_{i+1}) - \delta^2 (U_{i+1} U_{i+2} U_{i+1} - U_{i+1}) + U_i U_{i+1} U_i - U_i) \\
& = & \delta^2 F_i,
\end{eqnarray*}
where the penultimate equality follows from (\ref{eqn:quant_gp_antisymmetrizer2}).
\hfill
$\Box$

Note that if the condition (\ref{eqn:SU(3)q_condition}) is satisfied, $\mathrm{alg}(\mathbf{1}, F_i | i=1,\ldots,m-1)$ is not the Temperley-Lieb algebra, since although $F_i F_j = F_j F_i$ for $|i-j| > 2$, it is not the case for $|i-j|=2$, indeed $F_i F_{i+2} F_i = \delta F_i U_{i+3}$ so that $F_i$, $F_{i+2}$ do not commute.

We will now define a representation of the Hecke operators $U_k$ as elements of the path algebra for $\mathcal{ADE}$ graphs.
Let $\mathcal{G}$ be a finite $\mathcal{ADE}$ graph with Coxeter number $n < \infty$.
Let $M_0 = \mathbb{C}^{n_0}$ where $n_0$ is the number of 0-coloured vertices of $\mathcal{G}$, and let $M_0 \subset M_1 \subset M_2 \subset \cdots$ be finite dimensional von Neumann algebras, with the Bratteli diagram for the inclusion $M_j \subset M_{j+1}$ given by the graph $\mathcal{G}$, $j \geq 0$. Let $(\mu, \mu')$ be matrix units indexed by paths $\mu$, $\mu'$ on $\mathcal{G}$, and denote by $\mathfrak{E}^{\mathcal{G}}$, $\mathfrak{V}^{\mathcal{G}}$ the edges, vertices of $\mathcal{G}$ respectively. We define maps $s,r:\mathfrak{E}^{\mathcal{G}} \rightarrow \mathfrak{V}^{\mathcal{G}}$, where for an edge $\gamma \in \mathfrak{E}^{\mathcal{G}}$, $s(\gamma)$ denotes the source vertex of $\gamma$ and $r(\gamma)$ its range vertex. We define operators $U_k \in M_{k+1}$, for $k=1,2,\ldots$, by
\begin{equation} \label{Def:U_k}
U_{k} = \sum_{\sigma, \beta_i, \gamma_i} \mathcal{U}^{\beta_2,\gamma_2}_{\beta_1,\gamma_1}  (\sigma \cdot \beta_1 \cdot \gamma_1, \sigma \cdot \beta_2 \cdot \gamma_2),
\end{equation}
where the summation is over all paths $\sigma$ of length $k-1$ and edges $\beta_1, \beta_2, \gamma_1, \gamma_2$ of $\mathcal{G}$ such that  $r(\sigma) = s(\beta_1) = s(\beta_2)$, $s(\gamma_i) = r(\beta_i)$ for $i=1,2$, and $r(\gamma_1) = r(\gamma_2)$, and with $\mathcal{U}^{\beta_2,\gamma_2}_{\beta_1,\gamma_1}$ defined in (\ref{eqn:HeckeRep}). We will use the notation $W_{\rho_1, \rho_2, \rho_3}$ for $W(\triangle_{i_1, i_2, i_3}^{(\rho_1, \rho_2, \rho_3)})$, where $i_l = s(\rho_l)$, $l=1,2,3$.

\begin{Lemma} \label{Lemma:UiUi+1Ui-Ui=claw_crown}
With $U_k \in M_{k+1}$ given as in (\ref{Def:U_k}), the operator $F_k \in M_{k+2}$ defined in (\ref{Def:F_i}) is given by
\begin{equation} \label{F_i}
F_k = \sum_{\sigma, \beta_i, \gamma_i} \frac{1}{\phi_{r(\beta_3)}^2} W_{\gamma_1,\gamma_2,\gamma_3} \overline{W_{\beta_1,\beta_2,\beta_3}} \; (\sigma \cdot \beta_1 \cdot \beta_2 \cdot \beta_3, \sigma \cdot \gamma_1 \cdot \gamma_2 \cdot \gamma_3),
\end{equation}
where the summation is over all paths $\sigma$ of length $k-1$ and edges $\beta_i, \gamma_i$ of $\mathcal{G}$, $i=1,2,3$, such that $s(\beta_1) = s(\gamma_1) = r(\beta_3) = r(\gamma_3)$.
\end{Lemma}
\emph{Proof}

We have
\begin{eqnarray}
\lefteqn{U_k U_{k+1} U_k} \nonumber \\
& = & \sum_{\stackrel{\sigma_i, \beta_i,}{\scriptscriptstyle{\gamma_i, \mu_i}}} \mathcal{U}^{\beta_2,\gamma_2}_{\beta_1,\gamma_1} \mathcal{U}^{\beta_4,\gamma_4}_{\beta_3,\gamma_3} \mathcal{U}^{\beta_6,\gamma_6}_{\beta_5,\gamma_5} \; (\sigma_1 \cdot \beta_1 \cdot \gamma_1 \cdot \mu_1, \sigma_1 \cdot \beta_2 \cdot \gamma_2 \cdot \mu_1) \nonumber \\
& & \qquad \qquad \qquad \times (\sigma_2 \cdot \mu_2 \cdot \beta_3 \cdot \gamma_3, \sigma_2 \cdot \mu_2 \cdot \beta_4 \cdot \gamma_4) (\sigma_3 \cdot \beta_5 \cdot \gamma_5 \cdot \mu_3, \sigma_3 \cdot \beta_6 \cdot \gamma_6 \cdot \mu_3) \nonumber \\
& = & \sum_{\stackrel{\sigma_i, \beta_i,}{\scriptscriptstyle{\gamma_i, \mu_i}}} \mathcal{U}^{\beta_2,\gamma_2}_{\beta_1,\gamma_1} \mathcal{U}^{\beta_3,\mu_1}_{\beta_4,\mu_3} \mathcal{U}^{\beta_6,\gamma_6}_{\beta_2,\gamma_4} \; (\sigma_1 \cdot \beta_1 \cdot \gamma_1 \cdot \mu_1, \sigma_1 \cdot \beta_6 \cdot \gamma_6 \cdot \mu_3) \nonumber \\
& = & \sum_{\stackrel{\sigma_i, \beta_i, \gamma_i}{\scriptscriptstyle{\mu_i, \lambda_i}}} \frac{1}{\phi_{s(\beta_6)} \phi_{r(\gamma_6)} \phi_{s(\beta_4)} \phi_{r(\mu_1)} \phi_{s(\beta_1)} \phi_{r(\gamma_1)}} W_{\beta_6,\gamma_6,\lambda_1} \overline{W_{\beta_2,\beta_4,\lambda_1}} W_{\beta_4,\mu_3,\lambda_2} \overline{W_{\beta_3,\mu_1,\lambda_2}} W_{\beta_2,\beta_3,\lambda_3} \nonumber \\
& & \qquad \qquad \qquad \qquad \qquad \qquad \qquad \times \overline{W_{\beta_1,\gamma_1,\lambda_3}} \; (\sigma_1 \cdot \beta_1 \cdot \gamma_1 \cdot \mu_1, \sigma_1 \cdot \beta_6 \cdot \gamma_6 \cdot \mu_3) \nonumber \\
& = & \sum_{\stackrel{\sigma_i, \beta_i, \gamma_i}{\scriptscriptstyle{\mu_i, \lambda_i}}} \frac{1}{\phi_{s(\beta_6)} \phi_{r(\gamma_6)} \phi_{r(\mu_1)} \phi_{s(\beta_1)} \phi_{r(\gamma_1)}} W_{\beta_6,\gamma_6,\lambda_1} \overline{W_{\beta_1,\gamma_1,\lambda_3}} \; \Big(\delta_{\lambda_1,\mu_3} \delta_{\lambda_3,\mu_1} \phi_{s(\mu_3)} \phi_{r(\mu_3)} \phi_{s(\mu_1)} \nonumber \\
& & \qquad \qquad \qquad + \; \delta_{\lambda_1,\lambda_3} \delta_{\mu_1,\mu_3} \phi_{r(\lambda_1)} \phi_{s(\mu_3)} \phi_{r(\mu_3)} \Big) \; (\sigma_1 \cdot \beta_1 \cdot \gamma_1 \cdot \mu_1, \sigma_1 \cdot \beta_6 \cdot \gamma_6 \cdot \mu_3) \label{eqn:uses_TypeII_frame_again} \\
& = & \sum_{\stackrel{\sigma, \beta_i,}{\scriptscriptstyle{\gamma_i, \mu_i}}} \frac{1}{\phi_{r(\mu_1)}^2} W_{\beta_6,\gamma_6,\mu_3} \overline{W_{\beta_1,\gamma_1,\mu_1}} \; (\sigma_1 \cdot \beta_1 \cdot \gamma_1 \cdot \mu_1, \sigma_1 \cdot \beta_6 \cdot \gamma_6 \cdot \mu_3) \; + \; U_k, \nonumber
\end{eqnarray}
where we obtain (\ref{eqn:uses_TypeII_frame_again}) by Ocneanu's type II equation (\ref{eqn:typeII_frame}).
\hfill
$\Box$

Note that if $p$ is a minimal projection in $M_k$ corresponding to a vertex $(v,k)$ of the Bratteli diagram $\widehat{\mathcal{G}}$ of $\mathcal{G}$, then $F_{k+1} p$ is a projection in $M_{k+3}$ corresponding to the vertex $(v,k+3)$ of $\widehat{\mathcal{G}}$, since from (\ref{F_i}) we see that the last three edges in any pairs of paths in $F_{k+1}$ form a closed loop of length 3 and hence the pairs of paths in $F_{k+1} p \in M_{k+3}$ must have the same end vertex as $p \in M_k$.

\begin{Lemma}
The operators $U_k$ defined in (\ref{Def:U_k}) satisfy the $SU(3)$-Temperley-Lieb relations.
\end{Lemma}
\emph{Proof}

These operators satisfy the Hecke relations H1-H3 since the connection defined in (\ref{Connection_using_weights_W}) satisfies the Yang-Baxter equation. We are left to show that they satisfy (\ref{eqn:SU(3)q_condition}). By Lemma \ref{Lemma:anitsymmetrizer_condition}, we need only show that $F_k F_{k+1} F_k = [2]^2 F_k$. We have
\begin{eqnarray*}
\lefteqn{F_k F_{k+1} F_k} \\
& = & \sum_{\stackrel{\sigma_i, \beta_i,}{\scriptscriptstyle{\gamma_i, \mu_i}}} \frac{1}{\phi_{r(\beta_3)}^2 \phi_{r(\beta_6)}^2 \phi_{r(\beta_9)}^2} W_{\gamma_7,\gamma_8,\gamma_9} \overline{W_{\beta_7,\beta_8,\beta_9}} W_{\gamma_4,\gamma_5,\gamma_6} \overline{W_{\beta_4,\beta_5,\beta_6}} W_{\gamma_1,\gamma_2,\gamma_3} \overline{W_{\beta_1,\beta_2,\beta_3}} \\
& & \qquad (\sigma_1 \cdot \beta_1 \cdot \beta_2 \cdot \beta_3 \cdot \mu_1, \sigma_1 \cdot \gamma_1 \cdot \gamma_2 \cdot \gamma_3 \cdot \mu_1) (\sigma_2 \cdot \mu_2 \cdot \beta_4 \cdot \beta_5 \cdot \beta_6, \sigma_2 \cdot \mu_2 \cdot \gamma_4 \cdot \gamma_5 \cdot \gamma_6) \\
& & \qquad \qquad \times (\sigma_3 \cdot \beta_7 \cdot \beta_8 \cdot \beta_9 \cdot \mu_3, \sigma_3 \cdot \gamma_7 \cdot \gamma_8 \cdot \gamma_9 \cdot \mu_3) \\
& = & \sum_{\stackrel{\sigma_1, \beta_i,}{\scriptscriptstyle{\gamma_i, \mu_i}}} \frac{1}{\phi_{r(\beta_3)}^2 \phi_{r(\mu_1)}^2 \phi_{s(\mu_3)}^2} W_{\gamma_7,\gamma_8,\gamma_9} \overline{W_{\beta_7,\beta_8,\beta_9}} W_{\beta_8,\beta_9,\mu_3} \overline{W_{\beta_4,\beta_5,\mu_1}} W_{\beta_7,\beta_4,\beta_5} \overline{W_{\beta_1,\beta_2,\beta_3}} \\
& & \qquad \qquad (\sigma_1 \cdot \beta_1 \cdot \beta_2 \cdot \beta_3 \cdot \mu_1, \sigma_1 \cdot \gamma_7 \cdot \gamma_8 \cdot \gamma_9 \cdot \mu_3) \\
& = & [2]^2 \sum_{\stackrel{\sigma_1, \beta_i,}{\scriptscriptstyle{\gamma_i, \mu_i}}} \frac{\phi_{s(\mu_1)} \phi_{r(\mu_1)} \phi_{s(\mu_3)} \phi_{r(\mu_3)}}{\phi_{r(\beta_3)}^2 \phi_{r(\mu_1)}^2 \phi_{s(\mu_3)}^2} W_{\gamma_7,\gamma_8,\gamma_9} \overline{W_{\beta_1,\beta_2,\beta_3}} \delta_{\mu_1,\beta_7} \delta_{\mu_1,\mu_3} \\
& & \qquad \qquad (\sigma_1 \cdot \beta_1 \cdot \beta_2 \cdot \beta_3 \cdot \mu_1, \sigma_1 \cdot \gamma_7 \cdot \gamma_8 \cdot \gamma_9 \cdot \mu_3) \\
& = & [2]^2 \sum_{\stackrel{\sigma_1, \beta_i,}{\scriptscriptstyle{\gamma_i, \mu_1}}} \frac{1}{\phi_{r(\beta_3)}^2} W_{\gamma_7,\gamma_8,\gamma_9} \overline{W_{\beta_1,\beta_2,\beta_3}} \; (\sigma_1 \cdot \beta_1 \cdot \beta_2 \cdot \beta_3 \cdot \mu_1, \sigma_1 \cdot \gamma_7 \cdot \gamma_8 \cdot \gamma_9 \cdot \mu_1) \\
& = & [2]^2 F_k.
\end{eqnarray*}
\hfill
$\Box$

By \cite[Theorem 6.1]{effros:1981} there is a unique normalized faithful trace on $\bigcup_k M_k$, defined as in \cite{evans/kawahigashi:1994} by
\begin{equation} \label{eqn:trace}
\mathrm{tr}((\sigma_1,\sigma_2)) = \delta_{\sigma_1, \sigma_2} [3]^{-k} \phi_{r(\sigma_1)},
\end{equation}
for paths $\sigma_i$ of length $k$, $i=1,2$, $k=0,1,\ldots \;$. The conditional expectation of $M_k$ onto $M_{k-1}$ with respect to the trace is given by
$$E((\sigma_1 \cdot \sigma_1', \sigma_2 \cdot \sigma_2')) = \delta_{\sigma_1', \sigma_2'} [3]^{-1} \frac{\phi_{r(\sigma_1')}}{\phi_{r(\sigma_1)}} (\sigma_1, \sigma_2),$$
for paths $\sigma_i$ of length $k-1$, and $\sigma_i'$ of length 1, $i=1,2$, $k=1,2,\ldots \;$ (see e.g. \cite[Lemma 11.7]{evans/kawahigashi:1998}).

\begin{Lemma} \label{Lemma-Markov_trace_for_SU(3)}
For an $\mathcal{ADE}$ graph $\mathcal{G}$, let $M_0 = \mathbb{C}^{n_0}$ where $n_0$ is the number of 0-coloured vertices of $\mathcal{G}$. Let $M_0 \subset M_1 \subset M_2 \subset \cdots$ be a sequence of finite dimensional von Neumann algebras with normalized trace. Then for the operator $U_{k} \in M_{k+1}$ defined in (\ref{Def:U_k}), $\mathrm{tr}$ is a Markov trace in the sense that $\mathrm{tr}(x U_{k}) = [2] [3]^{-1} \mathrm{tr}(x)$ for any $x \in M_k$, $k \geq 1$.
\end{Lemma}
\emph{Proof}

Let $x \in M_k$ be the matrix unit $(\alpha_1 \cdot \alpha_1', \alpha_2 \cdot \alpha_2')$. Then
\begin{eqnarray*}
x U_{k} & = & \sum_{\sigma, \beta_i, \gamma_i, \mu} \mathcal{U}^{\beta_2,\gamma_2}_{\beta_1,\gamma_1}  (\alpha_1 \cdot \alpha_1' \cdot \mu, \alpha_2 \cdot \alpha_2' \cdot \mu) \cdot (\sigma \cdot \beta_1 \cdot \gamma_1, \sigma \cdot \beta_2 \cdot \gamma_2) \\
& = & \sum_{\sigma, \beta_i, \gamma_i, \mu} \mathcal{U}^{\beta_2,\gamma_2}_{\beta_1,\gamma_1}  \delta_{\alpha_2, \sigma} \delta_{\alpha_2', \beta_1} \delta_{\mu, \gamma_1} (\alpha_1 \cdot \alpha_1' \cdot \mu, \sigma \cdot \beta_2 \cdot \gamma_2) \\
& = & \sum_{\beta_2, \gamma_2, \mu} \mathcal{U}^{\beta_2,\gamma_2}_{\alpha_2',\mu}  (\alpha_1 \cdot \alpha_1' \cdot \mu, \alpha_2 \cdot \beta_2 \cdot \gamma_2),
\end{eqnarray*}
and
\begin{eqnarray*}
\mathrm{tr}(x U_{k}) & = & \sum_{\beta_2, \gamma_2, \mu} \mathcal{U}^{\beta_2,\gamma_2}_{\alpha_2',\mu}  \mathrm{tr}((\alpha_1 \cdot \alpha_1' \cdot \mu, \alpha_2 \cdot \beta_2 \cdot \gamma_2)) \\
& = & \sum_{\beta_2, \gamma_2, \mu} \mathcal{U}^{\beta_2,\gamma_2}_{\alpha_2',\mu}  \delta_{\alpha_1, \alpha_2} \delta_{\alpha_1', \beta_2} \delta_{\mu, \gamma_2} [3]^{-k+1} \phi_{r(\mu)} \;\; = \;\; \delta_{\alpha_1, \alpha_2} [3]^{-k+1} \sum_{\mu} \mathcal{U}^{\alpha_1',\mu}_{\alpha_2',\mu}  \phi_{r(\mu)} \\
& = & \delta_{\alpha_1, \alpha_2} [3]^{-k+1} \sum_{\mu} \frac{1}{\phi_{s(\alpha_1')} \phi_{r(\mu)}} W_{\lambda, \alpha_1', \mu} \overline{W_{\lambda, \alpha_2', \mu}} \phi_{r(\mu)} \\
& = & \delta_{\alpha_1, \alpha_2} [3]^{-k+1} \frac{1}{\phi_{s(\alpha_1')}} [2] \phi_{s(\alpha_1')} \phi_{r(\alpha_1')} \delta_{\alpha_1', \alpha_2'} \;\; = \;\; [2] [3]^{-1} \mathrm{tr}(x),
\end{eqnarray*}
where we have used Ocneanu's type I equation (\ref{eqn:typeI_frame}) in the penultimate equality. The result for any $x \in M_k$ follows by linearity of the trace.
\hfill
$\Box$

Then we have $\mathrm{tr}(U_k) = [2]/[3]$, and the conditional expectation of $U_k \in M_{k+1}$ onto $M_k$ is $E(U_k) = [2] \mathbf{1}_k /[3]$, for all $k \geq 1$. We will need the following result:

\begin{Lemma} \label{Lemma-tr(UUU-U)}
Let $F_i \in M_{i+2}$ be as above and $\mathrm{tr}$ a Markov trace on the $M_i$, $i=1,2,\ldots \;$, then $\mathrm{tr}(F_{k+1} x) = [2] [3]^{-2} \mathrm{tr}(x)$, for $x \in M_k$, $k \in \mathbb{N}$.
\end{Lemma}
\emph{Proof}

Now $\mathrm{tr}(U_{k+1} U_{k+2} U_{k+1}x) = \mathrm{tr}(U_{k+2} U_{k+1}x U_{k+1}) = [2][3]^{-1} \mathrm{tr}(U_{k+1}x U_{k+1})$, since $\mathrm{tr}$ is a Markov trace. Then $\mathrm{tr}(U_{k+1}x U_{k+1}) = \mathrm{tr}(U_{k+1}^2 x) = [2] \mathrm{tr}(U_{k+1}x) = [2]^2[3]^{-1} \mathrm{tr}(x)$. We also have $\mathrm{tr}(U_{k+1}x) = [2][3]^{-1} \mathrm{tr}(x)$, so that
\begin{eqnarray*}
\mathrm{tr}((U_{k+1} U_{k+2} U_{k+1} - U_{k+1})x) & = & \left( \frac{[2]^3}{[3]^2} - \frac{[2]}{[3]} \right) \mathrm{tr}(x) = \frac{[2]}{[3]^2} \mathrm{tr}(x).
\end{eqnarray*}
\hfill
$\Box$

\begin{Prop} \label{Prop-star}
With $U_k \in M_{k+1}$ as above and $x \in M_k$, $k=1,2,\ldots \;$, $x$ commutes with $U_k$ if and only if $x \in M_{k-1}$, i.e. $M_{k-1} = \{ U_k \}' \cap M_{k}$.
\end{Prop}
\emph{Proof}

Since $U_k \in M_{k-1}' \cap M_k$, it is clear that $x \in M_{k-1}$ commutes with $U_k$.

We now check the converse. Let $x = \sum_{\alpha_i, \alpha_i'} \lambda_{\alpha_1 \cdot \alpha_2, \alpha_1' \cdot \alpha_2'} (\alpha_1 \cdot \alpha_2, \alpha_1' \cdot \alpha_2') \in M_k$, where the summation is over all $|\alpha_i| = k-1$, $|\alpha_i'| = 1$, $i=1,2$. Assume that $x$ commutes with $U_k$. We have the inclusion of $x$ in $M_{k+1}$ given by $x = \sum_{\alpha_i, \alpha_i', \mu} \lambda_{\alpha_1 \cdot \alpha_2, \alpha_1' \cdot \alpha_2'} (\alpha_1 \cdot \alpha_2 \cdot \mu, \alpha_1' \cdot \alpha_2' \cdot \mu)$. Since $x$ commutes with $U_k$ we have $U_k^2 x = U_k x U_k$, and taking the conditional expectation onto $M_k$ we have
\begin{equation} \label{eqn:E(Ui(x))=E(Ui(x)U)}
[2] E(U_k x) = E(U_k x U_k).
\end{equation}
By the Markov property of the trace on the $M_k$, the left hand side gives  $[2] E(U_k x) = [2] E(U_k) x = [2]^2 x /[3]$, since $x \in M_k$. For the right hand side of (\ref{eqn:E(Ui(x))=E(Ui(x)U)}) we have
\begin{eqnarray*}
E(U_k x U_k) & = & E \Bigg( \sum_{\stackrel{\sigma, \beta_3, \beta_4,}{\scriptscriptstyle{\gamma_3, \gamma_4}}} \mathcal{U}^{\beta_4,\gamma_4}_{\beta_3,\gamma_3}  (\sigma \cdot \beta_3 \cdot \gamma_3, \sigma \cdot \beta_4 \cdot \gamma_4) \\
& & \qquad \qquad \times \sum_{\stackrel{\alpha_i, \alpha_i', \beta_2,}{\scriptscriptstyle{\gamma_2, \mu}}} \mathcal{U}^{\beta_2,\gamma_2}_{\alpha_2',\mu} \lambda_{\alpha_1 \cdot \alpha_2, \alpha_1' \cdot \alpha_2'} (\alpha_1 \cdot \alpha_2 \cdot \mu, \alpha_1' \cdot \beta_2 \cdot \gamma_2) \Bigg) \\
& = & E \left( \sum_{\stackrel{\sigma, \alpha_i, \alpha_i'}{\scriptscriptstyle{\beta_i, \gamma_i, \mu}}} \mathcal{U}^{\beta_4,\gamma_4}_{\beta_3,\gamma_3}  \mathcal{U}^{\beta_2,\gamma_2}_{\alpha_2',\mu}  \lambda_{\alpha_1 \cdot \alpha_2, \alpha_1' \cdot \alpha_2'} \delta_{\alpha_1, \sigma} \delta_{\beta_4, \alpha_2} \delta_{\gamma_4, \mu} (\sigma \cdot \beta_3 \cdot \gamma_3, \alpha_1' \cdot \beta_2 \cdot \gamma_2) \right) \\
& = & [3]^{-1} \sum_{\stackrel{\alpha_i, \alpha_i', \beta_i,}{\scriptscriptstyle{\gamma_i, \mu}}} \mathcal{U}^{\alpha_2,\mu}_{\beta_3,\gamma_3}  \mathcal{U}^{\beta_2,\gamma_2}_{\alpha_2',\mu}  \lambda_{\alpha_1 \cdot \alpha_2, \alpha_1' \cdot \alpha_2'} E( (\alpha_1 \cdot \beta_3 \cdot \gamma_3, \alpha_1' \cdot \beta_2 \cdot \gamma_2) ) \\
& = & [3]^{-1} \sum_{\stackrel{\alpha_1, \alpha_1'}{\scriptscriptstyle{\beta_3, \beta_2}}} \left( \sum_{\stackrel{\alpha_2, \alpha_2'}{\scriptscriptstyle{\gamma_i, \mu}}} \mathcal{U}^{\alpha_2,\mu}_{\beta_3,\gamma_3} \mathcal{U}^{\beta_2,\gamma_2}_{\alpha_2',\mu}  \lambda_{\alpha_1 \cdot \alpha_2, \alpha_1' \cdot \alpha_2'} \delta_{\gamma_2, \gamma_3} \frac{\phi_{r(\gamma_2)}}{\phi_{r(\beta_3)}} \right) (\alpha_1 \cdot \beta_3, \alpha_1' \cdot \beta_2) \\
& = & \alpha^{-1} \sum_{\stackrel{\alpha_1, \alpha_1'}{\scriptscriptstyle{\beta_1, \beta_2}}} b_{\alpha_1 \cdot \beta_1, \alpha_1' \cdot \beta_2} (\alpha_1 \cdot \beta_1, \alpha_1' \cdot \beta_2),
\end{eqnarray*}
where $b_{\alpha_1 \cdot \beta_1, \alpha_1' \cdot \beta_2} = \sum_{\stackrel{\alpha_2, \alpha_2'}{\scriptscriptstyle{\gamma, \mu}}} \mathcal{U}^{\alpha_2,\mu}_{\beta_1,\gamma}  \mathcal{U}^{\beta_2,\gamma}_{\alpha_2',\mu} \lambda_{\alpha_1 \cdot \alpha_2, \alpha_1' \cdot \alpha_2'} \frac{\phi_{r(\gamma)}}{\phi_{r(\beta_1)}}$. Then for any paths $\alpha_1$, $\alpha_1'$ and edges $\beta_1$, $\beta_2$ on $\mathcal{G}$ we have
\begin{eqnarray}
b_{\alpha_1 \cdot \beta_1, \alpha_1' \cdot \beta_2} & = & \sum_{\stackrel{\alpha_2, \alpha_2', \gamma}{\scriptscriptstyle{\mu, \zeta_i}}} \frac{1}{\phi_{s(\alpha_1)} \phi_{r(\gamma)}} W_{\beta_1 \gamma \zeta_1} \overline{W_{\alpha_2 \mu \zeta_1}} \frac{1}{\phi_{s(\alpha_2)} \phi_{r(\gamma)}} W_{\alpha_2' \mu \zeta_2} \overline{W_{\beta_2 \gamma \zeta_2}} \lambda_{\alpha_1 \cdot \alpha_2, \alpha_1' \cdot \alpha_2'} \frac{\phi_{r(\gamma)}}{\phi_{r(\beta_1)}} \nonumber \\
& = & \sum_{\alpha_2, \alpha_2'} \frac{1}{\phi_{s(\alpha_2)}^2 \phi_{r(\beta_1)}} \lambda_{\alpha_1 \cdot \alpha_2, \alpha_1' \cdot \alpha_2'} \left( \sum_{\gamma, \mu, \zeta_i} \frac{1}{\phi_{r(\gamma)}} W_{\beta_1 \gamma \zeta_1} \overline{W_{\alpha_2 \mu \zeta_1}} W_{\alpha_2' \mu \zeta_2} \overline{W_{\beta_2 \gamma \zeta_2}} \right) \nonumber \\
& = & \sum_{\alpha_2, \alpha_2'} \frac{1}{\phi_{s(\alpha_2)}^2 \phi_{r(\beta_1)}} \lambda_{\alpha_1 \cdot \alpha_2, \alpha_1' \cdot \alpha_2'} \Big(\phi_{r(\alpha_2)} \phi_{s(\alpha_2)} \phi_{r(\beta_1)} \delta_{\alpha_2, \alpha_2'} \delta_{\beta_1, \beta_2} \nonumber \\
& & \qquad \qquad \qquad \qquad \qquad \; \; + \phi_{s(\alpha_2)} \phi_{r(\beta_1)} \phi_{s(\alpha_2)} \delta_{\alpha_2, \beta_1} \delta_{\alpha_2', \beta_2} \Big) \label{uses_Ocneanu_typeII_eqn} \\
& = & \sum_{\alpha_2} \frac{\phi_{r(\alpha_2)}}{\phi_{s(\alpha_2)}} \lambda_{\alpha_1 \cdot \alpha_2, \alpha_1' \cdot \alpha_2'} \delta_{\beta_1, \beta_2} + \lambda_{\alpha_1 \cdot \beta_1, \alpha_1' \cdot \beta_2}, \label{eqn:here_beta1=beta2}
\end{eqnarray}
where equality (\ref{uses_Ocneanu_typeII_eqn}) follows by Ocneanu's type II equation (\ref{eqn:typeII_frame}). Since $\beta_1 = \beta_2$ in the first term in (\ref{eqn:here_beta1=beta2}), here $r(\alpha_1) = r(\alpha_1')$. We define
$$\lambda_{r(\alpha_1)} := \sum_{\beta'} \delta_{s(\beta'), r(\alpha_1)} \frac{\phi_{r(\beta')}}{\phi_{r(\alpha_1)}} \lambda_{\alpha_1 \cdot \beta', \alpha_1' \cdot \beta'},$$
which only depends on the range of the path $\alpha_1$ (which is equal to the range of $\alpha_1'$). Then we have for the right hand side of (\ref{eqn:E(Ui(x))=E(Ui(x)U)})
\begin{eqnarray*}
E(U_k x U_k) & = & [3]^{-1} \Bigg( \sum_{\stackrel{\beta_1, \beta_2,}{\scriptscriptstyle{\alpha_i, \alpha_1'}}} \frac{\phi_{r(\alpha_2)}}{\phi_{s(\alpha_2)}} \lambda_{\alpha_1 \cdot \alpha_2, \alpha_1' \cdot \alpha_2'} \delta_{\beta_1, \beta_2} (\alpha_1 \cdot \beta_1, \alpha_1' \cdot \beta_2) \\
& & \qquad \qquad + \sum_{\stackrel{\beta_1, \beta_2,}{\scriptscriptstyle{\alpha_1, \alpha_1'}}} \lambda_{\alpha_1 \cdot \beta_1, \alpha_1' \cdot \beta_2} (\alpha_1 \cdot \beta_1, \alpha_1' \cdot \beta_2) \Bigg) \\
& = & [3]^{-1} \left( \sum_{\beta, \alpha_1, \alpha_1'} \lambda_{s(\beta)} (\alpha_1 \cdot \beta, \alpha_1 \cdot \beta) + \sum_{\stackrel{\beta_1, \beta_2,}{\scriptscriptstyle{\alpha_1, \alpha_1'}}} \lambda_{\alpha_1 \cdot \beta_1, \alpha_1' \cdot \beta_2} (\alpha_1 \cdot \beta_1, \alpha_1' \cdot \beta_2) \right) \\
& = & [3]^{-1} (w + x),
\end{eqnarray*}
where $w = \sum_{\alpha_1, \alpha_1'} \lambda_{r(\alpha_1)} (\alpha_1, \alpha_1') \in M_{k-1}$.
Then (\ref{eqn:E(Ui(x))=E(Ui(x)U)}) gives $([2]^2 - 1) x = w$, so we see that $x \in M_{k-1}$.
\hfill
$\Box$

\emph{Remark.} The above proof was motivated by the following pictorial argument, which uses concepts which will be introduced in \cite{evans/pugh:2009iii}.

Let $\jmath$ be the inclusion of $M_{k-1}$ in $M_k$ and $\imath$ the inclusion of $M_k$ in $M_{k+1}$. For $x \in M_{k-1}$, we have the embedding $\imath \jmath (x)$ of $x$ into $M_{k+1}$, and $U_1 \in M_{k+1}$ given by the tangles:
\begin{center}
  \includegraphics[width=60mm]{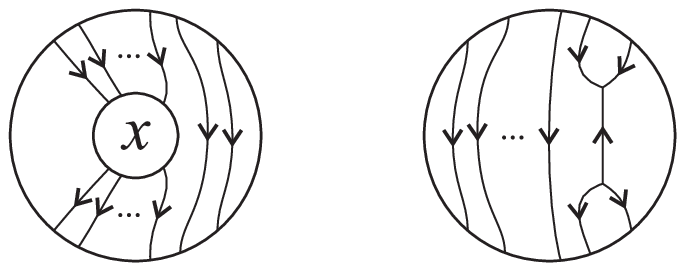}
\end{center}
Then inserting $x$ and $U_1$ into the discs of the multiplication tangle $M_{0,k+1}$, we have
\begin{center}
  \includegraphics[width=60mm]{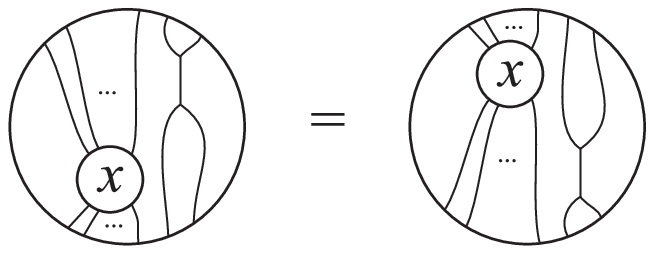}
\end{center}
and clearly $U_1 \imath \jmath (x) = \imath \jmath (x) U_1$.

\begin{figure}[bt]
\begin{center}
  \includegraphics[width=23mm]{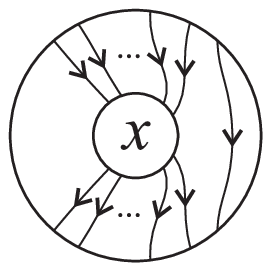}
 \caption{$\imath (x)$ for $x \in M_k$} \label{fig:i(x)}
\end{center}
\end{figure}

Conversely, if $x \in M_k$ we have $\imath(x) \in M_{k+1}$ as in Figure \ref{fig:i(x)}. Let $U_1 \imath (x) = \imath (x) U_1$, then we have the following equality of tangles:
\begin{center}
  \includegraphics[width=60mm]{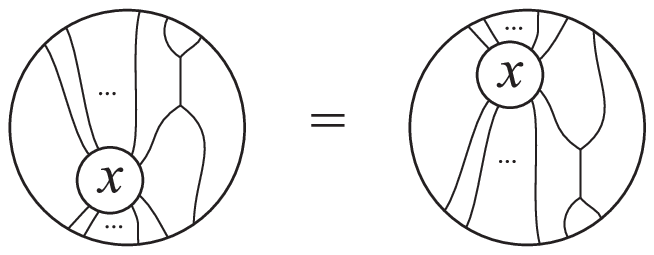}
\end{center}
Let $T$ be the tangle
\begin{center}
  \includegraphics[width=23mm]{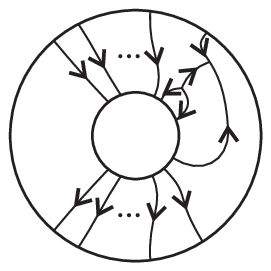}
\end{center}
We enclose both sides of $U_1 \imath (x) = \imath (x) U_1$ by the tangle $T$. Now $T (U_1 \imath(x)) = \delta^2 \imath(x)$, whilst $T (\imath(x) U_1)$ is
\begin{center}
  \includegraphics[width=90mm]{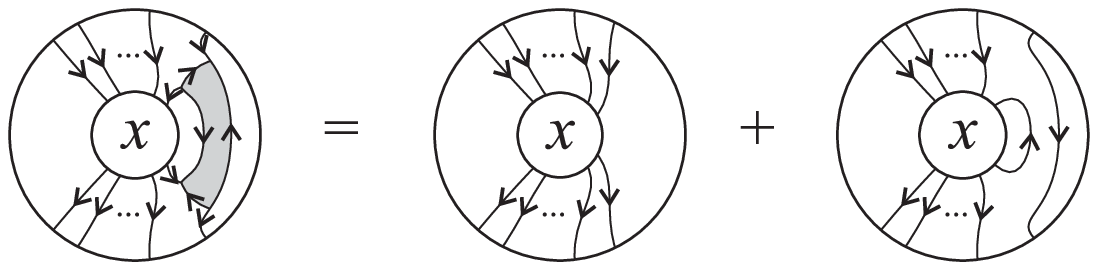}
\end{center}
i.e. $T (\imath(x) U_1) = x + \jmath(v)$, were $v = E_{M_{k-1}}(x) \in M_{k-1}$.
So $\delta^2 x = x + \jmath(v)$ which gives $x = (\delta^2 -1)^{-1} \jmath(v)$, i.e. $x \in M_{k-1}$.

We define the depth of the graph $\mathcal{G}$ to be $d_{\mathcal{G}} = \max_{v,v' \in \mathfrak{V}^{\mathcal{G}}} d_{v,v'}$, where $d_{v,v'}$ is the length of the shortest path between any two vertices $v, v'\in \mathfrak{V}^{\mathcal{G}}$.

\begin{Lemma} \label{Lemma:aU_kb}
Let $\mathcal{G}$ be an $SU(3)$ $\mathcal{ADE}$ graph $\mathcal{G}$ (except $\mathcal{D}^{(n)}$ for $n \not \equiv 0 \textrm{ mod } 3$, and $\mathcal{E}_4^{(12)}$). Then with $U_j \in M_{j+1}$ as above, any element of $M_{m+1}$ can be written as a linear combination of elements of the form $a U_m b$ and $c$ for $a,b,c \in M_m$, $m \geq d_{\mathcal{G}}+3$.
\end{Lemma}
\emph{Proof}

Let $a = (\lambda_1 \cdot \lambda_2, \zeta_1 \cdot \zeta_2), b = (\zeta_1 \cdot \zeta_2', \nu_1 \cdot \nu_2) \in M_m$ such that $\lambda_1$, $\zeta_1$, $\nu_1$ are paths of length $m-1$ on $\mathcal{G}$ starting from one of the 0-coloured vertices of $\mathcal{G}$, and $\lambda_2$, $\zeta_2$, $\zeta_2'$, $\nu_2$ are edges on $\mathcal{G}$. Then with $U_m$ as in (\ref{Def:U_k}), and embedding $a$, $b$ in $M_{m+1}$, we have
\begin{eqnarray}
a U_m b & = & \sum_{\sigma, \beta_i, \gamma_i, \mu, \mu'} \mathcal{U}^{\nu_2,\gamma_2}_{\nu_1,\gamma_1}  \delta_{\zeta_1,\sigma} \delta_{\zeta_2,\nu_1} \delta_{\mu,\gamma_1} \delta_{\nu_2,\zeta_2'} \delta_{\gamma_2,\mu'} \; (\lambda_1 \cdot \lambda_2 \cdot \mu, \nu_1 \cdot \nu_2 \cdot \mu') \nonumber \\
& = & \sum_{\mu, \mu'} \mathcal{U}^{\zeta_2',\mu'}_{\zeta_2,\mu} \; (\lambda_1 \cdot \lambda_2 \cdot \mu, \nu_1 \cdot \nu_2 \cdot \mu') \nonumber \\
& = & \sum_{\mu, \mu', \xi} \frac{1}{\phi_{s(\zeta_2)}\phi_{r(\mu)}} W(\triangle^{(\xi, \zeta_2, \mu)}) \overline{W(\triangle^{(\xi, \zeta_2', \nu)})} \; (\lambda_1 \cdot \lambda_2 \cdot \mu, \nu_1 \cdot \nu_2 \cdot \mu'). \label{eqn:commuting_with_U_k}
\end{eqnarray}

\begin{figure}[bt]
\begin{center}
  \includegraphics[width=40mm]{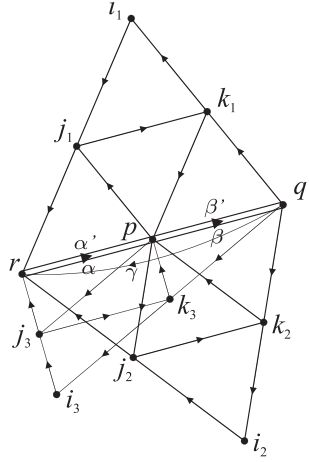}
 \caption{The $SU(3)$ graph $\mathcal{E}_1^{(12)}$} \label{fig:E1(12)}
\end{center}
\end{figure}

The proof for each graph is similar, so we illustrate the general method by considering the graph $\mathcal{E}_1^{(12)}$, illustrated in Figure \ref{fig:E1(12)}, which contains double edges. The proof for graphs without double edges is simpler. Let $m \geq d_{\mathcal{G}}+3$ be a fixed integer. We denote by $B$ the set of all linear combinations of elements of the form $a U_m b$ and $c$ for $a,b,c \in M_m$. We will write elements in $M_{m+1}$ in the form
\begin{equation} \label{x-Proof:Lemma5.1.7}
x = (\lambda_1 \cdot \lambda_2 \cdot \lambda_3, \nu_1 \cdot \nu_2 \cdot \nu_3)
\end{equation}
where $\lambda_1$, $\nu_1$ are paths of length $m-1$ on $\mathcal{G}$ with $s(\lambda_1) = s(\nu_1)$, and $\lambda_1$, $\lambda_2$, $\nu_1$, $\nu_2$ are edges of $\mathcal{G}$ with $r(\lambda_3) = r(\nu_3)$. Since the choice of the pair of paths $\lambda_1 \cdot \lambda_2$, $\nu_1 \cdot \nu_2$ in $a$, $b$ is arbitrary, the proof will depend on specific choices of $\zeta_2$, $\zeta_2'$ in (\ref{eqn:commuting_with_U_k}) in order to obtain the desired element. We label the vertices and some of the edges of $\mathcal{E}_1^{(12)}$ as in Figure \ref{fig:E1(12)}. For the other edges, let $\gamma_{v,v'}$ denote the edge on $\mathcal{E}_1^{(12)}$ from vertex $v$ to $v'$.

We first consider any element (\ref{x-Proof:Lemma5.1.7}) where $r(\lambda_2) = r(\nu_2)$. For any such pair $(\lambda_1 \cdot \lambda_2, \nu_1 \cdot \nu_2)$ with $r(\lambda_2) = i_l$, $l \in \{ 1,2,3 \}$, there is only one element $x$, which is given by the embedding of $x' = (\lambda_1 \cdot \lambda_2, \nu_1 \cdot \nu_2) \in M_m$ in $M_{m+1}$. If $r(\lambda_2) = i_l$, $l \in \{ 1,2,3 \}$, there are two possibilities for the edges $\lambda_3 = \nu_3$. If we choose $\zeta_2 = \zeta_2' = \gamma_{i_l,j_l}$ then (\ref{eqn:commuting_with_U_k}) gives $x^{(1)}_l = (\lambda_1 \cdot \lambda_2 \cdot \gamma_{j_l,k_l}, \nu_1 \cdot \nu_2 \cdot \gamma_{j_l,k_l})$, so that $x^{(1)}_l \in B$, $l=1,2,3$. Embedding $x'$ in $M_{m+1}$ we obtain $(\lambda_1 \cdot \lambda_2 \cdot \gamma_{j_l,r}, \nu_1 \cdot \nu_2 \cdot \gamma_{j_l,r}) = x' - x^{(1)}_l \in B$, for $l=1,2,3$. A similar method gives the result for the case when $r(\lambda_2) = r(\nu_2) = k_l$, $l=1,2,3$.

For any pair $(\lambda_1 \cdot \lambda_2, \nu_1 \cdot \nu_2)$ with $r(\lambda_2) = r(\nu_2) = p$, there are seven possibilities for $\lambda_3$, $\nu_3$. We denote these elements by $x^{(2)}_l$, $x_{(\xi,\xi')}$, for $l=1,2,3$, $\xi, \xi' \in \{ \beta, \beta' \}$, where $x^{(2)}_l = (\lambda_1 \cdot \lambda_2 \cdot \gamma_{p,j_l}, \nu_1 \cdot \nu_2 \cdot \gamma_{p,j_l})$, $x_{(\xi,\xi')} = (\lambda_1 \cdot \lambda_2 \cdot \xi, \nu_1 \cdot \nu_2 \cdot \xi')$. First, choosing $\zeta_2 = \zeta_2' = \alpha$, equation (\ref{eqn:commuting_with_U_k}) gives
\begin{eqnarray*}
y_0 & = & \frac{1}{\phi_r \phi_{j_1}} |W_{p,j_1,r(\alpha)}|^2 x_1^{(2)} + \frac{1}{\phi_r \phi_{j_2}} |W_{p,j_2,r(\alpha)}|^2 x_2^{(2)} + \frac{1}{\phi_r \phi_{j_3}} |W_{p,j_3,r(\alpha)}|^2 x_3^{(2)} \\
& & + \frac{1}{\phi_r \phi_q} |W_{p,q,r(\alpha \beta')}|^2 x_{(\beta',\beta')},
\end{eqnarray*}
where $y_0$ is an element in $B$. Using the solution $W^+$ for the cells of $\mathcal{E}_1^{(12)}$ given in \cite[Theorem 12.1]{evans/pugh:2009i}, we obtain
\begin{equation} \label{eqn:Lemma5.1.7-1}
y_1 = [2] r_1^+ \left(x_1^{(2)} + x_2^{(2)} + x_3^{(2)} \right) + [4] r_2^- \, x_{(\beta',\beta')},
\end{equation}
where $r_1^{\pm} = ([2][4] \pm \sqrt{[2][4]})$, $r_2^{\pm} = ([2]^2 \pm \sqrt{[2][4]})$ and $y_1 \in B$. Similarly, the choices $\zeta_2 = \zeta_2' = \alpha'$, $\; \zeta_2 = \alpha$, $\zeta_2' = \alpha'$ and $\zeta_2 = \alpha'$, $\zeta_2' = \alpha$ give
\begin{eqnarray}
y_2 & = & [2] r_1^- \left(x_1^{(2)} + x_2^{(2)} + x_3^{(2)} \right) + [4] r_2^+ \, x_{(\beta,\beta)}, \label{eqn:Lemma5.1.7-2} \\
y_3 & = & [2]\sqrt{r_1^+ r_1^-} \left(x_1^{(2)} + \overline{\omega} x_2^{(2)} + \omega x_3^{(2)} \right) + [4]\sqrt{r_2^+ r_2^-} \, x_{(\beta',\beta)}, \label{eqn:Lemma5.1.7-3} \\
y_4 & = & [2]\sqrt{r_1^+ r_1^-} \left(x_1^{(2)} + \omega x_2^{(2)} + \overline{\omega} x_3^{(2)} \right) + [4]\sqrt{r_2^+ r_2^-} \, x_{(\beta,\beta')}, \label{eqn:Lemma5.1.7-4}
\end{eqnarray}
where $\omega = e^{2 \pi i/3}$ and $y_j \in B$, $j=2,3,4$.
We can obtain three more equations by choosing $\zeta_2 = \zeta_2' = \gamma_{k_l,p}$ for $l=1,2,3$. Then (\ref{eqn:commuting_with_U_k}) gives
\begin{eqnarray}
y_5^{(l)} & = & x_1^{(2)} + x_2^{(2)} + x_3^{(2)} + \frac{[2]^2}{[3][4]^2} r_1^- \, x_{(\beta,\beta)} + \overline{\epsilon}_l \frac{[2]^2}{[3][4]^2}\sqrt{r_1^+ r_1^-} \, x_{(\beta,\beta')} \nonumber \\
& & + \epsilon_l \frac{[2]^2}{[3][4]^2}\sqrt{r_1^+ r_1^-} \, x_{(\beta',\beta)} + \frac{[2]^2}{[3][4]^2} r_1^+ \, x_{(\beta',\beta')}, \label{eqn:Lemma5.1.7-5}
\end{eqnarray}
where $\epsilon_l = \omega^{l-1}$ and $y_5^{(l)} \in B$, $l=1,2,3$.
Equations (\ref{eqn:Lemma5.1.7-1})-(\ref{eqn:Lemma5.1.7-5}) are linearly independent, and hence we can find $x_l^{(2)}$, $x_{(\xi,\xi')}$ in terms of $y_j$, $j=1,\ldots,4$, and $y_5^{(l)}$, for $l=1,2,3$, $\xi,\xi' \in \{ \beta, \beta' \}$; i.e. $x_l^{(2)}, x_{(\xi,\xi')} \in B$.

For any pair $(\lambda_1 \cdot \lambda_2, \nu_1 \cdot \nu_2)$ with $r(\lambda_2) = r(\nu_2) = q$, there are four possibilities for $\lambda_3$, $\nu_3$. We denote these elements by $x^{(3)}_l$, $x_r$, for $l=1,2,3$, where $x^{(3)}_l = (\lambda_1 \cdot \lambda_2 \cdot \gamma_{q,k_l}, \nu_1 \cdot \nu_2 \cdot \gamma_{q,k_l})$, $x_r = (\lambda_1 \cdot \lambda_2 \cdot \gamma, \nu_1 \cdot \nu_2 \cdot \gamma)$. Choosing $\zeta_2 = \zeta_2' = \beta$, equation (\ref{eqn:commuting_with_U_k}) gives
\begin{equation} \label{eqn:Lemma5.1.7-6}
y_6 = [2] r_1^- \left(x_1^{(3)} + x_2^{(3)} + x_3^{(3)} \right) + [4] r_2^+ \, x_r,
\end{equation}
where $y_6 \in B$. Similarly, the choices $\zeta_2 = \zeta_2' = \beta'$, $\; \zeta_2 = \beta$, $\zeta_2' = \beta'$ and $\zeta_2 = \beta'$, $\zeta_2' = \beta$ give
\begin{eqnarray}
y_7 & = & [2] r_1^+ \left(x_1^{(3)} + x_2^{(3)} + x_3^{(3)} \right) + [4] r_2^- \, x_r, \label{eqn:Lemma5.1.7-7} \\
y_8 & = & [2]\sqrt{r_1^+ r_1^-} \left(x_1^{(3)} + \overline{\omega} x_2^{(3)} + \omega x_3^{(3)} \right) + [4]\sqrt{r_2^+ r_2^-} \, x_r, \label{eqn:Lemma5.1.7-8} \\
y_9 & = & [2]\sqrt{r_1^+ r_1^-} \left(x_1^{(3)} + \omega x_2^{(3)} + \overline{\omega} x_3^{(3)} \right) + [4]\sqrt{r_2^+ r_2^-} \, x_r, \label{eqn:Lemma5.1.7-9}
\end{eqnarray}
where $y_j \in B$, $j=7,8,9$. Equations (\ref{eqn:Lemma5.1.7-6})-(\ref{eqn:Lemma5.1.7-9}) are linearly independent, and we find $x^{(3)}_l, x_r \in B$ for $l=1,2,3$.

For any pair $(\lambda_1 \cdot \lambda_2, \nu_1 \cdot \nu_2)$ with $r(\lambda_2) = r(\nu_2) = r$, there are four possibilities for $\lambda_3$, $\nu_3$, and we denote these elements by $x_{(\xi,\xi')} = (\lambda_1 \cdot \lambda_2 \cdot \xi, \nu_1 \cdot \nu_2 \cdot \xi')$, $\xi,\xi' \in \{ \alpha,\alpha' \}$. Choosing $\zeta_2 = \zeta_2' = \gamma$, equation (\ref{eqn:commuting_with_U_k}) gives
\begin{equation} \label{eqn:Lemma5.1.7-10}
y_{10} = r_2^- \, x_{(\alpha,\alpha)} + r_2^+ \, x_{(\alpha',\alpha')},
\end{equation}
where $y_{10} \in B$. We obtain three more equations by choosing $\zeta_2 = \zeta_2' = \gamma_{j_l,r}$, $l=1,2,3$:
\begin{equation} \label{eqn:Lemma5.1.7-11}
y_{11}^{(l)} = r_1^+ \, x_{(\alpha,\alpha)} + \overline{\epsilon}_l \sqrt{r_1^+ r_1^-} \, x_{(\alpha,\alpha')} + \epsilon_l \sqrt{r_1^+ r_1^-} \, x_{(\alpha',\alpha)} + r_1^- \, x_{(\alpha',\alpha')},
\end{equation}
where $y_{11}^{(l)} \in B$, $l=1,2,3$. So from (\ref{eqn:Lemma5.1.7-10}) and (\ref{eqn:Lemma5.1.7-11}) for $l=1,2,3$, we find that $x_{(\xi,\xi')} \in B$ for $\xi,\xi' \in \{ \alpha,\alpha' \}$.

We now consider any element $x$ in (\ref{x-Proof:Lemma5.1.7}) where $r(\lambda_2) \neq r(\nu_2)$. When $r(\lambda_2) = i_l$, $r(\nu_2) = p$, there is only one possibility for $\lambda_3$, $\nu_3$, which is $\lambda_3 = \gamma_{i_l,j_l}$, $\nu_3 = \gamma_{p,j_l}$, $l=1,2,3$, given by choosing $\zeta_2 = \gamma_{k_l,i_l}$, $\zeta_2' = \gamma_{k_l,p}$. Then $x = (\lambda_1 \cdot \lambda_2 \cdot \gamma_{i_l,j_l}, \nu_1 \cdot \nu_2 \cdot \gamma_{p,j_l}) \in B$. When $r(\lambda_2) = j_l$, $r(\nu_2) = j_{l+1}$, $l=1,2,3$, there is again only one possibility for $\lambda_3$, $\nu_3$. So $x \in B$. Similarly when $r(\lambda_2) = k_l$, $r(\nu_2) = k_{l+1}$, $l=1,2,3$.

Consider the pair $(\lambda_1 \cdot \lambda_2, \nu_1 \cdot \nu_2)$ where $r(\lambda_2) = j_l$, $l=1,2,3$, and $r(\nu_2) = q$. For each $l=1,2,3$, there are two possibilities for $\lambda_3$, $\nu_3$. We denote these by $x_l^{(4)} = (\lambda_1 \cdot \lambda_2 \cdot \gamma_{j_l,k_l}, \nu_1 \cdot \nu_2 \cdot \gamma_{q,k_l})$, $x_l^{(5)} = (\lambda_1 \cdot \lambda_2 \cdot \gamma_{j_l,r}, \nu_1 \cdot \nu_2 \cdot \gamma)$. Choosing $\zeta_2 = \gamma_{p,j_l}$, $\gamma_2' = \beta$, we obtain
\begin{equation} \label{eqn:Lemma5.1.7-12}
y_{12}^{(l)} = \sqrt{[3][4]} \, x_l^{(4)} - \sqrt{[2]} \sqrt{r_2^+} \, x_l^{(5)},
\end{equation}
where $y_{12}^{(l)} \in B$, $l=1,2,3$. Similarly, choosing $\zeta_2 = \gamma_{p,j_l}$, $\gamma_2' = \beta'$, we obtain
\begin{equation} \label{eqn:Lemma5.1.7-13}
y_{13}^{(l)} = \sqrt{[3][4]} \, x_l^{(4)} + \sqrt{[2]} \sqrt{r_2^-} \, x_l^{(5)},
\end{equation}
where $y_{13}^{(l)} \in B$, $l=1,2,3$. Then for each $l=1,2,3$, from (\ref{eqn:Lemma5.1.7-12}), (\ref{eqn:Lemma5.1.7-13}) we find that $x_l^{(4)}, x_l^{(5)} \in B$.

We now consider the pair $(\lambda_1 \cdot \lambda_2, \nu_1 \cdot \nu_2)$ where $r(\lambda_2) = k_l$, $l=1,2,3$, and $r(\nu_2) = r$. For each $l=1,2,3$, there are two possibilities for $\lambda_3$, $\nu_3$. We denote these by $x_{(\xi),l} = (\lambda_1 \cdot \lambda_2 \cdot \gamma_{k_l,p}, \nu_1 \cdot \nu_2 \cdot \xi)$, $\xi \in \{ \alpha, \alpha' \}$. Then for each $l = 1,2,3$, choosing $\zeta_2 = \gamma_{j_l,k_l}$, $\gamma_2' = \gamma_{j_l,r}$, we obtain
\begin{equation} \label{eqn:Lemma5.1.7-14}
y_{14}^{(l)} = \overline{\epsilon}_l \sqrt{r_1^+} x_{(\alpha),l} + \epsilon_l \sqrt{r_1^-} x_{(\alpha'),l},
\end{equation}
where $y_{14}^{(l)} \in B$, $l=1,2,3$. Similarly, choosing $\zeta_2 = \gamma_{q,k_l}$, $\gamma_2' = \gamma$, we obtain
\begin{equation} \label{eqn:Lemma5.1.7-15}
y_{15}^{(l)} = \sqrt{r_2^-} x_{(\alpha),l} - \sqrt{r_2^+} x_{(\alpha'),l},
\end{equation}
where $y_{15}^{(l)} \in B$, $l=1,2,3$. Then for each $l=1,2,3$, from (\ref{eqn:Lemma5.1.7-14}), (\ref{eqn:Lemma5.1.7-15}) we find that $x_{(\alpha),l}, x_{(\alpha'),l} \in B$.
All the other elements in $M_{m+1}$ are in $B$, since $y^{\ast} \in B$ if $y \in B$.
\hfill
$\Box$

The following lemma is an $SU(3)$ version of Skau's lemma. The proof is similar to the proof of Skau's lemma given in \cite[Theorem 4.4.3]{goodman/de_la_harpe/jones:1989}.

\begin{Lemma} \label{Lemma:SU(3)_Skau_lemma}
For an $\mathcal{ADE}$ graph $\mathcal{G}$, let $M_0 = \mathbb{C}^{n_0}$ where $n_0$ is the number of 0-coloured vertices of $\mathcal{G}$, and let $M_0 \subset M_1 \subset M_2 \subset \cdots$ be a tower of finite dimensional von Neumann algebras with Markov trace $\mathrm{tr}$ on the $M_i$, with the inclusions $M_j \subset M_{j+1}$ given by an $SU(3)$ $\mathcal{ADE}$ graph $\mathcal{G}$ (except $\mathcal{E}_4^{(12)}$), and operators $U_m \in M_{m+1}$, $m \geq 1$, which satisfy the relations H1-H3 for $\delta \leq 2$, and such that $U_m$ commutes with $M_{m-1}$. Let $M_{\infty}$ be the GNS-completion of $\bigcup_{j \geq 0} M_j$ with respect to the trace. Then $\{ U_1, U_2, \ldots \}' \cap M_{\infty} = M_0$.
\end{Lemma}
\emph{Proof}

The first inclusion $M_0 \subset \{ U_1, U_2, \ldots \}' \cap M_{\infty}$ is obvious, since $M_0$ commutes with $U_m$ for all $m \geq 1$.

We now show the opposite inclusion $M_0 \supset \{ U_1, U_2, \ldots \}' \cap M_{\infty}$. For each $k \geq 1$ let $F_k$ be the conditional expectation of $M_{\infty}$ onto $\{ U_k, U_{k+1}, \ldots \}' \cap M_{\infty}$ with respect to the trace. Note that $F_k F_l = F_{\min(k,l)}$. So we want to show $F_1 (M_{\infty}) \subset M_0$. We first show $F_2 (M_{\infty}) \subset M_m$ for some sufficiently large $m$. By \cite{goodman/de_la_harpe/jones:1989}, the diagram
$$\begin{array}{ccc}
  \{ U_{k+1}, U_{k+2}, \ldots \}' \cap M_{\infty} & \subset & M_{\infty} \\
   \cup & & \cup \\
   \{ U_{k+1}, U_{k+2}, \ldots \}' \cap \{ U_{k}, U_{k+1}, \ldots \}'' & \subset & \{ U_{k}, U_{k+1}, \ldots \}''
\end{array}$$
is a commuting square, for $k \geq 1$. Since $\{ U_{k+1}, U_{k+2}, \ldots \}'' \subset \{ U_{k}, U_{k+1}, \ldots \}''$ is isomorphic to $R_2 \subset R_1$, where $R_1 = \{ \mathbf{1}, U_1, U_2, \ldots \}''$, $R_2 = \{ \mathbf{1}, U_2, U_3, \ldots \}''$, we may write the commuting square as
$$\begin{array}{ccc}
   R_2' \cap M_{\infty} & \subset & M_{\infty} \\
   \cup & & \cup \\
   R_2' \cap R_1 & \subset & R_1.
\end{array}$$

Let $E$ denote the conditional expectation from $R_1$ onto $R_2' \cap R_1$ with respect to the trace. Since $F_{k+1}$ is the conditional expectation from $M_{\infty}$ onto $R_2' \cap M_{\infty}$ and $U_k \in R_1$, we have $F_{k+1}(U_k) = E(U_k)$.
Since by \cite[Cor. 3.4]{evans/kawahigashi:1994} the principal graph of $R_2 \subset R_1$ is the 01-part of $\mathcal{A}^{(n)}$, and there is only one vertex joined to the distinguished vertex $\ast$ of $\mathcal{A}^{(n)}$, the relative commutant $R_2' \cap R_1$ is trivial for $\alpha \leq 3$ (which corresponds to $\delta \leq 2$), and $E$ is just the trace.
Thus $F_{k+1}(U_k) \in \mathbb{C}$ for each $k \geq 1$. By Lemma \ref{Lemma:aU_kb}, for sufficiently large $m$, any element of $M_{m+1}$ can be written as $a U_m b$ for $a,b \in M_m$, and we have
\begin{eqnarray*}
F_2(a U_m b) & = & F_2 F_{m+1} (a U_m b) \;\; = \;\; F_2 (a F_{m+1} (U_m) b) \\
& = & F_{m+1} (U_m) F_2 (ab) \;\; = \;\; F_2 (\lambda ab) \;\; \in F_2 (M_m),
\end{eqnarray*}
where $\lambda \in \mathbb{C}$.
So $F_2 (M_{m+1}) \subset F_2 (M_m)$, for sufficiently large $m$, and by induction we have $F_2 (M_{\infty}) \subset F_2 (M_r)$, where $r$ is the smallest integer such that Lemma \ref{Lemma:aU_kb} holds. Then certainly $F_2 (M_{\infty}) \subset F_{r+1} (M_r)$, and by Proposition \ref{Prop-star}, with $k = r$, any element $x$ in $M_r$ commutes with $U_r$ if and only if $x \in M_{r-1}$, so $F_r F_{r+1} (M_r) \subset F_r (M_{r-1})$. Then by inductive use of Proposition \ref{Prop-star} we obtain $F_2 (M_{\infty}) \subset F_2 (M_1) = M_1$, and so $F_1 (M_{\infty}) = F_1 F_2 (M_{\infty}) \subset F_1 (M_1) = M_0$, by Proposition \ref{Prop-star}.
\hfill
$\Box$

We now construct the $SU(3)$-Goodman-de la Harpe-Jones subfactor for an $SU(3)$ $\mathcal{ADE}$ graph $\mathcal{G}$, following the idea of Goodman, de la Harpe and Jones for the $ADE$ Dynkin diagrams \cite{goodman/de_la_harpe/jones:1989}. Let $n$ be the Coxeter number for $\mathcal{G}$, $\ast_{\mathcal{G}}$ a distinguished vertex and let $n_0$ be the number of 0-coloured vertices of $\mathcal{G}$. Let $A_0$ be the von Neumann algebra $\mathbb{C}^{n_0}$, and form a sequence of finite dimensional von Neumann algebras $A_0 \subset A_1 \subset A_2 \subset \cdots$ such that the Bratteli diagram for the inclusion $A_{l-1} \subset A_l$ is given by (part of) the graph $\mathcal{G}$. There are operators $U_m \in A_{m+1}$ which satisfy the Hecke relations H1-H3. Let $\widetilde{C}$ be the GNS-completion of $\bigcup_{m \geq 0} A_m$ with respect to the trace, and $\widetilde{B}$ its von Neumann subalgebra generated by $\{ U_m \}_{m \geq 1}$. We have $\widetilde{B}' \cap \widetilde{C} = A_0$ by Lemma \ref{Lemma:SU(3)_Skau_lemma}. Then for $q$ the minimal projection in $A_0$ corresponding to the distinguished vertex $\ast_{\mathcal{G}}$ of $\mathcal{G}$, we have an $SU(3)$-Goodman-de la Harpe-Jones subfactor $B = q \widetilde{B} \subset q \widetilde{C} q = C$ for the graph $\mathcal{G}$. With $B_m = q \widetilde{B}_m$ and $C_m = q \widetilde{C}_m q$, the sequence $\{ B_m \subset C_m \}_m$ is a periodic sequence of commuting squares of period 3, in the sense of Wenzl in \cite{wenzl:1988}, that is, for large enough $m$ the Bratteli diagrams for the inclusions $B_m \subset B_{m+1}$, $C_m \subset C_{m+1}$ are the same as those for $B_{m+3} \subset B_{m+4}$, $C_{m+3} \subset C_{m+4}$, and the Bratteli diagrams for the inclusions $B_m \subset C_m$ and $B_{m+3} \subset C_{m+3}$ are the same. For such $m$ the graph of the Bratteli diagram for $B_{3m} \subset C_{3m}$ is the intertwining graph, given by the intertwining matrix $V$ computed in Proposition \ref{Prop:intertwinerV}, whose rows are indexed by the vertices of $\mathcal{G}$ and columns are indexed by the vertices of $\mathcal{A}^{(n)}$, such that $V \Delta_{\mathcal{A}} = \Delta_{\mathcal{G}} V$.
For sufficiently large $m$ we can make a basic construction $B_m \subset C_m \subset D_m$. Then with $D = \bigvee_m D_m$, $B \subset C \subset D$ is also a basic construction. The graph of the Bratteli diagram for $C_m \subset D_m$ is the reflection of the graph for $B_m \subset C_m$, which is the intertwining graph. Then we can extend the definition of $D_m$ to small $m$ so that the graph $C_m \subset D_m$ is still given by the reflection of the intertwining graph. We see that $D_0 = \bigoplus_{\mu \in \mathcal{A}^{(n)}} V V^{\ast}(\ast_{\mathcal{A}}, \mu) \mathbb{C}$, where $\ast_{\mathcal{A}}$ is the distinguished vertex $(0,0)$ of $\mathcal{A}^{(n)}$. The minimal projections in $D_0$ correspond to the vertices $\mu'$ of $\mathcal{A}^{(n)}$ such that
\begin{equation} \label{vertices_mu'}
V V^{\ast} (\ast, \mu') > 0,
\end{equation}
and the Bratteli diagram for the inclusion $D_{m-1} \subset D_m$ is given by (part of) the graph $\mathcal{A}^{(n)}$. Each algebra $B_m$ is generated by the $U_1, \ldots, U_{m-1}$ in $D_m$.

Now $\lambda_{(1,0)}(N) \subset N \cong P \subset Q$, where $P \subset Q$ is Wenzl's subfactor
with principal graph given by the 01-part $\mathcal{A}^{(n)}_{01}$ of $\mathcal{A}^{(n)}$ (see \cite[Cor. 3.4]{evans/kawahigashi:1994}). Then $(\lambda_{(1,0)}\overline{\lambda_{(1,0)}})^{d/2}(N) \cong P \subset Q_d$, where $P \subset Q \subset Q_1 \subset \cdots \;$ is the Jones tower.
For any 0-coloured vertex $\mu$ of $\mathcal{A}^{(n)}_{01}$ let $d_{\mu}$ be the minimum number of edges in any path from $(0,0)$ to $\mu$ on $\mathcal{A}^{(n)}_{01}$, and let $d = \mathrm{max} \{ d_{\mu}-2 | \; V V^{\ast}(\ast_{\mathcal{A}},\mu) > 0 \}$. Note that each $d_{\mu}$ is even since $\mu$ is a 0-coloured vertex.
Let $[\theta] = \bigoplus_{\mu \in \mathcal{A}^{(n)}} V V^{\ast}(\ast_{\mathcal{A}}, \mu) [\lambda_{\mu}]$. Now $[(\lambda_{(1,0)}\overline{\lambda_{(1,0)}})^{d/2}]$ decomposes into irreducibles as $\bigoplus_{\mu} n_{\mu} [\lambda_{\mu}]$, where $\mu$ are the 0-coloured vertices of $\mathcal{A}^{(n)}$ and $n_{\mu} \in \mathbb{N}$. Then $\theta(N) \subset N$ is a restricted version of $(\lambda_{(1,0)}\overline{\lambda_{(1,0)}})^{d/2}(N)$, so that $\theta(N) \subset N \cong q P \subset q (Q_d) q$ where $q \in P' \cap Q_d$ is a sum of minimal projections corresponding to the vertices $\mu'$ such that $[\theta] \supset [\lambda_{\mu'}]$. We will show that $q P \subset q (Q_d) q$ is isomorphic to a subfactor obtained by a basic construction.

Following the example in \cite[Lemma A.1]{bockenhauer/evans/kawahigashi:2000} for $E_7$ in the $SU(2)$ case, we now do the same construction for the graph $\mathcal{A}^{(n)}$, where $q$ is the projection corresponding to the distinguished vertex $\ast_{\mathcal{A}}$. We get a periodic sequence $\{ E_m \subset F_m \}_m$ of commuting squares of period 3. Then the resulting subfactor $E \subset F$, where $E = \bigvee_m E_m$, $F = \bigvee_m F_m$, is Wenzl's subfactor \cite{wenzl:1988}.

If we make basic constructions of $E_m \subset F_m$ for $d-1$ times then we get a periodic sequence $\{ E_m \subset G_m \}_m$ of commuting squares, and each $E_m$ is generated by the Hecke operators in $G_m$. Let $\widetilde{q}$ be a sum of the minimal projections corresponding to the vertices $\mu'$ in $G_0$ given by (\ref{vertices_mu'}). We set $\widetilde{E}_m = \widetilde{q} E_m$ and $\widetilde{G}_m = \widetilde{q} G_m \widetilde{q}$, and obtain a periodic sequence of commuting squares of period 3 such that the resulting subfactor is isomorphic to $qP \subset q(Q_d)q$. The Bratteli diagram for the sequence $\{ \widetilde{G}_m \}_m$ is the same as that for $\{ D_m \}_m$ since $D_0 = \widetilde{G}_0 = \mathbb{C}^r$ where the $r$ minimal projections correspond to the vertices $\mu'$ of (\ref{vertices_mu'}), where $r$ is the number of such vertices $\mu'$, and the rest of the Bratteli diagram is given by the 01-part of the graph $\mathcal{A}^{(n)}$. Each $\widetilde{E}_m$ is generated by the Hecke operators $U_1,\ldots,U_{m-1} \in \widetilde{G}_m$. Then the sequence of commuting squares $\{ B_m \subset D_m \}_m$ is isomorphic to the sequence of commuting squares $\{ \widetilde{E}_m \subset \widetilde{G}_m \}_m$, and so the subfactors $B \subset D$ and $qP \subset q(Q_d)q$ are also isomorphic. Since $B \subset D$ is a basic construction of $B \subset C$, then the subfactor $qP \subset q(Q_d)q$ is also the basic construction of some subfactor. Since $\theta(N) \subset N$ is isomorphic to $qP \subset q(Q_d)q$,
\begin{equation} \label{def:theta}
[\theta] = \bigoplus_{\mu \in \mathcal{A}^{(n)}} V V^{\ast}(\ast_{\mathcal{A}}, \mu) [\lambda_{\mu}]
\end{equation}
can be realised as the dual canonical endomorphism of some subfactor.

\subsection{Computing the intertwining graphs.}

Let $V(\mathcal{G})$ denote the free module over $\mathbb{Z}$ generated by the vertices of $\mathcal{G}$, identifying an element $a \in V(\mathcal{G})$ as $a = (a_v)$, $a_v \in \mathbb{Z}$, $v \in \mathfrak{V}^{\mathcal{G}}$. For graphs $\mathcal{G}_1$, $\mathcal{G}_2$, a map $V: V(\mathcal{G}_1) \longrightarrow V(\mathcal{G}_2)$ is \emph{positive} if $V_{ij} \geq 0$ for all $i \in \mathfrak{V}^{\mathcal{G}_2}$, $j \in \mathfrak{V}^{\mathcal{G}_1}$. Let $A(\mathcal{G})$ be the path algebra where the embeddings on the Bratteli diagram are given by the graph $\mathcal{G}$, and we will denote the finite dimensional algebra at the $k^{\mathrm{th}}$ level of the Bratteli diagram by $A(\mathcal{G})_k$.

The following lemma and proposition are the $SU(3)$ versions of Proposition 4.5 and Corollary 4.7 in \cite{evans/gould:1994i} (see also Lemma 11.26 and Proposition 11.27 in \cite{evans/kawahigashi:1998}).

\begin{Lemma} \label{Lemma-existence_of_intertwiner}
Suppose that $\mathcal{G}_1$, $\mathcal{G}_2$ are locally finite connected graphs with Coxeter number $n$, adjacency matrices $\Delta_{\mathcal{G}_1}$, $\Delta_{\mathcal{G}_2}$ respectively and distinguished vertices $\ast_1$, $\ast_2$ respectively. Let $(U_m)_{m \in \mathbb{N}}$, $(W_m)_{m \in \mathbb{N}}$ denote canonical families of operators in $A(\mathcal{G}_1)$ and $A(\mathcal{G}_2)$ respectively, which satisfy the $SU(3)$-Temperley-Lieb relations such that $U_m^2 = [2]_q U_m$, $W_m^2 = [2]_q W_m$ for all $m \in \mathbb{N}$, $q=e^{2 \pi i/n}$. Let $\pi : A(\mathcal{G}_1) \longrightarrow A(\mathcal{G}_2)$ be a unital embedding such that:
\begin{itemize}
\item[(a)] The diagram
    \begin{eqnarray*}
    A(\mathcal{G}_1)_m & \stackrel{\pi_m}{\longrightarrow} & A(\mathcal{G}_2)_m \\
    \iota_m \downarrow & & \downarrow \jmath_m \\
    A(\mathcal{G}_1)_{m+1} & \stackrel{\pi_{m+1}}{\longrightarrow} & A(\mathcal{G}_2)_{m+1} \\
    \end{eqnarray*}
commutes for all $m$, where $\pi_m = \pi|_{A(\mathcal{G}_1)_m}$, and $\iota_m$, $\jmath_m$ are standard inclusions.
\item[(b)] $\textrm{tr}_1 \cdot \pi_m = \textrm{tr}_2$, where $\textrm{tr}_i$ is a Markov trace on $A(\mathcal{G}_i)$, $i=1,2$.
\item[(c)] $\pi (U_m) = (W_m)$ for all $m \geq 1$ (so $\pi_{m+1} (U_m) = W_m$).
\end{itemize}
Then there exists a positive linear map $V: V(\mathcal{G}_1) \longrightarrow V(\mathcal{G}_2)$ such that:
\begin{itemize}
\item[(1)] $V \Delta_{\mathcal{G}_1} = \Delta_{\mathcal{G}_2} V$,
\item[(2)] $V$ has no zero rows or columns,
\item[(3)] $V \ast_1 = \ast_2$.
\end{itemize}
\end{Lemma}
\emph{Proof}

We denote by $\widehat{\mathcal{G}}$ the Bratteli diagram of $\mathcal{G}$. The vertex $(i,m)$ of $\widehat{\mathcal{G}}$ will be the vertex $i \in \mathfrak{V}^{\mathcal{G}}$ at level $m$ of the Bratteli diagram.
Let $p_i^m$ denote a minimal projection in $A(\mathcal{G}_1)_m$ corresponding to the vertex $(i,m)$ of the Bratteli diagram $\widehat{\mathcal{G}}_1$ of $\mathcal{G}_1$. Then $\pi_m(p_i^m)$ is a projection in $A(\mathcal{G}_2)_m$, and so there are families of equivalent minimal projections $\{ q_{j, k(j)}^m | k(j) = 1, \ldots, b_{ji}^m \}$ in $A(\mathcal{G}_2)_m$ corresponding to vertices $(j,m)$ in $\widehat{\mathcal{G}}_2$, such that
\begin{equation}\label{pi_n(p)}
\pi_m(p_i^m) = \sum_j \sum_{k(j)=1}^{b_{ji}^m} q_{j, k(j)}^m.
\end{equation}
The numbers $\{ b_{ji}^m \}_j$ are non-negative, are independent of the choice of $p_i^m$ and are not all zero, since $\pi_m$ is injective. Let $F_m^{(1)} = [2]^{-1} [3]^{-1}(U_m U_{m+1} U_m - U_m)$ in $A(\mathcal{G}_1)$, and $F_m^{(2)} = [2]^{-1} [3]^{-1}(W_m W_{m+1} W_m - W_m)$ in $A(\mathcal{G}_2)$. Now multiplying (\ref{pi_n(p)}) on the left by $F_{m+1}^{(2)}$, we have
$$F_{m+1}^{(2)} \pi_m(p_i^m) = \sum_j \sum_{k(j)=1}^{b_{ji}^m} F_{m+1}^{(2)} \, q_{j, k(j)}^m,$$
but by (a) and (c), $F_{m+1}^{(2)} \pi_m(p_i^m) = \pi_{m+3}(F_{m+1}^{(1)}) \pi_m(p_i^m) = \pi_{m+3} (F_{m+1}^{(1)} p_i^m)$, so we have
\begin{equation}\label{pi_n+2(p)}
\pi_{m+2}(F_{m+1}^{(1)} p_i^m) = \sum_j \sum_{k(j)=1}^{b_{ji}^m} F_{m+1}^{(2)} \, q_{j, k(j)}^m.
\end{equation}
Since $\textrm{tr}_1$ and $\textrm{tr}_2$ are Markov traces, by Lemma \ref{Lemma-tr(UUU-U)} we have $\textrm{tr}_1(F_{m+1}^{(1)} \, p_i^m) = [3]^{-3} \textrm{tr}_1(p_i^m)$, and $\textrm{tr}_2(F_{m+1}^{(2)} \, q_{j, k(j)}^m) = [3]^{-3} \textrm{tr}_2(q_{j, k(j)}^m)$. Since $p_i^m$, $q_{j, k(j)}^m$ are minimal projections, they have trace $[3]^{-k} \phi_i$, $[3]^{-k} \phi_j$ respectively. Then $F_{m+1}^{(1)} \, p_i^m$ has trace $[3]^{-k-3} \phi_i$, which shows that $F_{m+1}^{(1)} \, p_i^m$ is a minimal projection in $A(\mathcal{G}_1)_{m+3}$ corresponding to vertex $(i,m+3)$ of $\widehat{\mathcal{G}}_1$, and similarly $F_{m+1}^{(2)} \, q_{j, k(j)}^m$ is a minimal projection in $A(\mathcal{G}_2)_{m+3}$ corresponding to vertex $(j,m+3)$ of $\widehat{\mathcal{G}}_2$. It follows from (\ref{pi_n(p)}) and (\ref{pi_n+2(p)}) that the coefficients occurring in the decomposition of a minimal projection as in (\ref{pi_n(p)}) corresponding to vertex $(i,m)$ of $\widehat{\mathcal{G}}_1$, $m \geq 1$, is independent of the level $m$, i.e. $b_{ji}^m = b_{ji}^l =: b_{ji}$ for all $m,l \geq 0$.

Now put $V=(b_{ji})_{i \in \mathfrak{V}^{\mathcal{G}_1}, j \in \mathfrak{V}^{\mathcal{G}_2}}$, then since $A(\mathcal{G}_1)_0 \cong \mathbb{C} \cong A(\mathcal{G}_2)_0$, and $\pi_0 : A(\mathcal{G}_1)_0 \longrightarrow A(\mathcal{G}_2)_0$ we see that $V \ast_1 = \ast_2$. Note that since $\pi$ is unital, the rows of $V$ are non-zero. We need to show $V \Delta_{\mathcal{G}_1} = \Delta_{\mathcal{G}_2} V$.

Let $\Delta_{\mathcal{G}_k} (m)$, $k=1,2$, be the finite submatrix of $\Delta_{\mathcal{G}_k}$, whose rows and columns are labelled by the vertices $v \in \mathcal{G}_k^{(0)}$ with $d(v) \leq m+1$, where $d(v)$ is the distance of vertex $v$ from $\ast_k$, ie. the length of the shortest path on $\mathcal{G}_k$ from $\ast_k$ to $v$. Similarly let $V(m)$ denote the finite submatrix of $V$ whose rows are labelled by $j \in \mathfrak{V}^{\mathcal{G}_2}$ with $d(j) \leq m+1$, and whose columns are labelled by $i \in \mathfrak{V}^{\mathcal{G}_1}$ with $d(i) \leq m+1$. It follows from (a) that for each $m$ we have
\begin{equation}\label{K_0(intertwiner)}
K_0(\jmath_m) K_0(\pi_m) = K_0(\pi_{m+1}) K_0(\iota_m).
\end{equation}
Let $M_1$, $M_2$, be two multi-matrix algebras, with the embedding $\varphi$ of $M_1$ in $M_2$ given by a matrix $\Lambda$, with $p_1$ columns corresponding to the minimal central projections in $M_1$ and $p_2$ rows corresponding to the minimal central projections in $M_2$. Then $K_0(M_i)= \mathbb{Z}^{p_i}$, $i=1,2$, and $K_0(\varphi):\mathbb{Z}^{p_1} \rightarrow \mathbb{Z}^{p_2}$ is given by multiplication by the matrix $\Lambda$.
For $m$ of colour $\overline{j}$, we see that $K_0(\iota_m)$ is the submatrix of $\Delta_{\mathcal{G}_1}(m)$ mapping vertices of colour $\overline{j}$ to vertices of colour $\overline{j+1}$, and $K_0(\jmath_m)$ is the submatrix of $\Delta_{\mathcal{G}_2}(m)$ mapping vertices of colour $\overline{j}$ to vertices of colour $\overline{j+1}$. Similarly, $K_0(\pi_m)$ is the submatrix of $V(m)$ mapping vertices of $\mathcal{G}_1$ of colour $\overline{j}$ to vertices of $\mathcal{G}_2$ of colour $\overline{j}$. Then (\ref{K_0(intertwiner)}) implies $\Delta_{\mathcal{G}_2}(m) V(m-1) = V(m) \Delta_{\mathcal{G}_1}(m)$ holds for all $m$. Hence $V \Delta_{\mathcal{G}_1} = \Delta_{\mathcal{G}_2} V$.
\hfill
$\Box$

We define polynomials $S_{\nu}(x,y)$, for $\nu$ the vertices of $\mathcal{A}^{(n)}$, by $S_{(0,0)}(x,y) = 1$, and $x S_{\nu}(x,y) = \sum_{\mu} \Delta_{\mathcal{A}}(\nu, \mu) S_{\mu}(x,y)$, $\; y S_{\nu}(x,y) = \sum_{\mu} \Delta_{\mathcal{A}}^T(\nu, \mu) S_{\mu}(x,y)$. For concrete values of the first few $S_{\mu}(x,y)$ see \cite[p. 610]{evans/kawahigashi:1998}.

\begin{Prop} \label{Prop:intertwinerV}
Let $\mathcal{G}$ be a finite $SU(3)$-$\mathcal{ADE}$ graph with distinguished vertex $\ast_{\mathcal{G}}$ and Coxeter number $n < \infty$. Let $\{U_m \}_{m \geq 0}$, $\{W_m \}_{m \geq 0}$ be the canonical family of operators satisfying the Hecke relations in $A(\mathcal{A}^{(n)})$, $A(\mathcal{G})$ respectively. We can identify $A(\mathcal{A}^{(n)})$ with the algebra generated by $\{ \mathbf{1}, W_1, W_2, \ldots \}$. If we define $\pi: A(\mathcal{A}^{(n)}) \longrightarrow A(\mathcal{G})$ by $\pi(\mathbf{1}) = \mathbf{1}$, $\pi(U_m) = W_m$, then $\pi$ is a unital embedding, and there exists a positive linear map $V: V(\mathcal{A}^{(n)}) \longrightarrow V(\mathcal{G})$ such that:
\begin{itemize}
\item[(a)] $V \Delta_{\mathcal{A}} = \Delta_{\mathcal{G}} V$,
\item[(b)] $V$ has no zero rows or columns,
\item[(c)] $V \ast_{\mathcal{A}} = \ast_{\mathcal{G}}$, where $\ast_{\mathcal{A}}=(0,0)$ is the distinguished vertex of $\mathcal{A}^{(n)}$.
\end{itemize}
Let $V_{(0,0)}$ be the vector corresponding to the distinguished vertex $\ast_{\mathcal{G}}$, and for the other vertices define $V_{(\lambda_1, \lambda_2)} \in V(\mathcal{G})$ by $V_{(\lambda_1, \lambda_2)} = S_{(\lambda_1, \lambda_2)} \left( \Delta_{\mathcal{G}}^T, \Delta_{\mathcal{G}} \right) V_{(0,0)}$, for all vertices  $(\lambda_1, \lambda_2)$ of $\mathcal{A}^{(n)}$. Then $V = (V_{(0,0)}, V_{(1,0)}, V_{(0,1)}, V_{(2,0)}, \ldots , V_{(0,n-3)})$.
\end{Prop}
\emph{Proof}

Now $\pi: A(\mathcal{A}^{(n)}) \longrightarrow A(\mathcal{G})$ defined by $\pi(\mathbf{1}) = \mathbf{1}$, $\pi(U_m) = W_m$ is a unital embedding which satisfies the condition of Lemma \ref{Lemma-existence_of_intertwiner} with $\ast_1 = (0,0)$ and $\ast_2 = \ast_{\mathcal{G}}$. Hence when $m$ is finite there exists $V = (V_{(\lambda_1, \lambda_2)})$, for $(\lambda_1, \lambda_2)$ the vertices of $\mathcal{A}^{(n)}$, with the required properties. Now $V \Delta_{\mathcal{A}} = \left( V_{(\lambda_1 -1, \lambda_2)} + V_{(\lambda_1 +1, \lambda_2 -1)} + V_{(\lambda_1, \lambda_2 +1)} \right)_{(\lambda_1, \lambda_2)}$, where $V_{(\lambda_1, \lambda_2)}$ is understood to be zero if $(\lambda_1, \lambda_2)$ is off the graph $\mathcal{A}^{(n)}$. Thus $V \Delta_{\mathcal{A}} = \Delta_{\mathcal{G}} V$ implies that $\Delta_{\mathcal{G}} V_{(\lambda_1, \lambda_2)} = V_{(\lambda_1 -1, \lambda_2)} + V_{(\lambda_1 +1, \lambda_2 -1)} + V_{(\lambda_1, \lambda_2 +1)}$. Then $V_{(\lambda_1, \lambda_2)} = S_{(\lambda_1, \lambda_2)} \left( \Delta_{\mathcal{G}}^T, \Delta_{\mathcal{G}} \right) V_{(0,0)}$, since
\begin{eqnarray*}
\Delta_{\mathcal{G}} V_{(\lambda_1, \lambda_2)} & = & \Delta_{\mathcal{G}} S_{(\lambda_1, \lambda_2)} \left( \Delta_{\mathcal{G}}^T, \Delta_{\mathcal{G}} \right) V_{(0,0)} \\
& = & \sum_{(\mu_1, \mu_2)} \Delta_{\mathcal{A}}^T \left( (\lambda_1, \lambda_2), (\mu_1, \mu_2) \right) S_{(\mu_1, \mu_2)} \left( \Delta_{\mathcal{G}}^T, \Delta_{\mathcal{G}} \right) V_{(0,0)} \\
& = & V_{(\lambda_1 -1, \lambda_2)} + V_{(\lambda_1 +1, \lambda_2 -1)} + V_{(\lambda_1, \lambda_2 +1)},
\end{eqnarray*}
and $V_{(0,0)} = S_{(0,0)} \left( \Delta_{\mathcal{G}}^T, \Delta_{\mathcal{G}} \right) V_{(0,0)}$.
\hfill
$\Box$

\begin{figure}[htb]
\begin{center}
  \includegraphics[width=65mm]{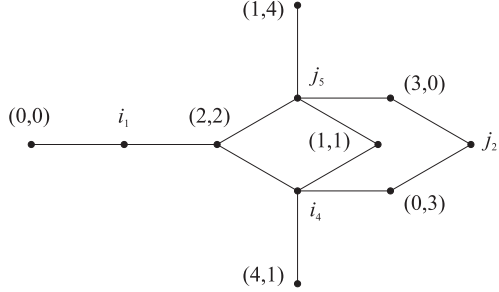}
 \caption{Principal graph for the $SU(3)$-Goodman-de la Harpe-Jones subfactor for $\mathcal{E}^{(8)}$}\label{fig:fig-GHJ-8}
\end{center}
\end{figure}

For any $\mathcal{ADE}$ graph $\mathcal{G}$ the matrix $V$ is the adjacency matrix of a (possibly disconnected) graph. By \cite[Theorem 4.2]{bockenhauer/evans:1999ii} the connected component of $\ast_{\mathcal{A}}$ of this graph gives the principal graph of the $SU(3)$-Goodman-de la Harpe-Jones subfactor. For the graph $\mathcal{E}^{(8)}$ with vertex $i_1$ chosen as the distinguished vertex this is the graph illustrated in Figure \ref{fig:fig-GHJ-8}, which was shown to be the principal graph for this subfactor in \cite{xu:1997}.

\section{Modular invariants associated to the dual canonical endomorphisms.} \label{Sect:Compute_N-M_graphs}

Let $N \subset M$ be the $SU(3)$-GHJ subfactor for the finite $\mathcal{ADE}$ graph $\mathcal{G}$, where the distinguished vertex $\ast_{\mathcal{G}}$ is the vertex with lowest Perron-Frobenius weight. Then the dual canonical endomorphism $\theta$ for $N \subset M$ is given by (\ref{def:theta}) where $V$ is now determined in Proposition \ref{Prop:intertwinerV}. We list these $\theta$'s below for the $\mathcal{ADE}$ graphs, where we use the same notation for the $\mathcal{ADE}$ graphs as in \cite{evans/pugh:2009i}. We must point out that as we have been unable to explicitly construct the Ocneanu cells $W$ for $\mathcal{E}_4^{(12)}$, the existence of the $SU(3)$-Goodman-de la Harpe-Jones subfactor which realizes the candidate for the dual canonical endomorphism for $\mathcal{E}_4^{(12)}$ is not shown here.
\begin{eqnarray}
\mathcal{A}^{(n)}: & & [\theta] = [\lambda_{(0,0)}], \label{theta-A(n)} \\
\mathcal{D}^{(n)}: & & [\theta] = [\lambda_{(0,0)}] \oplus [\lambda_{A(0,0)}] \oplus [\lambda_{A^2(0,0)}],  \label{theta-D(n)} \\
\mathcal{A}^{(n)\ast}: & & [\theta] = \bigoplus_{\mu \in \mathcal{A}^{(n)}} [\lambda_{\mu}], \label{theta-A(n)star} \\
\mathcal{D}^{(2k)\ast}: & & [\theta] = \bigoplus_{\stackrel{\mu \in \mathcal{A}^{(2k)}:}{\scriptscriptstyle{\tau(\mu)=0}}} [\lambda_{\mu}], \\
\mathcal{D}^{(2k+1)\ast}: & & [\theta] = \bigoplus_{\stackrel{\mu = (2 \mu_1, 2 \mu_2) \in \mathcal{A}^{(2k+1)}:}{\scriptscriptstyle{\tau(\mu)=0}}} [\lambda_{\mu}], \\
\mathcal{E}^{(8)}: & & [\theta] = [\lambda_{(0,0)}] \oplus [\lambda_{(2,2)}], \\
\mathcal{E}^{(8)\ast}: & & [\theta] = [\lambda_{(0,0)}] \oplus [\lambda_{(2,1)}] \oplus [\lambda_{(1,2)}] \oplus [\lambda_{(2,2)}] \oplus [\lambda_{(5,0)}] \oplus [\lambda_{(0,5)}], \\
\mathcal{E}_1^{(12)}: & & [\theta] = [\lambda_{(0,0)}] \oplus [\lambda_{(4,1)}] \oplus [\lambda_{(1,4)}] \oplus [\lambda_{(4,4)}] \oplus [\lambda_{(9,0)}] \oplus [\lambda_{(0,9)}], \\
\mathcal{E}_2^{(12)}: & & [\theta] = [\lambda_{(0,0)}] \oplus 2 [\lambda_{(2,2)}] \oplus [\lambda_{(4,1)}] \oplus [\lambda_{(1,4)}] \oplus 2 [\lambda_{(5,2)}] \oplus 2 [\lambda_{(2,5)}] \nonumber \\
& & \qquad \quad \oplus [\lambda_{(4,4)}] \oplus [\lambda_{(9,0)}] \oplus [\lambda_{(0,9)}], \label{theta-E2(12)} \\
\mathcal{E}_4^{(12)}: & & [\theta] = [\lambda_{(0,0)}] \oplus [\lambda_{(2,2)}] \oplus [\lambda_{(4,1)}] \oplus [\lambda_{(1,4)}] \oplus [\lambda_{(5,2)}] \oplus [\lambda_{(2,5)}] \oplus [\lambda_{(4,4)}] \nonumber \\
& & \qquad \quad \oplus [\lambda_{(9,0)}] \oplus [\lambda_{(0,9)}], \label{theta-E4(12)} \\
\mathcal{E}_5^{(12)}: & & [\theta] = [\lambda_{(0,0)}] \oplus [\lambda_{(3,3)}] \oplus [\lambda_{(9,0)}] \oplus [\lambda_{(0,9)}], \\
\mathcal{E}^{(24)}: & & [\theta] = [\lambda_{(0,0)}] \oplus [\lambda_{(4,4)}] \oplus [\lambda_{(10,1)}] \oplus [\lambda_{(1,10)}] \oplus [\lambda_{(6,6)}] \oplus [\lambda_{(9,6)}] \oplus [\lambda_{(6,9)}] \nonumber \\
& & \qquad \quad \oplus [\lambda_{(13,4)}] \oplus [\lambda_{(4,13)}] \oplus [\lambda_{(10,10)}] \oplus [\lambda_{(21,0)}] \oplus [\lambda_{(0,21)}]. \label{theta-E(24)}
\end{eqnarray}

Note that these dual canonical endomorphisms depend only on the existence of a cell system $W$ for each graph $\mathcal{G}$, but not on the choice of cell system since Lemma \ref{Lemma-existence_of_intertwiner} and Proposition \ref{Prop:intertwinerV} did not depend on this choice. Where we have found two inequivalent solutions, the computations below show that either choice will give the same $M$-$N$ graph, since the computations in these particular cases only depend on the dual canonical endomorphism $\theta$. Similarly, even if there exists other solutions for the cells $W$ for the $\mathcal{D}$, $\mathcal{D}^{\ast}$ and $\mathcal{E}_1^{(12)}$ graphs, these will not give any new $M$-$N$ graphs either. It is conceivable however that in certain situations, for $SU(n)$, $n>3$, the $M$-$N$ graph will depend on the connection and not just on the GHJ graph.

\emph{Remark.}
For $SU(2)$ it was shown in \cite{evans:2002} that the modular invariant $Z$ can be realized from a subfactor with a dual canonical endomorphism of the form
\begin{equation} \label{eqn:theta=Z(mu,mu-bar)}
[\theta] = \bigoplus_{\mu} Z_{\mu, \overline{\mu}} [\mu],
\end{equation}
where the direct summation is over all $\mu$ even.
This raises the question of whether all the $SU(3)$ modular invariants can be realized from some subfactor with dual canonical endomorphism $\theta$ of the form (\ref{eqn:theta=Z(mu,mu-bar)}), where now allow $\mu$ to be of any colour?
For the $\mathcal{A}^{(n)\ast}$ graphs the $\theta$ given in (\ref{theta-A(n)star}) is automatically in the form (\ref{eqn:theta=Z(mu,mu-bar)}), where $Z$ is the conjugate modular invariant $Z_{\mathcal{A}^{(n)\ast}} = C$.
For the $\mathcal{A}^{(n)}$ graphs, if we choose the $M$-$N$ morphism $[\overline{a}]$ to be $[\iota \lambda_{(p,0)}]$, where $p = \lfloor (n-3)/2 \rfloor$, the sector $[a \overline{a}]$ gives $[\lambda_{(0,0))}] \oplus [\lambda_{(1,1)}] \oplus [\lambda_{(2,2)}] \oplus \cdots \oplus [\lambda_{(p,p)}]$, and we obtain a dual canonical endomorphism $[\theta] = [a \overline{a}] = \bigoplus_{\mu} Z_{\mu, \overline{\mu}} [\mu]$, where the direct summation is over all $\mu$ (of any colour) and $Z$ is the identity modular invariant $Z_{\mathcal{A}^{(n)}} = I$. \\

For each of the $\mathcal{ADE}$ graphs (with the exception of $\mathcal{E}_4^{(12)}$) we have shown the existence of a braided subfactor $N \subset M$ with dual canonical endomorphisms $\theta$ given by (\ref{theta-A(n)})-(\ref{theta-E(24)}). By the $\alpha$-induction of \cite{bockenhauer/evans:1998, bockenhauer/evans:1999i, bockenhauer/evans:1999ii}, a matrix $Z$ can be defined by $Z_{\lambda,\mu} = \langle \alpha^+_{\lambda}, \alpha^-_{\mu} \rangle$, $\lambda, \mu \in {}_N \mathcal{X}_N$. If the braiding is non-degenerate, $Z$ is a modular invariant mass matrix.

For the dual canonical endomorphisms $\theta$ in (\ref{theta-A(n)})-(\ref{theta-E(24)}), what is the corresponding $M$-$N$ system or Cappelli-Itzykson-Zuber graph which classifies the modular invariant? And what is the corresponding modular invariant? For $\mathcal{A}^{(n)}$ the $M$-$M$, $M$-$N$ and $N$-$N$ systems are all equal since $N=M$. Subfactors given by conformal inclusions were considered in \cite{bockenhauer/evans:1999i, bockenhauer/evans:1999ii}. Those conformal inclusions which have $SU(3)$ invariants give identical dual canonical endomorphisms $\theta$ to those computed above. The $M$-$N$ system was computed for conformal inclusions with corresponding modular invariants associated to the graphs $\mathcal{D}^{(6)}$ and $\mathcal{E}^{(8)}$ in \cite{bockenhauer/evans:1999i}, and to $\mathcal{E}_1^{(12)}$ and $\mathcal{E}^{(24)}$ in \cite{bockenhauer/evans:1999ii}. The $M$-$N$ system was also computed in \cite{bockenhauer/evans:1999i} for the inclusion with the $\mathcal{D}^{(n)}$ dual canonical endomorphism (\ref{theta-D(n)}) for $n \equiv 0 \textrm{ mod } 3$, and in \cite{bockenhauer/evans:2001} for the inclusion with the $\mathcal{E}_2^{(12)}$ dual canonical endomorphism (\ref{theta-E2(12)}), which do not come from conformal inclusions. For each of these graphs, the graph of the $M$-$N$ system and the $\alpha$-graph can both be identified with the original graph itself, and the modular invariant is that associated with the original graph. We compute the $M$-$N$ graph for the remaining $\theta$'s. The proof for the case of $\mathcal{E}_2^{(12)}$ was not published in \cite{bockenhauer/evans:2001}, so we produce a proof using our method here. Knowledge of the dual canonical endomorphism $\theta$ is not usually sufficient to determine the $M$-$N$ graph, but we can utilize the fact that the list of $SU(3)$ modular invariants is complete. For an $\mathcal{ADE}$ graph $\mathcal{G}$ with Coxeter number $n$, the basic method is to compute $\langle \iota \lambda, \iota \mu \rangle$ for representations $\lambda$, $\mu$ on $\mathcal{A}^{(n)}$, and decompose into irreducibles. Sometimes there is an ambiguity about the decomposition, e.g. if $\langle \iota \lambda, \iota \lambda \rangle = 4$ then we could have $\iota \lambda = 2 \iota \lambda^{(1)}$ or $\iota \lambda = \iota \lambda^{(1)} + \iota \lambda^{(2)} + \iota \lambda^{(3)} + \iota \lambda^{(4)}$ where $\iota \lambda^{(i)}$, $i=1,2,3,4$, are irreducible sectors. By \cite[Cor. 6.13]{bockenhauer/evans/kawahigashi:1999}, $\sharp {}_M \mathcal{X}_N = \mathrm{tr}(Z)$ for some modular invariant $Z$, and therefore, since we have a complete list of $SU(3)$ modular invariants, we can eliminate any particular decomposition if the total number of irreducible sectors obtained does not agree with the trace of any of the modular invariants (\ref{Z(A)})-(\ref{Z(E24)}). We compute the trace for all the modular invariants at level $k$ in the following lemma:

\begin{Lemma}
The traces of the level $k$ modular invariants $Z$ are
\begin{eqnarray}
\mathrm{tr}(Z_{\mathcal{A}^{(k+3)}}) & = & \frac{1}{2} (k+1)(k+2), \\
\mathrm{tr}(Z_{\mathcal{D}^{(k+3)}}) & = & \frac{1}{6} (k+1)(k+2) + c_k, \\
\mathrm{tr}(Z_{\mathcal{A}^{(k+3)\ast}}) & = & \lfloor \frac{k+2}{2} \rfloor, \\
\mathrm{tr}(Z_{\mathcal{D}^{(k+3)\ast}}) & = & 3 \lfloor \frac{k+2}{2} \rfloor, \\
\mathrm{tr}(Z_{\mathcal{E}^{(8)}}) & = & 12, \label{tr(Z(E8))} \\
\mathrm{tr}(Z_{\mathcal{E}^{(8)\ast}}) & = & 4, \label{tr(Z(E8star))} \\
\mathrm{tr}(Z_{\mathcal{E}^{(12)}}) & = & 12, \label{tr(Z(E1-12))} \\
\mathrm{tr}(Z_{\mathcal{E}_{MS}^{(12)\ast}}) & = & 11, \label{tr(Z(E4-12))} \\
\mathrm{tr}(Z_{\mathcal{E}_{MS}^{(12)}}) & = & 17, \label{tr(Z(E5-12))} \\
\mathrm{tr}(Z_{\mathcal{E}^{(24)}}) & = & 24, \label{tr(Z(E24))}
\end{eqnarray}
where $c_k = 0$ if $k \not \equiv 0 \textrm{ mod } 3$, $c_{3m} = 2/3$ for $m \in \mathbb{N}$ and $\lfloor x \rfloor$ denotes the largest integer less than or equal to $x$.
\end{Lemma}
\emph{Proof}

For the $\mathcal{A}$ graphs, $\mathrm{tr}(Z_{\mathcal{A}^{(k+3)}})$ is given by the number of vertices of $\mathcal{A}^{(k+3)}$, which is $1 + 2 + 3 + \cdots + k+1 = (k+1)(k+2)/2$. For $k \not \equiv 0 \textrm{ mod } 3$, the diagonal terms in $Z_{\mathcal{D}^{(k+3)}}$ are given by the 0-coloured vertices of $\mathcal{A}^{(k+3)}$, so $\mathrm{tr}(Z_{\mathcal{D}^{(k+3)}})$ is $\mathrm{tr}(Z_{\mathcal{A}^{(k+3)}})/3$. For $k \equiv 0 \textrm{ mod } 3$ the 0-coloured vertices of $\mathcal{A}^{(k+3)}$ again give the diagonal terms in $Z_{\mathcal{D}^{(k+3)}}$ but the number of 0-coloured vertices of $\mathcal{A}^{(k+3)}$ is now one greater than the number of 1,2-coloured vertices. The trace of $Z_{\mathcal{A}^{(k+3)\ast}}$ is given by the number of ``diagonal'' elements $\mu = \overline{\mu}$ of $\mathcal{A}^{(k+3)}$, which is $\lfloor k+2/2 \rfloor$. For the $\mathcal{D}^{\ast}$ graphs, when $k \not \equiv 0 \textrm{ mod } 3$, the trace is given by the number of vertices $\mu = (\mu_1, \mu_2)$ of $\mathcal{A}^{(k+3)}$ such that $A^{(n-3)(\mu_1-\mu_2)} \mu = \overline{\mu}$. For the 0-coloured vertices this is the number of diagonal elements, whilst for the 1,2-coloured vertices this is where $A \mu = \overline{\mu}$ or $A^2 \mu = \overline{\mu}$, depending on the parity of $n$. In each case the number of such vertices is $\lfloor k+2/2 \rfloor$. For $k \equiv 0 \textrm{ mod } 3$ the trace is again given by a third of the number of vertices of $\mathcal{A}^{(k+3)}$ which satisfy each of the following $\mu = \overline{\mu}$, $A \mu = \overline{A^2 \mu}$, $A^2 = \overline{A \mu}$, $\mu = \overline{A \mu}$, $A \mu = \overline{\mu}$, $A^2 \mu = \overline{A^2 \mu}$, $\mu = \overline{A^2 \mu}$, $A \mu = \overline{A \mu}$ and $A^2 \mu = \overline{\mu}$. The first three equalities are satisfied when $\mu = \overline{\mu}$, the second three when $A \mu = \overline{\mu}$ and the last three when $A^2 \mu = \overline{\mu}$. So we have $\mathrm{tr}(Z_{\mathcal{D}^{(k+3)\ast}}) = 3 \lfloor k+2/2 \rfloor$ also.
The computations of $\mathrm{tr}(Z_{\mathcal{E}})$ for the exceptional invariants is clear from inspection of the modular invariant.
\hfill
$\Box$

\begin{Lemma}
The trace of the modular invariants at level $k$ are all different.
\end{Lemma}
\emph{Proof}
For level 5 we have $\mathrm{tr}(\mathcal{A}^{(8)}) = 21$, $\mathrm{tr}(\mathcal{D}^{(8)}) = 7$, $\mathrm{tr}(\mathcal{A}^{(8)\ast}) = 3$ and $\mathrm{tr}(\mathcal{D}^{(8)\ast}) = 9$, and compare these with (\ref{tr(Z(E8))}) and (\ref{tr(Z(E8star))}). For level 9, $\mathrm{tr}(\mathcal{A}^{(12)}) = 55$, $\mathrm{tr}(\mathcal{D}^{(12)}) = 19$, $\mathrm{tr}(\mathcal{A}^{(12)\ast}) = 5$ and $\mathrm{tr}(\mathcal{D}^{(12)\ast}) = 15$, and compare these with (\ref{tr(Z(E1-12))})-(\ref{tr(Z(E5-12))}). For level 21 we compare $\mathrm{tr}(\mathcal{A}^{(24)}) = 253$, $\mathrm{tr}(\mathcal{D}^{(24)}) = 85$, $\mathrm{tr}(\mathcal{A}^{(24)\ast}) = 11$ and $\mathrm{tr}(\mathcal{D}^{(24)\ast}) = 33$ with (\ref{tr(Z(E24))}). For all other levels we need to compare the modular invariants for the $\mathcal{A}$, $\mathcal{D}$, $\mathcal{A^{\ast}}$ and $\mathcal{D^{\ast}}$ graphs.

Comparing the $\mathcal{A}$ and $\mathcal{D}$ modular invariants, the traces can only be equal if $3(k+1)(k+2) = (k+1)(k+2) + 6c_k$. For $k \equiv 0 \textrm{ mod } 3$ this gives $k=0,-3$, whilst if $k \not \equiv 0 \textrm{ mod } 3$ we obtain $k=-1,-2$. So these traces cannot be equal except when $k = 0$, but the graphs $\mathcal{A}^{(3)}$ and $\mathcal{D}^{(3)}$ are both a single vertex. Comparing $\mathcal{A}$-$\mathcal{A}^{\ast}$, the traces are only equal if $(k+1)(k+2)=2 \lfloor (k+2)/2 \rfloor$. For even $k$ this gives solutions $k=0,-4$, but when $k = 0$ the graph $\mathcal{A}^{(3)\ast}$ is also just a single vertex, so identical to the graph $\mathcal{A}^{(3)}$. For $k$ odd we have $k=-1$. Next, comparing $\mathcal{A}$-$\mathcal{D}^{\ast}$, the traces are only equal if $(k+1)(k+2)=6 \lfloor (k+2)/2 \rfloor$. For $k$ even this gives solutions $k= \pm 2$, but for $k=2$ the graph $\mathcal{D}^{(5)\ast}$ is identical to $\mathcal{A}^{(5)}$. For $k$ odd we obtain solutions $k = -3,1$, but we again have for $k=1$ that the graphs $\mathcal{D}^{(4)\ast}$ and $\mathcal{A}^{(4)}$ are the same. We now compare $\mathcal{D}$-$\mathcal{A}^{\ast}$. When $k \equiv 0 \textrm{ mod } 3$, the traces are equal only if $(k+1)(k+2) + 4 = 6 \lfloor (k+2)/2 \rfloor = 6 \lfloor k/2 \rfloor + 6$, so we have the quadratic $k^2 + 3(k - 2 \lfloor k/2 \rfloor) = 0$. When $k$ is even we have only the solution $k=0$, whilst when $k$ is odd this gives $k^2 = -3$. When $k \not \equiv 0 \textrm{ mod } 3$, we obtain instead the quadratic $k^2 + 3(k - 2 \lfloor k/2 \rfloor) - 4 = 0$. For even $k$ this gives the solutions $k = \pm 2$, but we notice that the graphs $\mathcal{D}^{(5)}$ and $\mathcal{A}^{(5)\ast}$ are the same, whilst for odd $k$ we have the solutions $k = \pm 1$, but we again see that the graphs $\mathcal{D}^{(4)}$ and $\mathcal{A}^{(4)\ast}$ are the same. Comparing $\mathcal{D}$-$\mathcal{D}^{\ast}$ we now obtain the quadratic equations $k^2 + 3(k - 6 \lfloor k/2 \rfloor) - 14 = 0$, $k^2 + 3(k - 6 \lfloor k/2 \rfloor) - 18 = 0$ for $k \equiv 0 \textrm{ mod } 3$, $k \not \equiv 0 \textrm{ mod } 3$ respectively. Neither of these equations has integer solutions for odd or even $k$. Finally, comparing the $\mathcal{A}^{\ast}$ and $\mathcal{D}^{\ast}$ modular invariants, the traces are only equal if $\lfloor (k+2)/2 \rfloor = 3 \lfloor (k+2)/2 \rfloor$, giving $\lfloor (k+2)/2 \rfloor = 0$ which has solutions $k = -2,-3$.
\hfill
$\Box$

Since the traces of the modular invariants at any level are all different, once we have found the number of irreducible sectors, we can identify the corresponding modular invariant. There may however still be an ambiguity with regard to the fusion rules that these irreducible sectors satisfy, with different seemingly possible fusion rules giving different nimrep graphs for the $M$-$N$ system. However, we know that the nimrep must have spectrum $S_{\lambda,\nu}/S_{\lambda,0}$ with multiplicity determined by the diagonal part ${Z_{\lambda,\lambda}}$ of the modular invariant. It turns out that the consideration of the trace and the eigenvalues is sufficient to compute the $M$-$N$ graphs for $\mathcal{A}^{(12)\ast}$, $\mathcal{D}^{(12)\ast}$, $\mathcal{E}_2^{(12)}$, $\mathcal{E}_4^{(12)}$ and $\mathcal{E}_5^{(12)}$, and identify the corresponding modular invariant. The results are summarized in Table \ref{table:summary_of_M-N_system}. We will say that an irreducible sector $[\iota \lambda_{(\mu_1,\mu_2)}]$ such that $\mu_1 + \mu_2 = m$ appears at tier $m$.

\subsection{$\mathcal{E}^{(8)\ast}$}

\begin{figure}[tb]
\begin{center}
  \includegraphics[width=30mm]{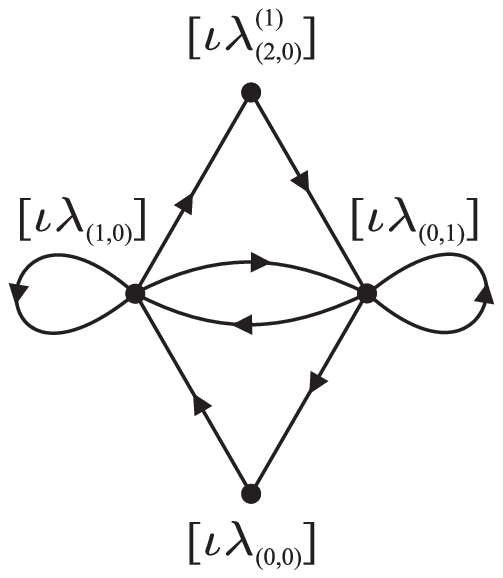}
 \caption{$M$-$N$ graph for the $\mathcal{E}^{(8)\ast}$ $SU(3)$-GHJ subfactor}\label{fig:fig-GHJ-14}
\end{center}
\end{figure}

For the graph $\mathcal{E}^{(8)\ast}$, we have $[\theta] = [\lambda_{(0,0)}] \oplus [\lambda_{(2,1)}] \oplus [\lambda_{(1,2)}] \oplus [\lambda_{(2,2)}] \oplus [\lambda_{(5,0)}] \oplus [\lambda_{(0,5)}]$. Then computing $\langle \iota \lambda, \iota \mu \rangle = \langle \lambda, \theta \mu \rangle$ (by Frobenius reciprocity) for $\lambda$, $\mu$ on $\mathcal{A}^{(8)}$, we find $\langle \iota \lambda, \iota \lambda \rangle = 1$ and $\langle \iota \lambda, \iota \mu \rangle = 0$ for $\lambda, \mu = \lambda_{(0,0)}, \lambda_{(1,0)}, \lambda_{(0,1)}$. At tier 2 we have $\langle \iota \lambda_{(2,0)}, \iota \lambda_{(2,0)} \rangle = 2$, $\langle \iota \lambda_{(2,0)}, \iota \lambda_{(1,0)} \rangle = 1$ and $\langle \iota \lambda_{(2,0)}, \iota \mu \rangle = 0$ for $\mu = \lambda_{(0,1)}, \lambda_{(0,0)}$. So $[\iota \lambda_{(2,0)}] = [\iota \lambda_{(1,0)}] \oplus [\iota \lambda_{(2,0)}^{(1)}]$. Since $\langle \iota \lambda_{(0,2)}, \iota \lambda_{(0,2)} \rangle = \langle \iota \lambda_{(0,2)}, \iota \lambda_{(2,0)} \rangle = 2$ we have $[\iota \lambda_{(0,2)}] = [\iota \lambda_{(2,0)}]$. Lastly at tier 2 we have $\langle \iota \lambda_{(1,1)}, \iota \lambda_{(1,1)} \rangle = 2$ and $\langle \iota \lambda_{(1,1)}, \iota \lambda_{(1,0)} \rangle = \langle \iota \lambda_{(1,1)}, \iota \lambda_{(0,1)} \rangle = 1$, giving $[\iota \lambda_{(1,1)}] = [\iota \lambda_{(1,0)}] \oplus [\iota \lambda_{(0,1)}]$. At tier 3 we have $\langle \iota \lambda_{(3,0)}, \iota \lambda_{(3,0)} \rangle = \langle \iota \lambda_{(3,0)}, \iota \lambda_{(0,2)} \rangle = 2$, so $[\iota \lambda_{(3,0)}] = [\iota \lambda_{(0,2)}]$. Similarly $[\iota \lambda_{(0,3)}] = [\iota \lambda_{(2,0)}]$. For $\iota \lambda_{(2,1)}$ we find $\langle \iota \lambda_{(2,1)}, \iota \lambda_{(2,1)} \rangle = 2$ and $\langle \iota \lambda_{(2,1)}, \iota \lambda_{(0,0)} \rangle = \langle \iota \lambda_{(2,1)}, \iota \lambda_{(1,0)} \rangle = 1$, giving $[\iota \lambda_{(2,1)}] = [\iota \lambda_{(0,0)}] \oplus [\iota \lambda_{(1,0)}]$ and similarly $[\iota \lambda_{(1,2)}] = [\iota \lambda_{(0,0)}] \oplus [\iota \lambda_{(0,1)}]$. So no new irreducibles appear at tier 3. No new irreducible sectors appear at the other tiers either, so we have 4 irreducible sectors $[\iota \lambda_{(0,0)}]$, $[\iota \lambda_{(1,0)}]$, $[\iota \lambda_{(0,1)}]$ and $[\iota \lambda_{(2,0)}^{(1)}]$.
We now compute the sector products of these irreducible sectors with the $M$-$N$ sector $[\rho] = [\lambda_{(1,0)}]$. It is easy to compute $[\iota \lambda_{(0,0)}] [\rho] = [\iota \lambda_{(1,0)}]$, $[\iota \lambda_{(1,0)}] [\rho] = [\iota \lambda_{(0,1)}] \oplus [\iota \lambda_{(2,0)}] = [\iota \lambda_{(0,1)}] \oplus [\iota \lambda_{(1,0)}] \oplus [\iota \lambda_{(2,0)}^{(1)}]$ and $[\iota \lambda_{(0,1)}] [\rho] = [\iota \lambda_{(0,0)}] \oplus [\iota \lambda_{(1,0)}] \oplus [\iota \lambda_{(0,1)}]$. We can invert these formula to obtain $[\iota \lambda_{(2,0)}^{(1)}] = [\iota \lambda_{(2,0)}] \ominus [\iota \lambda_{(1,0)}]$, and so $[\iota \lambda_{(2,0)}^{(1)}] [\rho] = [\iota \lambda_{(1,1)}] \oplus [\iota \lambda_{(3,0)}] \ominus ([\iota \lambda_{(2,0)}] \oplus [\iota \lambda_{(0,1)}]) = [\iota \lambda_{(0,1)}]$. Then we see that the multiplication graph for $[\rho]$ is the original graph $\mathcal{E}^{(8)\ast}$ itself, illustrated in Figure \ref{fig:fig-GHJ-14}, and the modular invariant associated to $\theta$ is $Z_{\mathcal{E}^{(8)\ast}}$.

\subsection{$\mathcal{E}_2^{(12)}$}

For the graph $\mathcal{E}_2^{(12)}$, we have $[\theta] = [\lambda_{(0,0)}] \oplus 2 [\lambda_{(2,2)}] \oplus [\lambda_{(4,1)}] \oplus [\lambda_{(1,4)}] \oplus 2 [\lambda_{(5,2)}] \oplus 2 [\lambda_{(2,5)}] \oplus [\lambda_{(4,4)}] \oplus [\lambda_{(9,0)}] \oplus [\lambda_{(0,9)}]$. We have $\langle \iota \lambda, \iota \lambda \rangle = 1$ and $\langle \iota \lambda, \iota \mu \rangle = 0$ for all $\lambda, \mu \in \{ \lambda_{(0,0)}, \lambda_{(1,0)}, \lambda_{(0,1)} \}$. At tier 2 we have $\langle \iota \lambda, \iota \lambda \rangle = 3$ and $\langle \iota \lambda, \iota \mu \rangle = 0$ for $\lambda = \lambda_{(2,0)}, \lambda_{(1,1)}, \lambda_{(0,2)}$, $\mu = \lambda_{(0,0)}, \lambda_{(1,0)}, \lambda_{(0,1)}$. Then $\lambda_{(2,0)}$, $\lambda_{(1,1)}$, $\lambda_{(0,2)}$ decompose into irreducibles as
\begin{eqnarray}
{} [\iota \lambda_{(2,0)}] & = & [\iota \lambda^{(1)}_{(2,0)}] \oplus [\iota \lambda^{(2)}_{(2,0)}] + [\iota \lambda^{(3)}_{(2,0)}], \label{eqn:E2(12)-lambda(2,0)} \\
{} [\iota \lambda_{(1,1)}] & = & [\iota \lambda^{(1)}_{(1,1)}] \oplus [\iota \lambda^{(2)}_{(1,1)}] \oplus [\iota \lambda^{(3)}_{(1,1)}], \label{eqn:E2(12)-lambda(1,1)} \\
{} [\iota \lambda_{(0,2)}] & = & [\iota \lambda^{(1)}_{(0,2)}] \oplus [\iota \lambda^{(2)}_{(0,2)}] \oplus [\iota \lambda^{(3)}_{(0,2)}].
\end{eqnarray}
At tier 3 we find $\langle \iota \lambda_{(3,0)}, \iota \lambda_{(3,0)} \rangle = \langle \iota \lambda_{(3,0)}, \iota \lambda_{(1,1)} \rangle = 3$ so that $[\iota \lambda_{(3,0)}] = [\iota \lambda_{(1,1)}]$, and similarly $[\iota \lambda_{(0,3)}] = [\iota \lambda_{(1,1)}]$. From $\langle \iota \lambda_{(2,1)}, \iota \lambda_{(2,1)} \rangle = 7$, $\langle \iota \lambda_{(2,1)}, \iota \lambda_{(1,0)} \rangle = 2$ and $\langle \iota \lambda_{(2,1)}, \iota \lambda_{(0,2)} \rangle = 3$, and similarly for $\iota \lambda_{(1,2)}$, we obtain
\begin{eqnarray}
{} [\iota \lambda_{(2,1)}] & = & 2 [\iota \lambda_{(1,0)}] \oplus [\iota \lambda^{(1)}_{(0,2)}] \oplus [\iota \lambda^{(2)}_{(0,2)}] + [\iota \lambda^{(3)}_{(0,2)}], \label{eqn:E2(12)-lambda(1,2)} \\
{} [\iota \lambda_{(1,2)}] & = & 2 [\iota \lambda_{(0,1)}] \oplus [\iota \lambda^{(1)}_{(2,0)}] \oplus [\iota \lambda^{(2)}_{(2,0)}] \oplus [\iota \lambda^{(3)}_{(2,0)}], \nonumber
\end{eqnarray}
and no new irreducible sectors appear at tier 3. Then we have twelve irreducible sectors $[\iota \lambda_{(0,0)}]$, $[\iota \lambda_{(1,0)}]$, $[\iota \lambda_{(0,1)}]$, $[\iota \lambda^{(i)}_{(2,0)}]$, $[\iota \lambda^{(i)}_{(1,1)}]$, $[\iota \lambda^{(i)}_{(0,2)}]$ for $i = 1,2,3$, and the corresponding modular invariant must be $Z_{\mathcal{E}^{(12)}}$ since $\mathrm{tr}(Z_{\mathcal{E}^{(12)}}) = 12$.

We now look at the fusion rules that these irreducible sectors satisfy. With $\rho = \lambda_{(1,0)}$, we have $[\iota \lambda_{(0,0)}] [\rho] = [\iota \lambda_{(1,0)}]$,
\begin{equation}\label{eqn:E2(12)-rho_x_lambda(1,0)}
[\iota \lambda_{(1,0)}] [\rho] = [\iota \lambda_{(0,1)}] \oplus [\iota \lambda_{(2,0)}] = [\iota \lambda_{(0,1)}] \oplus [\iota \lambda^{(1)}_{(2,0)}] \oplus [\iota \lambda^{(2)}_{(2,0)}] \oplus [\iota \lambda^{(3)}_{(2,0)}],
\end{equation}
and similarly $[\iota \lambda_{(0,1)}] [\rho] = [\iota \lambda_{(0,0)}] \oplus [\iota \lambda^{(1)}_{(1,1)}] \oplus [\iota \lambda^{(2)}_{(1,1)}] \oplus [ \iota \lambda^{(3)}_{(1,1)}]$. Since $[\iota \lambda_{(2,0)}] [\rho] = [\iota \lambda_{(1,1)}] \oplus [\iota \lambda_{(3,0)}] = 2 [\iota \lambda^{(1)}_{(1,1)}] \oplus 2 [\iota \lambda^{(2)}_{(1,1)}] \oplus 2 [ \iota \lambda^{(3)}_{(1,1)}]$, we obtain $([\iota \lambda^{(1)}_{(2,0)}] [\rho]) \oplus ([\iota \lambda^{(2)}_{(2,0)}] [\rho]) \oplus ([\iota \lambda^{(3)}_{(2,0)}] [\rho]) = 2 [\iota \lambda^{(1)}_{(1,1)}] \oplus 2 [\iota \lambda^{(2)}_{(1,1)}] \oplus+ 2 [\iota \lambda^{(3)}_{(1,1)}]$.

\begin{figure}[tb]
\begin{center}
  \includegraphics[width=80mm]{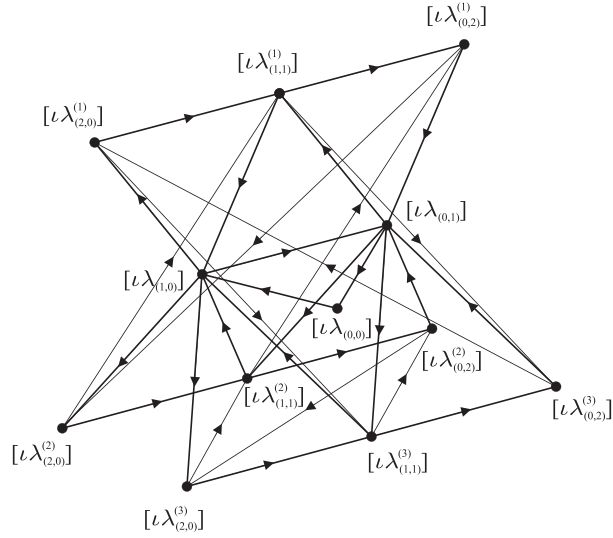}
 \caption{$M$-$N$ graph for the $\mathcal{E}_2^{(12)}$ $SU(3)$-GHJ subfactor}\label{fig:fig-GHJ-9}
\end{center}
\end{figure}

We now use a similar argument to that in \cite[$\S$2.4]{bockenhauer/evans:1999i}. The statistical dimension of the positive energy representation $(\mu_1, \mu_2)$ of $SU(3)_9$ is given by the Perron-Frobenius eigenvector for the graph $\mathcal{A}^{(12)}$: $d_{(\mu_1,\mu_2)} = [\mu_1 +1][\mu_2 +1][\mu_1 + \mu_2 +2]/[2]$. Then from (\ref{eqn:E2(12)-rho_x_lambda(1,0)}) we obtain $d_{(2,0)}^{(1)} + d_{(2,0)}^{(2)} + d_{(2,0)}^{(3)} = d_{(1,0)}^2 - d_{(1,0)} = [3]^3 - [3] = [3][4]/[2]$, where $d_{(2,0)}^{(i)} = d_{\iota \lambda_{(2,0)}^{(i)}}$. We may then assume without loss of generality that $d_{(2,0)}^{(1)} < [3][4]/(3[2]) = [2][3]/[4]$. Then since $([2][3]/[4])^2 \approx 2.488 < 3$, $[\overline{\iota \lambda_{(2,0)}^{(1)}}] [\iota \lambda_{(2,0)}^{(1)}]$ decomposes into at most two irreducible $N$-$N$ sectors. Then $\langle \iota \lambda_{(2,0)}^{(1)} \circ \rho, \iota \lambda_{(2,0)}^{(1)} \circ \rho \rangle = \langle \rho \circ \overline{\rho},  \overline{\iota \lambda_{(2,0)}^{(1)}} \circ \iota \lambda_{(2,0)}^{(1)} \rangle \leq 2$. So $[\iota \lambda_{(2,0)}^{(1)}] [\rho]$ cannot contain an irreducible sector with multiplicity greater than one. Since, by (\ref{eqn:E2(12)-lambda(2,0)}) and (\ref{eqn:E2(12)-lambda(1,2)}), $\langle \iota \lambda_{(2,0)}^{(1)} \circ \rho, \iota \lambda_{(1,1)} \rangle = \langle \iota \lambda_{(2,0)}^{(1)}, \iota \lambda_{(1,1)} \circ \overline{\rho} \rangle = \langle \iota \lambda_{(2,0)}^{(1)}, \iota \lambda_{(0,1)} + \iota \lambda_{(2,0)} + \iota \lambda_{(1,2)} \rangle = 2$, using (\ref{eqn:E2(12)-lambda(1,1)}) we may assume, again without loss of generality, that
$$[\iota \lambda_{(2,0)}^{(1)}] [\rho] = [\iota \lambda_{(1,1)}^{(1)}] \oplus [\iota \lambda_{(1,1)}^{(2)}].$$
Since $[\iota \lambda_{(1,0)}] [\rho] \supset [\iota \lambda_{(2,0)}^{(1)}]$ and $\langle \iota \lambda_{(1,0)}, \iota \lambda_{(2,0)}^{(1)} \circ \overline{\rho} \rangle = \langle \iota \lambda_{(1,0)} \circ \rho, \iota \lambda_{(2,0)}^{(1)} \rangle > 0$, then $[\iota \lambda_{(2,0)}^{(1)}] [\overline{\rho}] \supset [\iota \lambda_{(1,0)}]$. Then since $\langle \iota \lambda_{(2,0)}^{(1)} \circ \overline{\rho}, \iota \lambda_{(2,0)}^{(1)} \circ \overline{\rho} \rangle = \langle \iota \lambda_{(2,0)}^{(1)} \circ \rho, \iota \lambda_{(2,0)}^{(1)} \circ \rho \rangle = 2$, we have $[\iota \lambda_{(2,0)}^{(1)}] [\overline{\rho}] = [\iota \lambda_{(1,0)}] \oplus [\iota \lambda_{(0,2)}^{(j)}]$, for $j \in \{ 1,2,3 \}$. By a similar argument we may also assume that $[\iota \lambda_{(2,0)}^{(1)}]$ has statistical dimension $< [2][3]/[4]$, and using $[\overline{\rho}]$ instead of $[\rho]$, we find $[\iota \lambda_{(0,2)}^{(1)}] [\rho] = [\iota \lambda_{(0,1)}] \oplus [\iota \lambda_{(2,0)}^{(j')}]$, and have the freedom to set $j' = 3$. Then we also have $[\iota \lambda_{(0,2)}^{(j)}] [\rho] \supset [\iota \lambda_{(0,1)}]$ for $j=2,3$ and $([\iota \lambda_{(0,2)}^{(2)}] \oplus [\iota \lambda_{(0,2)}^{(3)}]) [\rho] = 2 [\iota \lambda_{(0,1)}] \oplus [\iota \lambda_{(2,0)}^{(1)}] \oplus [\iota \lambda_{(2,0)}^{(2)}]$.
From $[\iota \lambda_{(1,1)}] [\rho]$ we obtain $([\iota \lambda_{(1,1)}^{(1)}] \oplus [\iota \lambda_{(1,1)}^{(2)}] \oplus [\iota \lambda_{(1,1)}^{(3)}]) [\rho] = 3 [\iota \lambda_{(1,0)}] \oplus 2 [\iota \lambda_{(0,2)}^{(1)}] \oplus 2 [\iota \lambda_{(0,2)}^{(2)}] \oplus 2 [\iota \lambda_{(0,2)}^{(3)}]$ and since $[\iota \lambda_{(1,0)}] [\overline{\rho}] = [\iota \lambda_{(2,0)}] \oplus [\iota \lambda_{(1,1)}] = [\iota \lambda_{(2,0)}] \oplus [\iota \lambda_{(1,1)}^{(1)}] \oplus [\iota \lambda_{(1,1)}^{(2)}] \oplus [\iota \lambda_{(1,1)}^{(3)}]$ then $\langle \iota \lambda_{(1,1)}^{(j)} \circ \rho, \iota \lambda_{(1,0)} \rangle = \langle \iota \lambda_{(1,1)}^{(j)}, \iota \lambda_{(1,0)} \circ \overline{\rho} \rangle = 1$ and $[\iota \lambda_{(1,1)}^{(j)}] [\rho] \supset [\iota \lambda_{(1,0)}]$ for $j=1,2,3$.
There is still some ambiguity surrounding the decompositions of $[\iota \lambda_{(2,0)}^{(j)}] [\rho]$, $[\iota \lambda_{(1,1)}^{(j)}] [\rho]$ and $[\iota \lambda_{(0,2)}^{(j)}] [\rho]$, for $j=2,3$. Computing the eigenvalues of the nimrep graphs for the different possibilities, we find that the only nimrep graph which has eigenvalues $S_{\rho \mu}/S_{0 \mu}$ with multiplicities given by the diagonal entry $Z_{\mu,\mu}$ of the modular invariant is that for: $[\iota \lambda_{(2,0)}^{(j)}] [\rho] = [\iota \lambda_{(1,1)}^{(j)}] \oplus [\iota \lambda_{(1,1)}^{(j+1)}]$, $[\iota \lambda_{(1,1)}^{(j)}] [\rho] = [\iota \lambda_{(0,2)}^{(l)}] \oplus [\iota \lambda_{(0,2)}^{(l+1)}]$ and $[\iota \lambda_{(0,2)}^{(j)}] [\rho] = [\iota \lambda_{(0,1)}] \oplus [\iota \lambda_{(2,0)}^{(j+1)}]$ for $j=1,2,3$, $l \in \{ 1,2,3 \}$. The nimrep graph is the same for any choice of $l=1,2,3$, up to a relabeling of the irreducible representations $[\iota \lambda_{(2,0)}^{(j)}]$, $[\iota \lambda_{(1,1)}^{(j)}]$ and $[\iota \lambda_{(0,2)}^{(j)}]$, and the graph is just the graph $\mathcal{E}_2^{(12)}$ itself, illustrated in Figure \ref{fig:fig-GHJ-9}. The associated modular invariant is $Z_{\mathcal{E}^{(12)}}$.

\subsection{$\mathcal{E}_4^{(12)}$}

Warning: the existence of the $SU(3)$-Goodman-de la Harpe-Jones subfactor which gives the dual canonical endomorphism for $\mathcal{E}_4^{(12)}$ has not been shown yet by us.

For $\mathcal{E}_4^{(12)}$, we suppose $[\theta] = [\lambda_{(0,0)}] \oplus [\lambda_{(2,2)}] \oplus [\lambda_{(4,1)}] \oplus [\lambda_{(1,4)}] \oplus [\lambda_{(5,2)}] \oplus [\lambda_{(2,5)}] \oplus [\lambda_{(4,4)}] \oplus [\lambda_{(9,0)}] \oplus [\lambda_{(0,9)}]$. Then computing $\langle \iota \lambda, \iota \mu \rangle = \langle \lambda, \theta \mu \rangle$ for $\lambda$, $\mu$ on $\mathcal{A}^{(12)}$, we find $\langle \iota \lambda, \iota \lambda \rangle = 1$ for $\lambda = \lambda_{(0,0)}, \lambda_{(1,0)}, \lambda_{(0,1)}$. At tier 2 we have $\langle \iota \lambda, \iota \lambda \rangle = 2$ and $\langle \iota \lambda, \iota \mu \rangle = 0$ for $\lambda = \lambda_{(2,0)}, \lambda_{(1,1)}, \lambda_{(0,2)}$, $\mu = \lambda_{(0,0)}, \lambda_{(1,0)}, \lambda_{(0,1)}$. Then $[\lambda_{(2,0)}]$, $[\lambda_{(1,1)}]$, $[\lambda_{(0,2)}]$ decompose into irreducibles as
\begin{eqnarray}
{} [\iota \lambda_{(2,0)}] & = & [\iota \lambda^{(1)}_{(2,0)}] \oplus [\iota \lambda^{(2)}_{(2,0)}], \label{eqn:E4(12)-lambda(2,0)} \\
{} [\iota \lambda_{(1,1)}] & = & [\iota \lambda^{(1)}_{(1,1)}] \oplus [\iota \lambda^{(2)}_{(1,1)}], \label{eqn:E4(12)-lambda(1,1)} \\
{} [\iota \lambda_{(0,2)}] & = & [\iota \lambda^{(1)}_{(0,2)}] \oplus [\iota \lambda^{(2)}_{(0,2)}].  \label{eqn:E4(12)-lambda(0,2)}
\end{eqnarray}
At tier 3, $\langle \iota \lambda_{(3,0)}, \iota \lambda_{(3,0)} \rangle = \langle \iota \lambda_{(3,0)}, \iota \lambda_{(1,1)} \rangle = 2$ and similarly for $\iota \lambda_{(0,3)}$, so that $[\iota \lambda_{(3,0)}] = [\iota \lambda_{(0,3)}] = [\iota \lambda_{(1,1)}]$. From $\langle \iota \lambda_{(2,1)}, \iota \lambda_{(2,1)} \rangle = 5$, $\langle \iota \lambda_{(2,1)}, \iota \lambda_{(1,0)} \rangle = 1$ and $\langle \iota \lambda_{(2,1)}, \iota \lambda_{(0,2)} \rangle = 2$, we have two possibilities for the decomposition of $[\iota \lambda_{(2,1)}]$:
\begin{equation} \label{eqn:E4(12)-lambda(2,1)}
[\iota \lambda_{(2,1)}] = \left\{ \begin{array}{ll}
{} [\iota \lambda_{(1,0)}] \oplus 2 [\iota \lambda_{(0,2)}^{(j)}] & \textrm{case } (i), \\
{} [\iota \lambda_{(1,0)}] \oplus [\iota \lambda_{(0,2)}^{(1)}] \oplus [\iota \lambda_{(0,2)}^{(2)}] \oplus [\iota \lambda_{(2,1)}^{(1)}] \oplus [\iota \lambda_{(2,1)}^{(2)}] & \textrm{case } (ii),
\end{array} \right.
\end{equation}
where we may assume $j=1$ without loss of generality. Similarly,
\begin{equation} \label{eqn:E4(12)-lambda(1,2)}
[\iota \lambda_{(1,2)}] = \left\{ \begin{array}{ll}
{} [\iota \lambda_{(0,1)}] \oplus 2 [\iota \lambda_{(2,0)}^{(1)}] & \textrm{case } (i'), \\
{} [\iota \lambda_{(0,1)}] \oplus [\iota \lambda_{(2,0)}^{(1)}] \oplus [\iota \lambda_{(2,0)}^{(2)}] \oplus [\iota \lambda_{(1,2)}^{(1)}] \oplus [\iota \lambda_{(1,2)}^{(2)}] & \textrm{case } (ii'),
\end{array} \right.
\end{equation}
At tier 4 we have $\langle \iota \lambda_{(4,0)}, \iota \lambda_{(4,0)} \rangle = 3$, $\langle \iota \lambda_{(4,0)}, \iota \lambda_{(1,0)} \rangle = 1$ and $\langle \iota \lambda_{(4,0)}, \iota \lambda_{(0,2)} \rangle = 2$, and similarly for $\iota \lambda_{(0,4)}$, giving
\begin{eqnarray*}
{} [\iota \lambda_{(4,0)}] & = & [\iota \lambda_{(1,0)}] \oplus [\iota \lambda^{(1)}_{(0,2)}] \oplus [\iota \lambda^{(2)}_{(0,2)}], \\
{} [\iota \lambda_{(0,4)}] & = & [\iota \lambda_{(0,1)}] \oplus [\iota \lambda^{(1)}_{(2,0)}] \oplus [\iota \lambda^{(2)}_{(2,0)}].
\end{eqnarray*}
From $\langle \iota \lambda_{(3,1)}, \iota \lambda_{(3,1)} \rangle = 8$, $\langle \iota \lambda_{(3,1)}, \iota \lambda_{(0,1)} \rangle = 2$, $\langle \iota \lambda_{(3,1)}, \iota \lambda_{(2,0)} \rangle = 2$ and $\langle \iota \lambda_{(3,1)}, \iota \lambda_{(1,2)} \rangle = 6$ we have
\begin{eqnarray*}
{} [\iota \lambda_{(3,1)}] & = & \left\{ \begin{array}{ll}
{} 2 [\iota \lambda_{(0,1)}] \oplus 2 [\iota \lambda_{(2,0)}^{(1)}] & \textrm{for case } (i'), \\
{} 2 [\iota \lambda_{(0,1)}] \oplus [\iota \lambda_{(2,0)}^{(1)}] \oplus [\iota \lambda_{(2,0)}^{(2)}] \oplus [\iota \lambda_{(1,2)}^{(1)}] \oplus [\iota \lambda_{(1,2)}^{(2)}] & \textrm{for case } (ii'),
\end{array} \right. \\
{} [\iota \lambda_{(1,3)}] & = & \left\{ \begin{array}{ll}
{} 2 [\iota \lambda_{(1,0)}] \oplus 2 [\iota \lambda_{(0,2)}^{(1)}] & \textrm{for case } (i), \\
{} 2 [\iota \lambda_{(1,0)}] \oplus [\iota \lambda_{(0,2)}^{(1)}] \oplus [\iota \lambda_{(0,2)}^{(2)}] \oplus [\iota \lambda_{(2,1)}^{(1)}] \oplus [\iota \lambda_{(2,1)}^{(2)}] & \textrm{for case } (ii).
\end{array} \right. \\
\end{eqnarray*}
We have $\langle \iota \lambda_{(2,2)}, \iota \lambda_{(2,2)} \rangle = 11$, $\langle \iota \lambda_{(2,2)}, \iota \lambda_{(0,0)} \rangle = 1$ and $\langle \iota \lambda_{(2,2)}, \iota \lambda_{(1,1)} \rangle = 4$, giving
\begin{equation} \label{eqn:E4(12)-lambda(2,2)}
[\iota \lambda_{(2,2)}] = \left\{ \begin{array}{ll}
{} [\iota \lambda_{(0,0)}] {} \oplus 3 [\iota \lambda_{(1,1)}^{(j)}] \oplus [\iota \lambda_{(1,1)}^{(3-j)}] & \textrm{case I}, \\
{} [\iota \lambda_{(0,0)}] {} \oplus 2 [\iota \lambda_{(1,1)}^{(1)}] \oplus 2 [\iota \lambda_{(1,1)}^{(2)}] \oplus [\iota \lambda_{(2,2)}^{(1)}] \oplus [\iota \lambda_{(2,2)}^{(2)}] & \textrm{case II},
\end{array} \right.
\end{equation}
where $j \in \{ 1,2 \}$. Again, without loss of generality, we may assume that $j=1$, and we see that for case I nothing new appears at tier 4. For case II, at tier 5 we find $[\iota \lambda_{(5,0)}] = [\iota \lambda_{(0,4)}]$, $[\iota \lambda_{(0,5)}] = [\iota \lambda_{(4,0)}]$, $[\iota \lambda_{(4,1)}] = [\iota \lambda_{(1,4)}] = [\iota \lambda_{(0,0)}] \oplus 2 [\iota \lambda_{(1,1)}^{(1)}] \oplus 2 [\iota \lambda_{(1,1)}^{(2)}]$ and
\begin{eqnarray*}
{} [\iota \lambda_{(3,2)}] & = & \left\{ \begin{array}{ll}
{} 2 [\iota \lambda_{(1,0)}] \oplus 3 [\iota \lambda_{(0,2)}^{(1)}] \oplus [\iota \lambda_{(0,2)}^{(2)}] & \textrm{for case } (i), \\
{} 2 [\iota \lambda_{(1,0)}] \oplus 2 [\iota \lambda_{(0,2)}^{(1)}] \oplus 2 [\iota \lambda_{(0,2)}^{(2)}] \oplus [\iota \lambda_{(2,1)}^{(1)}] \oplus [\iota \lambda_{(2,1)}^{(2)}] & \textrm{for case } (ii),
\end{array} \right. \\
{} [\iota \lambda_{(2,3)}] & = & \left\{ \begin{array}{ll}
{} 2 [\iota \lambda_{(0,1)}] \oplus 3 [\iota \lambda_{(2,0)}^{(1)}] \oplus [\iota \lambda_{(2,0)}^{(2)}] & \textrm{for case } (i'), \\
{} 2 [\iota \lambda_{(0,1)}] \oplus 2 [\iota \lambda_{(2,0)}^{(1)}] \oplus 2 [\iota \lambda_{(2,0)}^{(2)}] \oplus [\iota \lambda_{(1,2)}^{(1)}] \oplus [\iota \lambda_{(1,2)}^{(2)}] & \textrm{for case } (ii'),
\end{array} \right. \\
\end{eqnarray*}
and nothing new appears at tier 5.
Then the total number of irreducible sectors for case I$(i)(i')$ is 9, for cases I$(i)(ii')$, I$(ii)(i')$, II$(i)(i')$ we have 11, for cases I$(ii)(ii')$, II$(i)(ii')$, II$(ii)(i')$ we have 13 and for case II$(ii)(ii')$ we have 15. The values of $\mathrm{tr}(Z)$ at level 12 are $\mathrm{tr}(Z_{\mathcal{A}^{(12)}}) = 55$, $\mathrm{tr}(Z_{\mathcal{D}^{(12)}}) = 19$, $\mathrm{tr}(Z_{\mathcal{A}^{(12)\ast}}) = 5$, $\mathrm{tr}(Z_{\mathcal{D}^{(12)\ast}}) = 15$, $\mathrm{tr}(Z_{\mathcal{E}^{(12)}}) = 12$, $\mathrm{tr}(Z_{\mathcal{E}_{MS}^{(12)\ast}}) = 11$ and $\mathrm{tr}(Z_{\mathcal{E}_{MS}^{(12)}}) = 17$. So we see that the only possible cases are I$(i)(ii')$, I$(ii)(i')$, II$(i)(i')$ which have corresponding modular invariant $Z_{\mathcal{E}_{MS}^{(12)\ast}}$, and II$(ii)(ii')$ associated with the modular invariant $Z_{\mathcal{D}^{(12)\ast}}$. For case II$(i)(i')$, where we again use the notation $\rho = \lambda_{(1,0)}$, we have $[\iota \lambda_{(1,2)}] [\rho] = [\iota \lambda_{(1,1)}] \oplus [\iota \lambda_{(0,3)}] \oplus [\iota \lambda_{(2,2)}]$ and $[\iota \lambda_{(1,2)}] [\rho] = ([\iota \lambda_{(0,1)}] \oplus 2 [\iota \lambda_{(2,0)}^{(1)}]) [\rho] = [\iota \lambda_{(0,0)}] \oplus [\iota \lambda_{(1,1)}] \oplus 2 ([\iota \lambda_{(2,0)}^{(1)}] [\rho])$, giving $2 [\iota \lambda_{(2,0)}^{(1)}] [\rho] = 3 [\iota \lambda_{(1,1)}^{(1)}] \oplus 3 [\iota \lambda_{(1,1)}^{(2)}] \oplus [\iota \lambda_{(2,2)}^{(1)}] \oplus [\iota \lambda_{(2,2)}^{(2)}]$, which is impossible since $[\iota \lambda_{(2,0)}^{(1)}] [\rho]$ must have integer coefficients. Note that case II$(ii)(i')$ is the conjugate of case II$(i)(ii')$, where we replace $\iota \lambda_{(\mu_1, \mu_2)} \leftrightarrow \iota \lambda_{(\mu_2, \mu_1)}$. So we need to only consider cases I$(i)(ii')$ and II$(ii)(ii')$.

Consider first the case I$(i)(ii')$. From $[\iota \lambda_{(2,1)}] [\rho] = [\iota \lambda_{(2,0)}] \oplus [\iota \lambda_{(1,2)}] \oplus [\iota \lambda_{(3,1)}]$ and (\ref{eqn:E4(12)-lambda(2,1)}) we find $[\iota \lambda_{(2,1)}^{(1)}] [\rho] = [\iota \lambda_{(2,0)}] \oplus [\iota \lambda_{(1,2)}] \oplus [\iota \lambda_{(3,1)}] \ominus ([\iota \lambda_{(0,1)}] \oplus [\iota \lambda_{(2,0)}]) = [\iota \lambda_{(0,1)}] \oplus [\iota \lambda_{(2,1)}^{(1)}] \oplus [\iota \lambda_{(2,1)}^{(2)}] \oplus [\iota \lambda_{(1,2)}^{(1)}] \oplus [\iota \lambda_{(1,2)}^{(2)}]$. Then by $[\iota \lambda_{(0,2)}] [\rho] = [\iota \lambda_{(0,1)}] \oplus [\iota \lambda_{(1,2)}]$ and (\ref{eqn:E4(12)-lambda(0,2)}), $[\iota \lambda_{(0,2)}^{(2)}] [\rho] = [\iota \lambda_{(0,1)}]$. From $[\iota \lambda_{(1,1)}] [\rho] = [\iota \lambda_{(1,0)}] \oplus [\iota \lambda_{(0,2)}] \oplus [\iota \lambda_{(2,1)}]$ and (\ref{eqn:E4(12)-lambda(1,1)}) we obtain
\begin{equation} \label{eqn:E4(12)-rho_x_lambda(1,1)}
([\iota \lambda_{(1,1)}^{(1)}] \oplus [\iota \lambda_{(1,1)}^{(2)}]) [\rho] = 2 [\iota \lambda_{(1,0)}] \oplus 3 [\iota \lambda_{(0,2)}^{(1)}] \oplus [\iota \lambda_{(0,2)}^{(2)}],
\end{equation}
whilst from $[\iota \lambda_{(2,2)}] [\rho] = [\iota \lambda_{(2,1)}] \oplus [\iota \lambda_{(1,3)}] \oplus [\iota \lambda_{(3,2)}]$ and (\ref{eqn:E4(12)-lambda(2,2)}) we have
\begin{equation} \label{eqn:E4(12)-rho_x_lambda(2,2)}
(3 [\iota \lambda_{(1,1)}^{(1)}] \oplus [\iota \lambda_{(1,1)}^{(2)}]) [\rho] = 4 [\iota \lambda_{(1,0)}] \oplus 7 [\iota \lambda_{(0,2)}^{(1)}] \oplus [\iota \lambda_{(0,2)}^{(2)}].
\end{equation}
Then from (\ref{eqn:E4(12)-rho_x_lambda(1,1)}) and (\ref{eqn:E4(12)-rho_x_lambda(2,2)}) we find
$$[\iota \lambda_{(1,1)}^{(1)}] [\rho] = [\iota \lambda_{(1,0)}] \oplus 2 [\iota \lambda_{(0,2)}^{(1)}], \qquad [\iota \lambda_{(1,1)}^{(2)}] [\rho] = [\iota \lambda_{(1,0)}] \oplus [\iota \lambda_{(0,2)}^{(1)}] \oplus [\iota \lambda_{(0,2)}^{(2)}].$$
In the same manner, by considering $[\iota \lambda_{(2,0)}] [\rho] = [\iota \lambda_{(1,1)}] \oplus [\iota \lambda_{(3,0)}]$ and $[\iota \lambda_{(1,2)}] [\rho] = [\iota \lambda_{(1,1)}] \oplus [\iota \lambda_{(0,3)}] \oplus [\iota \lambda_{(2,2)}]$, and using (\ref{eqn:E4(12)-lambda(2,0)}) and (\ref{eqn:E4(12)-lambda(1,2)}), we have
\begin{eqnarray}
{} ([\iota \lambda_{(2,0)}^{(1)}] \oplus [\iota \lambda_{(2,0)}^{(2)}]) [\rho] & = & 2 [\iota \lambda_{(1,1)}^{(1)}] \oplus 2 [\iota \lambda_{(1,1)}^{(2)}], \label{eqn:E4(12)-rho_x_lambda(2,0)} \\
{} ([\iota \lambda_{(2,0)}^{(1)}] \oplus [\iota \lambda_{(2,0)}^{(2)}] \oplus [\iota \lambda_{(1,2)}^{(1)}] \oplus [\iota \lambda_{(1,2)}^{(2)}]) [\rho] & = & [\iota \lambda_{(0,0)}] \oplus 5 [\iota \lambda_{(1,1)}^{(1)}] \oplus 3 [\iota \lambda_{(1,1)}^{(2)}] \nonumber \\
& & \quad \ominus ([\iota \lambda_{(0,0)}] \oplus [\iota \lambda_{(1,1)}]). \label{eqn:E4(12)-rho_x_lambda(1,2)}
\end{eqnarray}
Then from (\ref{eqn:E4(12)-rho_x_lambda(2,0)}), (\ref{eqn:E4(12)-rho_x_lambda(1,1)}) and (\ref{eqn:E4(12)-lambda(1,1)}), we have $([\iota \lambda_{(1,2)}^{(1)}] [\rho]) \oplus ([\iota \lambda_{(1,2)}^{(2)}] [\rho]) = 2 [\iota \lambda_{(1,1)}^{(1)}]$ giving $[\iota \lambda_{(1,2)}^{(j)}] [\rho] = [\iota \lambda_{(1,1)}^{(1)}]$ for $j=1,2$.
From $[\iota \lambda_{(2,2)}] [\overline{\rho}] = [\iota \lambda_{(1,2)}] \oplus [\iota \lambda_{(3,1)}] \oplus [\iota \lambda_{(2,3)}]$ and (\ref{eqn:E4(12)-lambda(2,2)}) we have
\begin{equation} \label{eqn:E4(12)-rho(bar)_x_lambda(2,2)}
([\iota \lambda_{(1,1)}] \oplus 2 [\iota \lambda_{(1,1)}^{(1)}]) [\overline{\rho}] = 4 [\iota \lambda_{(0,1)}] \oplus 4 [\iota \lambda_{(2,0)}^{(1)}] \oplus 4 [\iota \lambda_{(2,0)}^{(2)}] \oplus 3 [\iota \lambda_{(1,2)}^{(1)}] \oplus 3 [\iota \lambda_{(1,2)}^{(2)}],
\end{equation}
giving $2 [\iota \lambda_{(1,1)}^{(1)}] [\overline{\rho}] = 2 [\iota \lambda_{(0,1)}] \oplus 2 [\iota \lambda_{(2,0)}^{(1)}] \oplus 2 [\iota \lambda_{(2,0)}^{(2)}] \oplus 2 [\iota \lambda_{(1,2)}^{(1)}] \oplus 2 [\iota \lambda_{(1,2)}^{(2)}]$. Then $\langle \iota \lambda_{(2,0)}^{(j)} \circ \rho, \lambda_{(1,1)}^{(1)} \rangle = \langle \iota \lambda_{(2,0)}^{(j)}, \lambda_{(1,1)}^{(1)} \circ \overline{\rho} \rangle = 1$ for $j=1,2$, and the decompositions of $[\iota \lambda_{(2,0)}^{(1)}] [\rho]$ and $[\iota \lambda_{(2,0)}^{(2)}] [\rho]$ both contain the irreducible sector $[\iota \lambda_{(1,1)}^{(1)}]$. Then $[\iota \lambda_{(1,1)}^{(2)}] [\overline{\rho}] = ([\iota \lambda_{(1,1)}] [\overline{\rho}]) \ominus ([\iota \lambda_{(1,1)}^{(1)}] [\overline{\rho}]) = [\iota \lambda_{(0,1)}] \oplus [\iota \lambda_{(2,0)}^{(1)}] \oplus [\iota \lambda_{(2,0)}^{(2)}]$ and $[\iota \lambda_{(2,0)}^{(1)}] [\rho]$ and $[\iota \lambda_{(2,0)}^{(2)}] [\rho]$ both also contain $[\iota \lambda_{(1,1)}^{(2)}]$. Then from (\ref{eqn:E4(12)-rho_x_lambda(2,0)}) we have $[\iota \lambda_{(2,0)}^{(j)}] [\rho] = [\iota \lambda_{(1,1)}^{(1)}] \oplus [\iota \lambda_{(1,1)}^{(2)}]$. The nimrep graph for multiplication by $[\rho]$ for the case I$(i)(ii')$ is then seen to be just the graph $\mathcal{E}_4^{(12)}$.

\begin{figure}[tb]
\begin{center}
  \includegraphics[width=75mm]{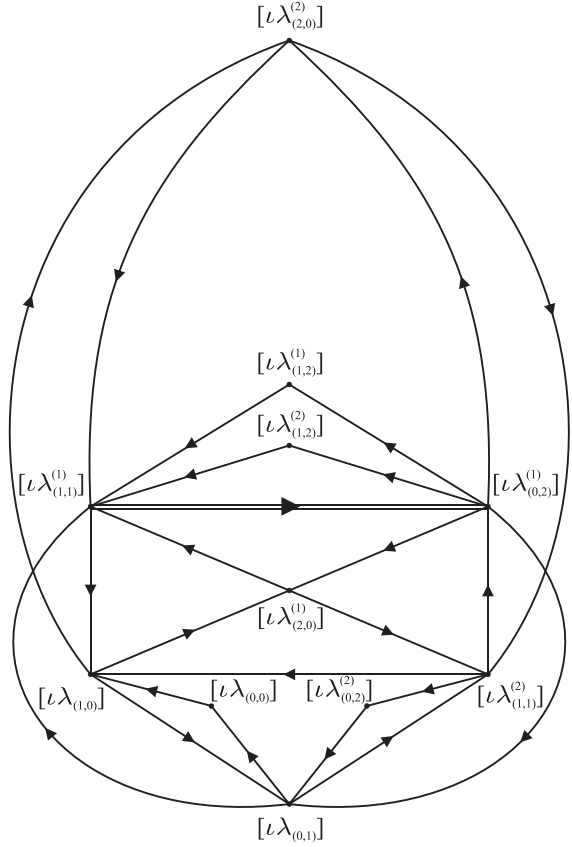}
 \caption{$M$-$N$ graph for the $\mathcal{E}_4^{(12)}$ $SU(3)$-GHJ subfactor}\label{fig:fig-GHJ-10}
\end{center}
\end{figure}

Now consider the case II$(ii)(ii')$, which has corresponding modular invariant $Z_{\mathcal{D}^{(12)\ast}}$. We obtain the following sector products:
\begin{eqnarray*}
{} ([\iota \lambda_{(2,0)}^{(1)}] \oplus [\iota \lambda_{(2,0)}^{(2)}]) [\rho] & = & 2 [\iota \lambda_{(1,1)}^{(1)}] \oplus 2 [\iota \lambda_{(1,1)}^{(2)}], \\
{} ([\iota \lambda_{(1,1)}^{(1)}] \oplus [\iota \lambda_{(1,1)}^{(2)}]) [\rho] & = & 2 [\iota \lambda_{(1,0)}] \oplus 2 [\iota \lambda_{(0,2)}^{(1)}] \oplus 2 [\iota \lambda_{(0,2)}^{(2)}] \oplus [\iota \lambda_{(2,1)}^{(1)}] \oplus [\iota \lambda_{(2,1)}^{(2)}], \\
{} ([\iota \lambda_{(0,2)}^{(1)}] \oplus [\iota \lambda_{(0,2)}^{(2)}]) [\rho] & = & 2 [\iota \lambda_{(0,1)}] \oplus [\iota \lambda_{(2,0)}^{(1)}] \oplus [\iota \lambda_{(2,0)}^{(2)}] \oplus [\iota \lambda_{(1,2)}^{(1)}] \oplus [\iota \lambda_{(1,2)}^{(2)}], \\
{} ([\iota \lambda_{(2,1)}^{(1)}] \oplus [\iota \lambda_{(2,1)}^{(2)}]) [\rho] & = & [\iota \lambda_{(2,0)}^{(1)}] \oplus [\iota \lambda_{(2,0)}^{(2)}] \oplus [\iota \lambda_{(1,2)}^{(1)}] \oplus [\iota \lambda_{(1,2)}^{(2)}], \\
{} ([\iota \lambda_{(1,2)}^{(1)}] \oplus [\iota \lambda_{(1,2)}^{(2)}]) [\rho] & = & [\iota \lambda_{(1,1)}^{(1)}] \oplus [\iota \lambda_{(1,1)}^{(2)}] \oplus [\iota \lambda_{(2,2)}^{(1)}] \oplus [\iota \lambda_{(2,2)}^{(2)}],
\end{eqnarray*}
and from $([\iota \lambda_{(2,2)}^{(1)}] \oplus [\iota \lambda_{(2,2)}^{(2)}]) [\rho] = [\iota \lambda_{(2,1)}^{(1)}] \oplus [\iota \lambda_{(2,1)}^{(2)}]$ we may choose without loss of generality $[\iota \lambda_{(2,2)}^{(j)}] [\rho] = [\iota \lambda_{(2,1)}^{(j)}]$ for $j=1,2$. Then there are four different possibilities for $[\iota \lambda_{(1,1)}^{(j)}] [\rho]$, three for $[\iota \lambda_{(2,0)}^{(j)}] [\rho]$, six for $[\iota \lambda_{(0,2)}^{(j)}] [\rho]$ and six for $[\iota \lambda_{(2,1)}^{(j)}] [\rho]$, $j=1,2$. From these, the only nimrep graph which has eigenvalues $S_{\rho \mu}/S_{0 \mu}$ with multiplicities given by the diagonal entry $Z_{\mu,\mu}$ of the modular invariant for $\mathcal{D}^{(12)\ast}$ is that for the following sector products:
\begin{eqnarray*}
{} [\iota \lambda_{(2,0)}^{(j)}] [\rho] & = & 2 [\iota \lambda_{(1,1)}^{(j)}], \\
{} [\iota \lambda_{(1,1)}^{(1)}] [\rho] & = & [\iota \lambda_{(1,0)}] \oplus 2 [\iota \lambda_{(0,2)}^{(j)}] \oplus [\iota \lambda_{(2,1)}^{(j)}], \\
{} [\iota \lambda_{(0,2)}^{(j)}] [\rho] & = & [\iota \lambda_{(0,1)}] \oplus [\iota \lambda_{(2,0)}^{(j)}] \oplus [\iota \lambda_{(1,2)}^{(j)}], \\
{} [\iota \lambda_{(2,1)}^{(j)}] [\rho] & = & [\iota \lambda_{(2,0)}^{(j)}] \oplus [\iota \lambda_{(1,2)}^{(j)}], \\
{} [\iota \lambda_{(1,2)}^{(j)}] [\rho] & = & [\iota \lambda_{(1,1)}^{(j)}] \oplus [\iota \lambda_{(2,2)}^{(3-j)}],
\end{eqnarray*}
for $j=1,2$. For any $\lambda \in {}_M \mathcal{X}_N$, let $[\lambda] [\rho] = \bigoplus_{\mu \in {}_M \mathcal{X}_N} a_{\mu} [\mu]$, $a_{\mu} \in \mathbb{C}$. Then $\langle \mu \circ \overline{\rho}, \lambda \rangle = \langle \mu, \lambda \circ \rho \rangle = a_{\mu}$ for all $\mu \in {}_M \mathcal{X}_N$, so $[\mu] [\overline{\rho}] \supset a_{\mu} [\lambda]$. Then if $G$ is the multiplication matrix for $[\rho]$, $G^T$ is the multiplication matrix for $[\overline{\rho}]$. This graph cannot be the nimrep graph since $GG^T \neq G^TG$, which means $[\iota \lambda] [\rho] [\overline{\rho}] \neq [\iota \lambda] [\overline{\rho}] [\rho]$. Then the only possibility for the nimrep graph for the $M$-$N$ system is the graph $\mathcal{E}_4^{(12)}$, illustrated in Figure \ref{fig:fig-GHJ-10}, and the associated modular invariant is $Z_{\mathcal{E}_{MS}^{(12)\ast}}$, assuming that $\theta$ is as expressed in (\ref{theta-E4(12)}).

\subsection{$\mathcal{E}_5^{(12)}$}

For the graph $\mathcal{E}_5^{(12)}$, we have $[\theta] = [\lambda_{(0,0)}] \oplus [\lambda_{(3,3)}] \oplus [\lambda_{(9,0)}] \oplus [\lambda_{(0,9)}]$. Then computing $\langle \iota \lambda, \iota \mu \rangle = \langle \lambda, \theta \mu \rangle$ for $\lambda$, $\mu$ on $\mathcal{A}^{(12)}$, we find $\langle \iota \lambda, \iota \lambda \rangle = 1$ for $\lambda = \lambda_{(\mu_1, \mu_2)}$ such that $\mu_1 + \mu_2 \leq 2$. At tier 3 we have $\langle \iota \lambda, \iota \lambda \rangle = 2$ and $\langle \iota \lambda, \iota \mu \rangle = 0$ for $\lambda = \lambda_{(3,0)}, \lambda_{(2,1)}, \lambda_{(1,2)}, \lambda_{(0,3)}$, $\mu = \lambda_{(\mu_1, \mu_2)}$ such that $\mu_1 + \mu_2 \leq 2$. We also have $\langle \iota \lambda_{(3,0)}, \iota \lambda_{(0,3)} \rangle = 0$. Then $\lambda_{(3,0)}$, $\lambda_{(2,1)}$, $\lambda_{(1,2)}$, $\lambda_{(0,3)}$ decompose into irreducibles as
\begin{eqnarray}
{} [\iota \lambda_{(3,0)}] & = & [\iota \lambda^{(1)}_{(3,0)}] \oplus [\iota \lambda^{(2)}_{(3,0)}], \label{eqn:E5(12)-lambda(3,0)} \\
{} [\iota \lambda_{(2,1)}] & = & [\iota \lambda^{(1)}_{(2,1)}] \oplus [\iota \lambda^{(2)}_{(2,1)}], \label{eqn:E5(12)-lambda(2,1)} \\
{} [\iota \lambda_{(1,2)}] & = & [\iota \lambda^{(1)}_{(1,2)}] \oplus [\iota \lambda^{(2)}_{(1,2)}], \label{eqn:E5(12)-lambda(1,2)} \\
{} [\iota \lambda_{(0,3)}] & = & [\iota \lambda^{(1)}_{(0,3)}] \oplus [\iota \lambda^{(2)}_{(0,3)}].  \label{eqn:E5(12)-lambda(0,3)}
\end{eqnarray}
At tier 4 we have $\langle \iota \lambda_{(4,0)}, \iota \lambda_{(4,0)} \rangle = 2$, $\langle \iota \lambda_{(4,0)}, \iota \lambda_{(2,1)} \rangle = 1$ and $\langle \iota \lambda_{(4,0)}, \iota \mu \rangle = 0$ for $\mu = \lambda_{(1,0)}, \lambda_{(0,2)}$. Then $[\iota \lambda_{(4,0)}] = [\iota \lambda_{(2,1)}^{(j)}] \oplus [\iota \lambda_{(4,0)}^{(1)}]$ for $j \in \{ 1,2 \}$. We have the freedom to choose $j=1$ without loss of generality. Similarly for $\iota \lambda_{(0,4)}$. Then
\begin{eqnarray}
{} [\iota \lambda_{(4,0)}] & = & [\iota \lambda_{(2,1)}^{(1)}] \oplus [\iota \lambda_{(4,0)}^{(1)}], \label{eqn:E5(12)-lambda(4,0)} \\
{} [\iota \lambda_{(0,4)}] & = & [\iota \lambda_{(1,2)}^{(1)}] \oplus [\iota \lambda_{(0,4)}^{(1)}].  \label{eqn:E5(12)-lambda(0,4)}
\end{eqnarray}
From $\langle \iota \lambda_{(3,1)}, \iota \lambda_{(3,1)} \rangle = 3$, $\langle \iota \lambda_{(3,1)}, \iota \lambda_{(2,0)} \rangle = 1$, $\langle \iota \lambda_{(3,1)}, \iota \lambda_{(1,2)} \rangle = 1$ and $\langle \iota \lambda_{(3,1)}, \iota \lambda_{(0,4)} \rangle = 1$, we have two possibilities for the decomposition of $[\iota \lambda_{(3,1)}]$:
\begin{equation} \label{eqn:E5(12)-lambda(3,1)}
[\iota \lambda_{(3,1)}] = \left\{ \begin{array}{ll}
{} [\iota \lambda_{(2,0)}] \oplus [\iota \lambda_{(1,2)}^{(1)}] \oplus [\iota \lambda_{(3,1)}^{(1)}] & \textrm{case } (i), \\
{} [\iota \lambda_{(2,0)}] \oplus [\iota \lambda_{(1,2)}^{(2)}] \oplus [\iota \lambda_{(0,4)}^{(1)}] & \textrm{case } (ii).
\end{array} \right.
\end{equation}
Similarly,
\begin{equation} \label{eqn:E5(12)-lambda(1,3)}
[\iota \lambda_{(1,3)}] = \left\{ \begin{array}{ll}
{} [\iota \lambda_{(0,2)}] \oplus [\iota \lambda_{(2,1)}^{(1)}] \oplus [\iota \lambda_{(1,3)}^{(1)}] & \textrm{case } (i'), \\
{} [\iota \lambda_{(0,2)}] \oplus [\iota \lambda_{(2,1)}^{(2)}] \oplus [\iota \lambda_{(4,0)}^{(1)}] & \textrm{case } (ii'),
\end{array} \right.
\end{equation}
Since $\langle \iota \lambda_{(2,2)}, \iota \lambda_{(2,2)} \rangle = 3$, $\langle \iota \lambda_{(2,2)}, \iota \lambda_{(1,1)} \rangle = 1$, $\langle \iota \lambda_{(2,2)}, \iota \lambda_{(3,0)} \rangle = 1$ and $\langle \iota \lambda_{(2,2)}, \iota \lambda_{(0,3)} \rangle = 1$, we have $[\iota \lambda_{(2,2)}] = [\iota \lambda_{(1,1)}] \oplus [\iota \lambda_{(3,0)}^{(j_1)}] \oplus [\iota \lambda_{(0,3)}^{(j_2)}]$ for $j_1,j_2 \in \{ 1,2 \}$. We again have the freedom to choose, without loss of generality, $j_1 = j_2 = 1$, so that
\begin{equation} \label{eqn:E5(12)-lambda(2,2)}
[\iota \lambda_{(2,2)}] = [\iota \lambda_{(1,1)}] \oplus [\iota \lambda_{(3,0)}^{(1)}] \oplus [\iota \lambda_{(0,3)}^{(1)}].
\end{equation}
At tier 5, $\langle \iota \lambda_{(5,0)}, \iota \lambda_{(5,0)} \rangle = \langle \iota \lambda_{(5,0)}, \iota \lambda_{(0,4)} \rangle = 2$ giving $[\iota \lambda_{(5,0)}] = [\iota \lambda_{(0,4)}]$, and similarly $[\iota \lambda_{(0,5)}] = [\iota \lambda_{(4,0)}]$. Since $\langle \iota \lambda_{(3,2)}, \iota \lambda_{(3,2)} \rangle = 4$, $\langle \iota \lambda_{(3,2)}, \iota \lambda_{(1,0)} \rangle = 1$, $\langle \iota \lambda_{(3,2)}, \iota \lambda_{(0,2)} \rangle = 1$ and $\langle \iota \lambda_{(3,2)}, \iota \lambda_{(2,1)} \rangle = 2$, we have $[\iota \lambda_{(3,2)}] = [\iota \lambda_{(1,0)}] \oplus [\iota \lambda_{(0,2)}] \oplus [\iota \lambda_{(2,1)}^{(1)}] \oplus [\iota \lambda_{(2,1)}^{(2)}]$, and similarly $[\iota \lambda_{(2,3)}] = [\iota \lambda_{(0,1)}] \oplus [\iota \lambda_{(2,0)}] \oplus [\iota \lambda_{(1,2)}^{(1)}] \oplus [\iota \lambda_{(1,2)}^{(2)}]$. We have $\langle \iota \lambda_{(4,1)}, \iota \lambda_{(4,1)} \rangle = \langle \iota \lambda_{(4,1)}, \iota \lambda_{(1,4)} \rangle = \langle \iota \lambda_{(1,4)}, \iota \lambda_{(1,4)} \rangle = 3$ so that $[\iota \lambda_{(4,1)}] = [\iota \lambda_{(1,4)}]$. Since $\langle \iota \lambda_{(4,1)}, \iota \lambda_{(1,1)} \rangle = 1$, $\langle \iota \lambda_{(4,1)}, \iota \lambda_{(2,2)} \rangle = 2$, $\langle \iota \lambda_{(4,1)}, \iota \lambda_{(3,0)} \rangle = 1$ and $\langle \iota \lambda_{(4,1)}, \iota \lambda_{(0,3)} \rangle = 1$, we have two possibilities for the decomposition of $[\iota \lambda_{(4,1)}]$:
\begin{equation} \label{eqn:E5(12)-lambda(4,1)}
[\iota \lambda_{(4,1)}] = \left\{ \begin{array}{ll}
{} [\iota \lambda_{(1,1)}] {} \oplus [\iota \lambda_{(3,0)}^{(1)}] \oplus [\iota \lambda_{(0,3)}^{(2)}] & \textrm{case I}, \\
{} [\iota \lambda_{(1,1)}] {} \oplus [\iota \lambda_{(3,0)}^{(2)}] \oplus [\iota \lambda_{(0,3)}^{(1)}] & \textrm{case II}.
\end{array} \right.
\end{equation}
Then we see that no new irreducible sectors appear at tier 5. We also have at tier 6, $\langle \iota \lambda_{(5,1)}, \iota \lambda_{(5,1)} \rangle = \langle \iota \lambda_{(5,1)}, \iota \lambda_{(1,3)} \rangle = 3$ giving $[\iota \lambda_{(5,1)}] = [\iota \lambda_{(1,3)}]$. Case $(i)(i')$ gives 16 irreducible sectors, whilst case $(ii)(ii')$ gives 18 irreducibles, and therefore by looking at $\mathrm{tr}(Z)$ for the level 12 modular invariants $Z$ we see that neither of these cases is possible. Case $(ii)(i')$ is the `conjugate' of case $(i)(ii')$, that is, we replace each irreducible sector $[\iota \lambda]$ in case $(i)(ii')$ by $[\overline{\iota \lambda}]$ in case $(ii)(i')$. We therefore only need to consider case $(i)(ii')$, which has seventeen irreducible sectors: $[\lambda_{(0,0)}]$, $[\lambda_{(1,0)}]$, $[\lambda_{(0,1)}]$, $[\lambda_{(2,0)}]$, $[\lambda_{(1,1)}]$, $[\lambda_{(0,2)}]$, $[\lambda_{(3,0)}^{(1)}]$, $[\lambda_{(3,0)}^{(2)}]$, $[\lambda_{(0,3)}^{(1)}]$, $[\lambda_{(0,3)}^{(2)}]$, $[\lambda_{(2,1)}^{(1)}]$, $[\lambda_{(2,1)}^{(2)}]$, $[\lambda_{(1,2)}^{(1)}]$, $[\lambda_{(1,2)}^{(2)}]$, $[\lambda_{(4,0)}^{(1)}]$, $[\lambda_{(0,4)}^{(1)}]$ and $[\lambda_{(3,1)}^{(1)}]$.

We now consider the sector products for these irreducible sectors, where we again denote by $[\rho]$ the irreducible $N$-$N$ sector $[\lambda_{(1,0)}]$. The products $[\iota \lambda] [\rho]$ are inherited from those for the $N$-$N$ system for $\lambda = \lambda_{(\mu_1, \mu_2)}$ such that $\mu_1 + \mu_2 \leq 2$, and we use (\ref{eqn:E5(12)-lambda(3,0)})-(\ref{eqn:E5(12)-lambda(0,3)}) to decompose into irreducibles where necessary, e.g.
\begin{equation} \label{eqn:E5(12)-rho_x_lambda(0,2)}
[\iota \lambda_{(0,2)}] [\rho] = [\iota \lambda_{(0,1)}] \oplus [\iota \lambda_{(1,2)}]  = [\iota \lambda_{(0,1)}] \oplus [\iota \lambda_{(1,2)}^{(1)}] \oplus [\iota \lambda_{(1,2)}^{(2)}].
\end{equation}
From $[\iota \lambda_{(2,1)}] [\rho] = [\iota \lambda_{(2,0)}] \oplus [\iota \lambda_{(1,2)}] \oplus [\iota \lambda_{(3,1)}]$ and (\ref{eqn:E5(12)-lambda(2,1)}) we obtain
\begin{equation} \label{eqn:E5(12)-rho_x_lambda(2,1)}
([\iota \lambda_{(2,1)}^{(1)}] \oplus [\iota \lambda_{(2,1)}^{(2)}]) [\rho] = 2 [\iota \lambda_{(2,0)}] \oplus 2 [\iota \lambda_{(1,2)}^{(1)}] \oplus [\iota \lambda_{(1,2)}^{(2)}] \oplus [\iota \lambda_{(3,1)}^{(1)}].
\end{equation}
Similarly, by considering $[\iota \lambda_{(1,3)}] [\rho]$ and $[\iota \lambda_{(4,0)}] [\rho]$, and using (\ref{eqn:E5(12)-lambda(1,3)}) and (\ref{eqn:E5(12)-lambda(4,0)}) we have
\begin{eqnarray}
{} ([\iota \lambda_{(2,1)}^{(2)}] \oplus [\iota \lambda_{(4,0)}^{(1)}]) [\rho] & = & [\iota \lambda_{(2,0)}] \oplus 2 [\iota \lambda_{(1,2)}^{(1)}] \oplus [\iota \lambda_{(1,2)}^{(2)}] \oplus [\iota \lambda_{(0,4)}^{(1)}], \label{eqn:E5(12)-rho_x_lambda(1,3)} \\
{} ([\iota \lambda_{(2,1)}^{(1)}] \oplus [\iota \lambda_{(4,0)}^{(1)}]) [\rho] & = & [\iota \lambda_{(2,0)}] \oplus 2 [\iota \lambda_{(1,2)}^{(1)}] \oplus [\iota \lambda_{(3,1)}^{(1)}] \oplus [\iota \lambda_{(0,4)}^{(1)}]. \label{eqn:E5(12)-rho_x_lambda(4,0)}
\end{eqnarray}
Then from (\ref{eqn:E5(12)-rho_x_lambda(2,1)})-(\ref{eqn:E5(12)-rho_x_lambda(4,0)}) we find
\begin{eqnarray}
{} [\iota \lambda_{(2,1)}^{(1)}] [\rho] & = & [\iota \lambda_{(2,0)}] \oplus [\iota \lambda_{(1,2)}^{(1)}] \oplus [\iota \lambda_{(3,1)}^{(1)}], \label{eqn:E5(12)-rho_x_lambda(2,1)(1)} \\
{} [\iota \lambda_{(2,1)}^{(2)}] [\rho] & = & [\iota \lambda_{(2,0)}] \oplus [\iota \lambda_{(1,2)}^{(1)}] \oplus [\iota \lambda_{(1,2)}^{(2)}], \label{eqn:E5(12)-rho_x_lambda(2,1)(2)} \\
{} [\iota \lambda_{(4,0)}^{(1)}] [\rho] & = & [\iota \lambda_{(1,2)}^{(1)}] \oplus [\iota \lambda_{(0,4)}^{(1)}]. \label{eqn:E5(12)-rho_x_lambda(4,0)(1)}
\end{eqnarray}

Now we focus on case I. From $[\iota \lambda_{(3,0)}] [\rho] = [\iota \lambda_{(4,0)}] \oplus [\iota \lambda_{(2,1)}]$ and (\ref{eqn:E5(12)-lambda(3,0)}) we obtain
\begin{equation} \label{eqn:E5(12)-rho_x_lambda(3,0)}
([\iota \lambda_{(3,0)}^{(1)}] \oplus [\iota \lambda_{(3,0)}^{(2)}]) [\rho] = 2 [\iota \lambda_{(2,1)}^{(1)}] \oplus [\iota \lambda_{(2,1)}^{(2)}] \oplus [\iota \lambda_{(4,0)}^{(1)}].
\end{equation}
Similarly by considering $[\iota \lambda_{(0,3)}] [\rho]$ we have
\begin{equation} \label{eqn:E5(12)-rho_x_lambda(0,3)}
([\iota \lambda_{(0,3)}^{(1)}] \oplus [\iota \lambda_{(0,3)}^{(2)}]) [\rho] = 2 [\iota \lambda_{(0,2)}] \oplus [\iota \lambda_{(2,1)}^{(2)}] \oplus [\iota \lambda_{(4,0)}^{(1)}].
\end{equation}
From $[\iota \lambda_{(2,2)}] [\rho] = [\iota \lambda_{(2,1)}] \oplus [\iota \lambda_{(1,3)}] \oplus [\iota \lambda_{(3,2)}]$ and (\ref{eqn:E5(12)-lambda(2,2)}) we find
\begin{equation} \label{eqn:E5(12)-rho_x_lambda(2,2)}
([\iota \lambda_{(3,0)}^{(1)}] \oplus [\iota \lambda_{(0,3)}^{(1)}]) [\rho] = [\iota \lambda_{(0,2)}] \oplus [\iota \lambda_{(2,1)}^{(1)}] \oplus 2 [\iota \lambda_{(2,1)}^{(2)}] \oplus [\iota \lambda_{(4,0)}^{(1)}],
\end{equation}
whilst from $[\iota \lambda_{(4,1)}] [\rho] = [\iota \lambda_{(4,0)}] \oplus [\iota \lambda_{(3,2)}] \oplus [\iota \lambda_{(5,1)}]$ and (\ref{eqn:E5(12)-lambda(4,1)}) we find
\begin{equation} \label{eqn:E5(12)-rho_x_lambda(4,1)}
([\iota \lambda_{(3,0)}^{(1)}] \oplus [\iota \lambda_{(0,3)}^{(2)}]) [\rho] = [\iota \lambda_{(1,0)}] \oplus 2 [\iota \lambda_{(0,2)}] \oplus 2 [\iota \lambda_{(2,1)}^{(1)}] \oplus 2 [\iota \lambda_{(2,1)}^{(2)}] \oplus 2 [\iota \lambda_{(4,0)}^{(1)}].
\end{equation}
Then from (\ref{eqn:E5(12)-rho_x_lambda(3,0)})-(\ref{eqn:E5(12)-rho_x_lambda(4,1)}) we obtain
\begin{eqnarray}
{} [\iota \lambda_{(3,0)}^{(1)}] [\rho] & = & [\iota \lambda_{(2,1)}^{(1)}] \oplus [\iota \lambda_{(2,1)}^{(2)}] \oplus [\iota \lambda_{(4,0)}^{(1)}], \label{eqn:E5(12)-rho_x_lambda(3,0)(1)} \\
{} [\iota \lambda_{(3,0)}^{(2)}] [\rho] & = & [\iota \lambda_{(2,1)}^{(1)}], \label{eqn:E5(12)-rho_x_lambda(3,0)(2)} \\
{} [\iota \lambda_{(0,3)}^{(1)}] [\rho] & = & [\iota \lambda_{(0,2)}] \oplus [\iota \lambda_{(2,1)}^{(2)}], \label{eqn:E5(12)-rho_x_lambda(0,3)(1)} \\
{} [\iota \lambda_{(0,3)}^{(2)}] [\rho] & = & [\iota \lambda_{(0,2)}] \oplus [\iota \lambda_{(4,0)}^{(1)}]. \label{eqn:E5(12)-rho_x_lambda(0,3)(2)}
\end{eqnarray}
Next, by considering $[\iota \lambda] [\rho]$ for $\lambda = \lambda_{(1,2)}, \lambda_{(3,1)}, \lambda_{(0,4)}$, and (\ref{eqn:E5(12)-lambda(1,2)}), (\ref{eqn:E5(12)-lambda(3,1)}) and (\ref{eqn:E5(12)-lambda(0,4)}) we obtain
\begin{eqnarray}
{} ([\iota \lambda_{(1,2)}^{(1)}] \oplus [\iota \lambda_{(1,2)}^{(2)}]) [\rho] & = & 2 [\iota \lambda_{(1,1)}] \oplus [\iota \lambda_{(3,0)}^{(1)}] \oplus 2 [\iota \lambda_{(0,3)}^{(1)}] \oplus [\iota \lambda_{(0,3)}^{(2)}], \label{eqn:E5(12)-rho_x_lambda(1,2)} \\
{} ([\iota \lambda_{(1,2)}^{(1)}] \oplus [\iota \lambda_{(3,1)}^{(1)}]) [\rho] & = & [\iota \lambda_{(1,1)}] \oplus 2 [\iota \lambda_{(3,0)}^{(1)}] \oplus [\iota \lambda_{(0,3)}^{(1)}] \oplus [\iota \lambda_{(0,3)}^{(2)}], \label{eqn:E5(12)-rho_x_lambda(3,1)} \\
{} ([\iota \lambda_{(1,2)}^{(1)}] \oplus [\iota \lambda_{(0,4)}^{(1)}]) [\rho] & = & [\iota \lambda_{(1,1)}] \oplus [\iota \lambda_{(3,0)}^{(1)}] \oplus [\iota \lambda_{(0,3)}^{(1)}] \oplus 2 [\iota \lambda_{(0,3)}^{(2)}]. \label{eqn:E5(12)-rho_x_lambda(0,4)}
\end{eqnarray}
We see from (\ref{eqn:E5(12)-rho_x_lambda(1,2)})-(\ref{eqn:E5(12)-rho_x_lambda(0,4)}) that $[\iota \lambda_{(1,2)}^{(1)}] [\rho] \subset [\iota \lambda_{(1,1)}] \oplus [\iota \lambda_{(3,0)}^{(1)}] \oplus [\iota \lambda_{(0,3)}^{(1)}] \oplus [\iota \lambda_{(0,3)}^{(2)}]$. From (\ref{eqn:E5(12)-rho_x_lambda(0,2)}) and (\ref{eqn:E5(12)-rho_x_lambda(2,1)(1)})-(\ref{eqn:E5(12)-rho_x_lambda(4,0)(1)}) we see that $[\iota \lambda_{(1,2)}^{(1)}] [\overline{\rho}] = [\iota \lambda_{(0,2)}] \oplus [\iota \lambda_{(2,1)}^{(1)}] \oplus [\iota \lambda_{(2,1)}^{(2)}] \oplus [\iota \lambda_{(4,0)}^{(1)}]$, since $\langle \iota \lambda_{(1,2)}^{(1)} \circ \overline{\rho}, \iota \lambda \rangle = \langle \iota \lambda_{(1,2)}^{(1)}, \iota \lambda \circ \rho \rangle = 1$ for $\lambda = \lambda_{(0,2)}, \lambda_{(2,1)}^{(1)}, \lambda_{(2,1)}^{(2)}, \lambda_{(4,0)}^{(1)}$. Then $\langle \iota \lambda_{(1,2)}^{(1)} \circ \rho, \iota \lambda_{(1,2)}^{(1)} \circ \rho \rangle = \langle \iota \lambda_{(1,2)}^{(1)} \circ \overline{\rho}, \iota \lambda_{(1,2)}^{(1)} \circ \overline{\rho} \rangle = 4$ implies that we must have $[\iota \lambda_{(1,2)}^{(1)}] [\rho] = [\iota \lambda_{(1,1)}] \oplus [\iota \lambda_{(3,0)}^{(1)}] \oplus [\iota \lambda_{(0,3)}^{(1)}] \oplus [\iota \lambda_{(0,3)}^{(2)}]$. Then from (\ref{eqn:E5(12)-rho_x_lambda(1,2)})-(\ref{eqn:E5(12)-rho_x_lambda(0,4)}) we obtain
\begin{eqnarray*}
{} [\iota \lambda_{(1,2)}^{(2)}] [\rho] & = & [\iota \lambda_{(1,1)}] \oplus [\iota \lambda_{(0,3)}^{(1)}], \label{eqn:E5(12)-rho_x_lambda(1,2)(2)} \\
{} [\iota \lambda_{(3,0)}^{(2)}] [\rho] & = & [\iota \lambda_{(0,3)}^{(2)}], \label{eqn:E5(12)-rho_x_lambda(0,4)(1)} \\
{} [\iota \lambda_{(0,3)}^{(1)}] [\rho] & = & [\iota \lambda_{(3,0)}^{(1)}]. \label{eqn:E5(12)-rho_x_lambda(3,1)(1)}
\end{eqnarray*}
It is easy to check that the nimrep graph for multiplication by $[\rho]$ obtained in case I is just the graph $\mathcal{E}_5^{(12)}$.

For case II, we again have (\ref{eqn:E5(12)-rho_x_lambda(1,2)}), and by considering $[\iota \lambda_{(3,1)}] [\rho] = [\iota \lambda_{(3,0)}] \oplus [\iota \lambda_{(2,2)}] \oplus [\iota \lambda_{(4,1)}]$ and (\ref{eqn:E5(12)-lambda(3,0)}), (\ref{eqn:E5(12)-lambda(2,2)}) and (\ref{eqn:E5(12)-lambda(4,1)}) we obtain
\begin{equation} \label{eqn:E5(12)-rho_x_lambda(3,1)II}
([\iota \lambda_{(1,2)}^{(1)}] \oplus [\iota \lambda_{(3,1)}^{(1)}]) [\rho] = [\iota \lambda_{(1,1)}] \oplus [\iota \lambda_{(3,0)}^{(1)}] \oplus 2 [\iota \lambda_{(3,0)}^{(2)}] \oplus 2 [\iota \lambda_{(0,3)}^{(1)}],
\end{equation}
and similarly from $[\iota \lambda_{(0,4)}] [\rho]$, (\ref{eqn:E5(12)-lambda(0,4)}), (\ref{eqn:E5(12)-lambda(0,3)}) and (\ref{eqn:E5(12)-lambda(4,1)}) we obtain
\begin{equation} \label{eqn:E5(12)-rho_x_lambda(0,4)II}
{} ([\iota \lambda_{(1,2)}^{(1)}] \oplus [\iota \lambda_{(0,4)}^{(1)}]) [\rho] = [\iota \lambda_{(1,1)}] \oplus [\iota \lambda_{(3,0)}^{(2)}] \oplus 2 [\iota \lambda_{(0,3)}^{(1)}] \oplus [\iota \lambda_{(0,3)}^{(2)}].
\end{equation}
Then from (\ref{eqn:E5(12)-rho_x_lambda(1,2)}), (\ref{eqn:E5(12)-rho_x_lambda(3,1)II}) and (\ref{eqn:E5(12)-rho_x_lambda(0,4)II}) we see that $[\iota \lambda_{(1,2)}^{(1)}] [\rho] \subset [\iota \lambda_{(1,1)}] \oplus 2 [\iota \lambda_{(0,3)}^{(1)}]$. Since $\langle \iota \lambda_{(1,2)}^{(1)} \circ \rho, \iota \lambda_{(1,2)}^{(1)} \circ \rho \rangle = \langle \iota \lambda_{(1,2)}^{(1)} \circ \overline{\rho}, \iota \lambda_{(1,2)}^{(1)} \circ \overline{\rho} \rangle = 4$, we must have $[\iota \lambda_{(1,2)}^{(1)}] [\rho] = 2 [\iota \lambda_{(0,3)}^{(1)}]$. Then from (\ref{eqn:E5(12)-rho_x_lambda(1,2)}) we obtain $[\iota \lambda_{(1,2)}^{(2)}] [\rho] = 2 [\iota \lambda_{(1,1)}] \oplus [\iota \lambda_{(3,0)}^{(1)}] \oplus [\iota \lambda_{(0,3)}^{(2)}]$, and we have $\langle \iota \lambda_{(1,2)}^{(2)} \circ \overline{\rho}, \iota \lambda_{(1,2)}^{(2)} \circ \overline{\rho} \rangle = \langle \iota \lambda_{(1,2)}^{(1)} \circ \rho, \iota \lambda_{(1,2)}^{(1)} \circ \rho \rangle = 6$. From (\ref{eqn:E5(12)-rho_x_lambda(0,2)}) and (\ref{eqn:E5(12)-rho_x_lambda(2,1)(1)})-(\ref{eqn:E5(12)-rho_x_lambda(4,0)(1)}) we see that $[\iota \lambda_{(1,2)}^{(2)}] [\overline{\rho}] = [\iota \lambda_{(0,2)}] \oplus [\iota \lambda_{(2,1)}^{(2)}]$, giving $\langle \iota \lambda_{(1,2)}^{(2)} \circ \overline{\rho}, \iota \lambda_{(1,2)}^{(2)} \circ \overline{\rho} \rangle =2 \neq 6$, which is a contradiction. Then we reject case II.

\begin{figure}[tb]
\begin{center}
  \includegraphics[width=110mm]{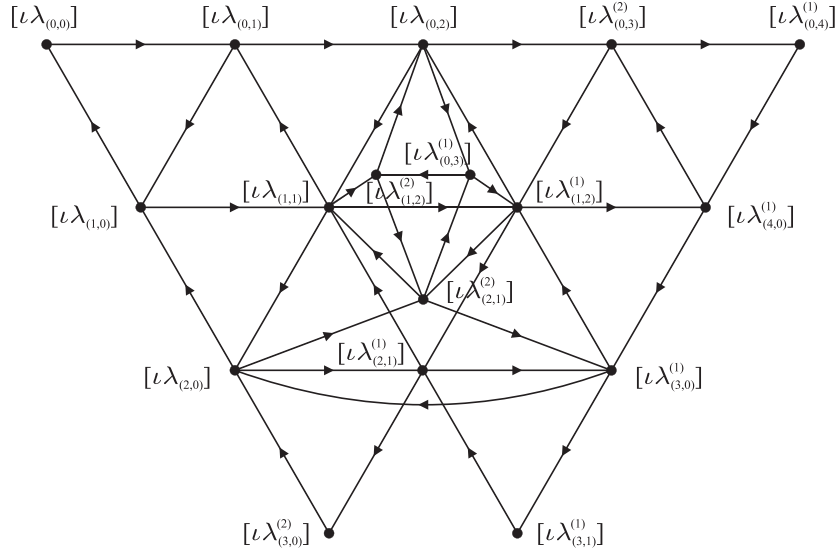}
 \caption{$M$-$N$ graph for the $\mathcal{E}_5^{(12)}$ $SU(3)$-GHJ subfactor}\label{fig:fig-GHJ-11}
\end{center}
\end{figure}

Then the only possibility for the graph of the $M$-$N$ system is $\mathcal{E}_5^{(12)}$, illustrated in Figure \ref{fig:fig-GHJ-11}, and the modular invariant for $\theta$ is $Z_{\mathcal{E}_{MS}^{(12)}}$.

\subsection{$\mathcal{A}^{(n)\ast}$}

We compute the nimrep graph for the case $n=12$. It appears that the results will carry over to all other $n$, however we have not been able to show this in general. For the graph $\mathcal{A}^{(12)\ast}$, we have $[\theta] = \bigoplus_{\mu} [\lambda_{\mu}]$, where the direct sum is over all representations $\mu$ on $\mathcal{A}^{(12)}$. Then computing $\langle \iota \lambda, \iota \mu \rangle = \langle \lambda, \theta \mu \rangle$ for $\lambda$, $\mu$ on $\mathcal{A}^{(12)}$, we find that $\langle \iota \lambda_{(\mu_2, \mu_1)}, \iota \lambda_{(\mu_2, \mu_1)} \rangle = \langle \iota \lambda_{(\mu_2, \mu_1)}, \iota \lambda_{(\mu_1, \mu_2)} \rangle$ so we have $[\iota \lambda_{(\mu_2, \mu_1)}] = [\iota \lambda_{(\mu_1, \mu_2)}]$ for all $(\mu_1, \mu_2)$ on $\mathcal{A}^{(12)}$. At tier 0 we have $\langle \iota \lambda_{(0,0)}, \iota \lambda_{(0,0)} \rangle = 1$. At tier 1, $\langle \iota \lambda_{(1,0)}, \iota \lambda_{(1,0)} \rangle = 2$ and $\langle \iota \lambda_{(1,0)}, \iota \lambda_{(0,0)} \rangle = 1$, giving
\begin{equation} \label{eqn:A(12star)-lambda(1,0)}
[\iota \lambda_{(1,0)}] = [\iota \lambda_{(0,0)}] \oplus [\iota \lambda_{(1,0)}^{(1)}].
\end{equation}
At tier 2 we have $\langle \iota \lambda_{(2,0)}, \iota \lambda_{(2,0)} \rangle = 3$ and $\langle \iota \lambda_{(2,0)}, \iota \lambda_{(1,0)} \rangle = 2$, so $[\iota \lambda_{(2,0)}] = [\iota \lambda_{(0,0)}] \oplus [\iota \lambda_{(1,0)}^{(1)}] \oplus [\iota \lambda_{(2,0)}^{(1)}]$. We also have $\langle \iota \lambda_{(1,1)}, \iota \lambda_{(1,1)} \rangle = 6$, $\langle \iota \lambda_{(1,1)}, \iota \lambda_{(0,0)} \rangle = 1$, $\langle \iota \lambda_{(1,1)}, \iota \lambda_{(1,0)} \rangle = 3$ and $\langle \iota \lambda_{(1,1)}, \iota \lambda_{(2,0)} \rangle = 4$, giving $[\iota \lambda_{(1,1)}] = [\iota \lambda_{(0,0)}] \oplus 2 [\iota \lambda_{(1,0)}^{(1)}] \oplus [\iota \lambda_{(2,0)}^{(1)}]$. At tier 3 we have $\langle \iota \lambda_{(3,0)}, \iota \lambda_{(3,0)} \rangle = 4$ and $\langle \iota \lambda_{(3,0)}, \iota \lambda_{(2,0)} \rangle = 3$, so $[\iota \lambda_{(3,0)}] = [\iota \lambda_{(0,0)}] \oplus [\iota \lambda_{(1,0)}^{(1)}] \oplus [\iota \lambda_{(2,0)}^{(1)}] \oplus [\iota \lambda_{(3,0)}^{(1)}]$. We also have $\langle \iota \lambda_{(2,1)}, \iota \lambda_{(2,1)} \rangle = 10$, $\langle \iota \lambda_{(2,1)}, \iota \lambda_{(0,0)} \rangle = 1$, $\langle \iota \lambda_{(2,1)}, \iota \lambda_{(1,0)} \rangle = 3$, $\langle \iota \lambda_{(2,1)}, \iota \lambda_{(2,0)} \rangle = 5$ and $\langle \iota \lambda_{(2,1)}, \iota \lambda_{(3,0)} \rangle = 6$, giving $[\iota \lambda_{(2,1)}] = [\iota \lambda_{(0,0)}] \oplus 2 [\iota \lambda_{(1,0)}^{(1)}] \oplus 2 [\iota \lambda_{(2,0)}^{(1)}] \oplus [\iota \lambda_{(3,0)}^{(1)}]$. Similarly, at tier 4 we find
\begin{eqnarray*}
{} [\iota \lambda_{(4,0)}] & = & [\iota \lambda_{(0,0)}] \oplus [\iota \lambda_{(1,0)}^{(1)}] \oplus [\iota \lambda_{(2,0)}^{(1)}] \oplus [\iota \lambda_{(3,0)}^{(1)}] \oplus [\iota \lambda_{(4,0)}^{(1)}], \\
{} [\iota \lambda_{(3,1)}] & = & [\iota \lambda_{(0,0)}] \oplus 2 [\iota \lambda_{(1,0)}^{(1)}] \oplus 2 [\iota \lambda_{(2,0)}^{(1)}] \oplus 2 [\iota \lambda_{(3,0)}^{(1)}] \oplus [\iota \lambda_{(4,0)}^{(1)}], \\
{} [\iota \lambda_{(2,2)}] & = & [\iota \lambda_{(0,0)}] \oplus 2 [\iota \lambda_{(1,0)}^{(1)}] \oplus 3 [\iota \lambda_{(2,0)}^{(1)}] \oplus 2 [\iota \lambda_{(3,0)}^{(1)}] \oplus [\iota \lambda_{(4,0)}^{(1)}],
\end{eqnarray*}
and at tier 5:
\begin{eqnarray*}
{} [\iota \lambda_{(5,0)}] & = & [\iota \lambda_{(4,0)}], \\
{} [\iota \lambda_{(4,1)}] & = & [\iota \lambda_{(0,0)}] \oplus 2 [\iota \lambda_{(1,0)}^{(1)}] \oplus 2 [\iota \lambda_{(2,0)}^{(1)}] \oplus 2 [\iota \lambda_{(3,0)}^{(1)}] \oplus 2 [\iota \lambda_{(4,0)}^{(1)}] \oplus [\iota \lambda_{(5,0)}^{(1)}], \\
{} [\iota \lambda_{(3,2)}] & = & [\iota \lambda_{(0,0)}] \oplus 2 [\iota \lambda_{(1,0)}^{(1)}] \oplus 3 [\iota \lambda_{(2,0)}^{(1)}] \oplus 3 [\iota \lambda_{(3,0)}^{(1)}] \oplus 2 [\iota \lambda_{(4,0)}^{(1)}] \oplus [\iota \lambda_{(5,0)}^{(1)}].
\end{eqnarray*}
Then we have six irreducible sectors $[\iota \lambda_{(0,0)}]$, $[\iota \lambda_{(1,0)}^{(1)}]$, $[\iota \lambda_{(2,0)}^{(1)}]$, $[\iota \lambda_{(3,0)}^{(1)}]$, $[\iota \lambda_{(4,0)}^{(1)}]$ and $[\iota \lambda_{(5,0)}^{(1)}]$.

\begin{figure}[tb]
\begin{center}
  \includegraphics[width=55mm]{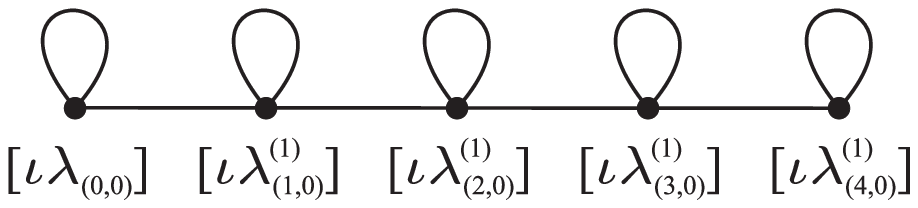}
 \caption{$M$-$N$ graph for the $\mathcal{A}^{(12)\ast}$ $SU(3)$-GHJ subfactor}\label{fig:fig-GHJ-12}
\end{center}
\end{figure}

We now compute the sector products. We have $[\iota \lambda_{(0,0)}] [\rho] = [\iota \lambda_{(1,0)}] = [\iota \lambda_{(0,0)}] \oplus [\iota \lambda_{(1,0)}^{(1)}]$. From $[\iota \lambda_{(1,0)}] [\rho] = [\iota \lambda_{(2,0)}] \oplus [\iota \lambda_{(0,1)}] = 2 [\iota \lambda_{(0,0)}] \oplus 2 [\iota \lambda_{(1,0)}^{(1)}] \oplus [\iota \lambda_{(2,0)}^{(1)}]$ and (\ref{eqn:A(12star)-lambda(1,0)}) we find $[\iota \lambda_{(1,0)}^{(1)}] [\rho] = 2 [\iota \lambda_{(0,0)}] \oplus 2 [\iota \lambda_{(1,0)}^{(1)}] \oplus [\iota \lambda_{(2,0)}^{(1)}] \ominus ([\iota \lambda_{(0,0)}] \oplus [\iota \lambda_{(1,0)}^{(1)}]) = [\iota \lambda_{(0,0)}] \oplus [\iota \lambda_{(1,0)}^{(1)}] \oplus [\iota \lambda_{(2,0)}^{(1)}]$. Similarly, we find
\begin{eqnarray*}
{} [\iota \lambda_{(2,0)}^{(1)}] [\rho] & = & [\iota \lambda_{(1,0)}^{(1)}] \oplus [\iota \lambda_{(2,0)}^{(1)}] \oplus [\iota \lambda_{(3,0)}^{(1)}], \\
{} [\iota \lambda_{(3,0)}^{(1)}] [\rho] & = & [\iota \lambda_{(2,0)}^{(1)}] \oplus [\iota \lambda_{(3,0)}^{(1)}] \oplus [\iota \lambda_{(4,0)}^{(1)}], \\
{} [\iota \lambda_{(4,0)}^{(1)}] [\rho] & = & [\iota \lambda_{(3,0)}^{(1)}] \oplus [\iota \lambda_{(4,0)}^{(1)}],
\end{eqnarray*}
and the nimrep graph is $\mathcal{A}^{(12)\ast}$. The labelled nimrep graph is illustrated in Figure \ref{fig:fig-GHJ-12}. The associated modular invariant is $Z_{\mathcal{A}^{(12)\ast}}$.

In the case above, since $n = 12$ is even, we have $[\iota \lambda_{(5,0)}] = [\iota \lambda_{(4,0)}]$ and so $[\iota \lambda_{(4,0)}] [\rho] = [\iota \lambda_{(5,0)}] \oplus [\iota \lambda_{(3,1)}] = [\iota \lambda_{(4,0)}] \oplus [\iota \lambda_{(3,1)}]$. This leads to $[\iota \lambda_{(4,0)}^{(1)}] [\rho] \supset [\iota \lambda_{(4,0)}^{(1)}]$, and there is a loop from $[\iota \lambda_{(4,0)}^{(1)}]$ to itself in the nimrep graph. However, when $n$ is odd, e.g. for $n = 11$, we have instead $[\iota \lambda_{(5,0)}] = [\iota \lambda_{(3,0)}]$ so $[\iota \lambda_{(4,0)}] [\rho] = [\iota \lambda_{(5,0)}] \oplus [\iota \lambda_{(3,1)}] = [\iota \lambda_{(3,0)}] \oplus [\iota \lambda_{(3,1)}]$. This causes $[\iota \lambda_{(4,0)}^{(1)}] [\rho] \not \supset [\iota \lambda_{(4,0)}^{(1)}]$, hence there is no loop from $[\iota \lambda_{(4,0)}^{(1)}]$ to itself in the nimrep graph for the $n=11$ case.

\subsection{$\mathcal{D}^{(n)\ast}$}

We compute the nimrep graph for the case $n=12$. For the graph $\mathcal{D}^{(12)\ast}$, we have $[\theta] = \bigoplus_{\mu} [\lambda_{\mu}]$, where the direct sum is over all representations $\mu$ of colour 0 on $\mathcal{A}^{(12)}$. At tier 0 we have $\langle \iota \lambda_{(0,0)}, \iota \lambda_{(0,0)} \rangle = 1$. At tier 1, $\langle \iota \lambda_{(1,0)}, \iota \lambda_{(1,0)} \rangle = 2$ and $\langle \iota \lambda_{(1,0)}, \iota \lambda_{(0,0)} \rangle = 0$, and similarly for $\iota \lambda_{(0,1)}$, giving
\begin{eqnarray}
{} [\iota \lambda_{(1,0)}] & = & [\iota \lambda_{(1,0)}^{(1)}] \oplus [\iota \lambda_{(1,0)}^{(2)}], \label{eqn:D(12star)-lambda(1,0)} \\
{} [\iota \lambda_{(0,1)}] & = & [\iota \lambda_{(0,1)}^{(1)}] \oplus [\iota \lambda_{(0,1)}^{(2)}]. \label{eqn:D(12star)-lambda(0,1)}
\end{eqnarray}
At tier 2 we have $\langle \iota \lambda_{(2,0)}, \iota \lambda_{(2,0)} \rangle = 3$ and $\langle \iota \lambda_{(2,0)}, \iota \lambda_{(0,1)} \rangle = 1$, and similarly for $\iota \lambda_{(0,2)}$, so we have
\begin{eqnarray}
{} [\iota \lambda_{(2,0)}] & = & [\iota \lambda_{(0,1)}^{(1)}] \oplus [\iota \lambda_{(0,1)}^{(2)}] \oplus [\iota \lambda_{(2,0)}^{(1)}], \label{eqn:D(12star)-lambda(2,0)} \\
{} [\iota \lambda_{(0,2)}] & = & [\iota \lambda_{(1,0)}^{(1)}] \oplus [\iota \lambda_{(1,0)}^{(2)}] \oplus [\iota \lambda_{(0,2)}^{(1)}]. \label{eqn:D(12star)-lambda(0,2)}
\end{eqnarray}
For $\iota \lambda_{(1,1)}$ we have $\langle \iota \lambda_{(1,1)}, \iota \lambda_{(1,1)} \rangle = 6$ and $\langle \iota \lambda_{(1,1)}, \iota \lambda_{(0,0)} \rangle = 1$, so there are two possibilities for the decomposition of $[\iota \lambda_{(1,1)}]$ as irreducible sectors:
\begin{equation} \label{eqn:D(12star)-lambda(1,1)}
[\iota \lambda_{(1,1)}] = \left\{ \begin{array}{ll}
{} [\iota \lambda_{(0,0)}] \oplus 2 [\iota \lambda_{(1,1)}^{(1)}] \oplus [\iota \lambda_{(1,1)}^{(2)}] & \textrm{case I}, \\
{} [\iota \lambda_{(0,0)}] \oplus [\iota \lambda_{(1,1)}^{(1)}] \oplus [\iota \lambda_{(1,1)}^{(2)}] \oplus [\iota \lambda_{(1,1)}^{(3)}] \oplus [\iota \lambda_{(1,1)}^{(4)}] \oplus [\iota \lambda_{(1,1)}^{(5)}] & \textrm{case II}.
\end{array} \right.
\end{equation}
At tier 3 we have $\langle \iota \lambda_{(3,0)}, \iota \lambda_{(3,0)} \rangle = 4$, $\langle \iota \lambda_{(3,0)}, \iota \lambda_{(1,1)} \rangle = 4$ and $\langle \iota \lambda_{(3,0)}, \iota \lambda_{(0,0)} \rangle = 1$, giving
\begin{equation} \label{eqn:D(12star)-lambda(3,0)}
[\iota \lambda_{(3,0)}] = \left\{ \begin{array}{ll}
{} [\iota \lambda_{(0,0)}] \oplus [\iota \lambda_{(1,1)}^{(1)}] \oplus [\iota \lambda_{(1,1)}^{(2)}] \oplus [\iota \lambda_{(3,0)}^{(1)}] & \textrm{for case I}, \\
{} [\iota \lambda_{(0,0)}] \oplus [\iota \lambda_{(1,1)}^{(1)}] \oplus [\iota \lambda_{(1,1)}^{(2)}] \oplus [\iota \lambda_{(1,1)}^{(3)}] & \textrm{for case II}.
\end{array} \right.
\end{equation}
Then we see that for case II $[\iota \lambda_{(1,1)}] \supset [\iota \lambda_{(3,0)}]$. However, this contradicts the following values of the inner-products at tier 6, $\langle \iota \lambda_{(3,3)}, \iota \lambda_{(1,1)} \rangle = 8$ and $\langle \iota \lambda_{(3,3)}, \iota \lambda_{(3,0)} \rangle = 10$. So we reject case II.

Continuing at tier 3 we have $\langle \iota \lambda_{(0,3)}, \iota \lambda_{(0,3)} \rangle = \langle \iota \lambda_{(0,3)}, \iota \lambda_{(3,0)} \rangle = 4$, so that $[\iota \lambda_{(0,3)}] = [\iota \lambda_{(3,0)}]$. From $\langle \iota \lambda_{(2,1)}, \iota \lambda_{(2,1)} \rangle = 10$, $\langle \iota \lambda_{(2,1)}, \iota \lambda_{(1,0)} \rangle = 3$ and $\langle \iota \lambda_{(2,1)}, \iota \lambda_{(0,2)} \rangle = 5$, and similarly for $\iota \lambda_{(1,2)}$, we have
\begin{eqnarray}
{} [\iota \lambda_{(2,1)}] & = & 2 [\iota \lambda_{(1,0)}^{(1)}] \oplus [\iota \lambda_{(1,0)}^{(2)}] \oplus 2 [\iota \lambda_{(0,2)}^{(1)}] \oplus [\iota \lambda_{(2,1)}^{(1)}], \label{eqn:D(12star)-lambda(2,1)} \\
{} [\iota \lambda_{(1,2)}] & = & 2 [\iota \lambda_{(0,1)}^{(1)}] \oplus [\iota \lambda_{(0,1)}^{(2)}] \oplus 2 [\iota \lambda_{(2,0)}^{(1)}] \oplus [\iota \lambda_{(1,2)}^{(1)}]. \label{eqn:D(12star)-lambda(1,2)}
\end{eqnarray}
Next, at tier 4, we have $\langle \iota \lambda_{(4,0)}, \iota \lambda_{(4,0)} \rangle = 5$, $\langle \iota \lambda_{(4,0)}, \iota \lambda_{(1,0)} \rangle = 2$, $\langle \iota \lambda_{(4,0)}, \iota \lambda_{(0,2)} \rangle = 3$ and $\langle \iota \lambda_{(4,0)}, \iota \lambda_{(2,1)} \rangle = 6$, so there are two possibilities for the decomposition of $[\iota \lambda_{(4,0)}]$, and similarly for $[\iota \lambda_{(0,4)}]$:
\begin{eqnarray}
{} [\iota \lambda_{(4,0)}] & = & \left\{ \begin{array}{ll}
{} [\iota \lambda_{(1,0)}^{(1)}] \oplus [\iota \lambda_{(1,0)}^{(2)}] \oplus [\iota \lambda_{(0,2)}^{(1)}] \oplus [\iota \lambda_{(0,2)}^{(2)}] \oplus [\iota \lambda_{(4,0)}^{(1)}] & \textrm{case } (i), \\
{} 2 [\iota \lambda_{(1,0)}^{(1)}] \oplus [\iota \lambda_{(0,2)}^{(1)}] & \textrm{case } (ii),
\end{array} \right. \label{eqn:D(12star)-lambda(4,0)} \\
{} [\iota \lambda_{(0,4)}] & = & \left\{ \begin{array}{ll}
{} [\iota \lambda_{(0,1)}^{(1)}] \oplus [\iota \lambda_{(0,1)}^{(2)}] \oplus [\iota \lambda_{(2,0)}^{(1)}] \oplus [\iota \lambda_{(2,0)}^{(2)}] \oplus [\iota \lambda_{(0,4)}^{(1)}] & \textrm{case } (i'), \\
{} 2 [\iota \lambda_{(0,1)}^{(1)}] \oplus [\iota \lambda_{(2,0)}^{(1)}] & \textrm{case } (ii').
\end{array} \right. \label{eqn:D(12star)-lambda(0,4)}
\end{eqnarray}
Since $\langle \iota \lambda_{(3,1)}, \iota \lambda_{(3,1)} \rangle = 14$, $\langle \iota \lambda_{(3,1)}, \iota \lambda_{(0,1)} \rangle = 3$, $\langle \iota \lambda_{(3,1)}, \iota \lambda_{(2,0)} \rangle = 5$, $\langle \iota \lambda_{(3,1)}, \iota \lambda_{(1,2)} \rangle = 11$ and $\langle \iota \lambda_{(3,1)}, \iota \lambda_{(0,4)} \rangle = 8$, then
\begin{equation} \label{eqn:D(12star)-lambda(3,1)}
[\iota \lambda_{(3,1)}] = \left\{ \begin{array}{ll}
{} 2 [\iota \lambda_{(0,1)}^{(1)}] \oplus [\iota \lambda_{(0,1)}^{(2)}] \oplus 2 [\iota \lambda_{(2,0)}^{(1)}] \oplus 2 [\iota \lambda_{(1,2)}^{(1)}] \oplus [\iota \lambda_{(0,4)}^{(1)}] & \textrm{for case } (i'), \\
{} 3 [\iota \lambda_{(0,1)}^{(1)}] \oplus 2 [\iota \lambda_{(2,0)}^{(1)}] \oplus [\iota \lambda_{(1,2)}^{(1)}] & \textrm{for case } (ii').
\end{array} \right.
\end{equation}
Similarly, for $[\iota \lambda_{(1,3)}]$,
\begin{equation} \label{eqn:D(12star)-lambda(1,3)}
[\iota \lambda_{(1,3)}] = \left\{ \begin{array}{ll}
{} 2 [\iota \lambda_{(1,0)}^{(1)}] \oplus [\iota \lambda_{(1,0)}^{(2)}] \oplus 2 [\iota \lambda_{(0,2)}^{(1)}] \oplus 2 [\iota \lambda_{(2,1)}^{(1)}] \oplus [\iota \lambda_{(4,0)}^{(1)}] & \textrm{for case } (i), \\
{} 3 [\iota \lambda_{(1,0)}^{(1)}] \oplus 2 [\iota \lambda_{(0,2)}^{(1)}] \oplus [\iota \lambda_{(2,1)}^{(1)}] & \textrm{for case } (ii).
\end{array} \right.
\end{equation}
From $\langle \iota \lambda_{(2,2)}, \iota \lambda_{(2,2)} \rangle = 19$, $\langle \iota \lambda_{(2,2)}, \iota \lambda_{(0,0)} \rangle = 1$, $\langle \iota \lambda_{(2,2)}, \iota \lambda_{(1,1)} \rangle = 8$ and $\langle \iota \lambda_{(2,2)}, \iota \lambda_{(3,0)} \rangle = 8$, we must have
\begin{equation} \label{eqn:D(12star)-lambda(2,2)}
[\iota \lambda_{(1,3)}] = [\iota \lambda_{(0,0)}] \oplus 2 [\iota \lambda_{(1,1)}^{(1)}] \oplus 3 [\iota \lambda_{(1,1)}^{(2)}] \oplus 2 [\iota \lambda_{(3,0)}^{(1)}] \oplus [\iota \lambda_{(2,2)}^{(1)}].
\end{equation}
At tier 5 we have $\langle \iota \lambda_{(5,0)}, \iota \lambda_{(5,0)} \rangle = \langle \iota \lambda_{(5,0)}, \iota \lambda_{(0,4)} \rangle = 5$, giving $[\iota \lambda_{(5,0)}] = [\iota \lambda_{(0,4)}]$, and similarly $[\iota \lambda_{(0,5)}] = [\iota \lambda_{(4,0)}]$. From $\langle \iota \lambda_{(3,2)}, \iota \lambda_{(3,2)} \rangle = 27$, $\langle \iota \lambda_{(3,2)}, \iota \lambda_{(1,0)} \rangle = 3$, $\langle \iota \lambda_{(3,2)}, \iota \lambda_{(0,2)} \rangle = 6$, $\langle \iota \lambda_{(3,2)}, \iota \lambda_{(2,1)} \rangle = 14$ and $\langle \iota \lambda_{(3,2)}, \iota \lambda_{(1,3)} \rangle = 19$ we must have
\begin{equation} \label{eqn:D(12star)-lambda(3,2)}
[\iota \lambda_{(3,2)}] = \left\{ \begin{array}{ll}
{} 2 [\iota \lambda_{(1,0)}^{(1)}] \oplus [\iota \lambda_{(1,0)}^{(2)}] \oplus 3 [\iota \lambda_{(0,2)}^{(1)}] \oplus 3 [\iota \lambda_{(2,1)}^{(1)}] \oplus 2 [\iota \lambda_{(4,0)}^{(1)}] & \textrm{for case } (i), \\
{} 3 [\iota \lambda_{(1,0)}^{(1)}] \oplus 3 [\iota \lambda_{(0,2)}^{(1)}] \oplus 2 [\iota \lambda_{(2,1)}^{(1)}] \oplus 2 [\iota \lambda_{(4,0)}^{(1)}] \oplus [\iota \lambda_{(3,2)}^{(1)}]  & \textrm{for case } (ii).
\end{array} \right.
\end{equation}
However, case $(ii)$ does not satisfy $\langle \iota \lambda_{(3,2)}, \iota \lambda_{(4,0)} \rangle = 11$, and hence we discard it. Similarly we discard case $(ii')$ since no possible decomposition of $[\iota \lambda_{(2,3)}]$ exists for that case. Then we are left with only the one case $(i)(i')$. We have
\begin{equation} \label{eqn:D(12star)-lambda(2,3)}
[\iota \lambda_{(2,3)}] = 2 [\iota \lambda_{(0,1)}^{(1)}] \oplus [\iota \lambda_{(0,1)}^{(2)}] \oplus 3 [\iota \lambda_{(2,0)}^{(1)}] \oplus 3 [\iota \lambda_{(1,2)}^{(1)}] \oplus 2 [\iota \lambda_{(0,4)}^{(1)}].
\end{equation}
From $\langle \iota \lambda_{(4,1)}, \iota \lambda_{(4,1)} \rangle = 17$, $\langle \iota \lambda_{(4,1)}, \iota \lambda_{(0,0)} \rangle = 1$, $\langle \iota \lambda_{(4,1)}, \iota \lambda_{(1,1)} \rangle = 7$, $\langle \iota \lambda_{(4,1)}, \iota \lambda_{(3,0)} \rangle = 7$ and $\langle \iota \lambda_{(4,1)}, \iota \lambda_{(2,2)} \rangle = 17$, we have
\begin{equation} \label{eqn:D(12star)-lambda(4,1)}
[\iota \lambda_{(4,1)}] = [\iota \lambda_{(0,0)}] \oplus 2 [\iota \lambda_{(1,1)}^{(1)}] \oplus 2 [\iota \lambda_{(1,1)}^{(2)}] \oplus 2 [\iota \lambda_{(3,0)}^{(1)}] \oplus 2 [\iota \lambda_{(2,2)}^{(1)}],
\end{equation}
and since $\langle \iota \lambda_{(1,4)}, \iota \lambda_{(1,4)} \rangle = \langle \iota \lambda_{(1,4)}, \iota \lambda_{(4,1)} \rangle = 17$, $[\iota \lambda_{(1,4)}] = [\iota \lambda_{(4,1)}]$. We see that no new irreducible sectors appear at tier 5, so the $M$-$N$ system contains 15 irreducible sectors.
We also have the following decompositions at tier 6:
\begin{eqnarray}
{} [\iota \lambda_{(6,0)}] & = & [\iota \lambda_{(0,6)}] \;\; = \;\; [\iota \lambda_{(3,0)}], \label{eqn:D(12star)-lambda(6,0)} \\
{} [\iota \lambda_{(5,1)}] & = & 2 [\iota \lambda_{(1,0)}^{(1)}] \oplus [\iota \lambda_{(1,0)}^{(2)}] \oplus 2 [\iota \lambda_{(0,2)}^{(1)}] \oplus 2 [\iota \lambda_{(2,1 )}^{(1)}] \oplus [\iota \lambda_{(4,0)}^{(1)}], \label{eqn:D(12star)-lambda(5,1)} \\
{} [\iota \lambda_{(4,2)}] & = & 2 [\iota \lambda_{(0,1)}^{(1)}] \oplus [\iota \lambda_{(0,1)}^{(2)}] \oplus 3 [\iota \lambda_{(2,0)}^{(1)}] \oplus 3 [\iota \lambda_{(1,2)}^{(1)}] \oplus 2 [\iota \lambda_{(0,4)}^{(1)}], \label{eqn:D(12star)-lambda(4,2)} \\
{} [\iota \lambda_{(1,5)}] & = & 2 [\iota \lambda_{(0,1)}^{(1)}] \oplus [\iota \lambda_{(0,1)}^{(2)}] \oplus 2 [\iota \lambda_{(2,0)}^{(1)}] \oplus 2 [\iota \lambda_{(1,2)}^{(1)}] \oplus [\iota \lambda_{(0,4)}^{(1)}]. \label{eqn:D(12star)-lambda(1,5)}
\end{eqnarray}

We now find the sector products of the irreducible sectors with the $N$-$N$ sector $[\rho] = [\lambda_{(1,0)}]$. We have $[\iota \lambda_{(0,0)}] [\rho] = [\iota \lambda_{(1,0)}] = [\iota \lambda_{(1,0)}^{(1)}] \oplus [\iota \lambda_{(1,0)}^{(2)}]$. From $[\iota \lambda_{(1,1)}] [\rho] = [\iota \lambda_{(1,0)}] \oplus [\iota \lambda_{(0,2)}] \oplus [\iota \lambda_{(2,1)}]$ and (\ref{eqn:D(12star)-lambda(1,1)}) we have
\begin{eqnarray}
(2 [\iota \lambda_{(1,1)}^{(1)}] \oplus [\iota \lambda_{(1,1)}^{(2)}]) [\rho] & = & 4 [\iota \lambda_{(1,0)}^{(1)}] \oplus 3 [\iota \lambda_{(1,0)}^{(2)}] \oplus 3 [\iota \lambda_{(0,2)}^{(1)}] \oplus [\iota \lambda_{(2,1)}^{(1)}] \ominus ([\iota \lambda_{(0,0)}][\rho]) \nonumber \\
& = & 3 [\iota \lambda_{(1,0)}^{(1)}] \oplus 2 [\iota \lambda_{(1,0)}^{(2)}] \oplus 3 [\iota \lambda_{(0,2)}^{(1)}] \oplus [\iota \lambda_{(2,1)}^{(1)}]. \label{eqn:D(12star)-rho_x_lambda(1,1)}
\end{eqnarray}
Similarly, by considering $[\iota \lambda_{(3,0)}] [\rho]$, $[\iota \lambda_{(2,2)}] [\rho]$ and $[\iota \lambda_{(4,1)}] [\rho]$, and using (\ref{eqn:D(12star)-lambda(3,0)}), (\ref{eqn:D(12star)-lambda(2,2)}) and (\ref{eqn:D(12star)-lambda(4,1)}), we have the following:
\begin{eqnarray}
([\iota \lambda_{(1,1)}^{(1)}] \oplus [\iota \lambda_{(1,1)}^{(2)}] \oplus [\iota \lambda_{(3,0)}^{(1)}]) [\rho] & = & 2 [\iota \lambda_{(1,0)}^{(1)}] \oplus [\iota \lambda_{(1,0)}^{(2)}] \oplus 3 [\iota \lambda_{(0,2)}^{(1)}] \nonumber \\
& & \; \oplus 2 [\iota \lambda_{(2,1)}^{(1)}] \oplus [\iota \lambda_{(4,0)}^{(1)}], \label{eqn:D(12star)-rho_x_lambda(3,0)} \\
(2 [\iota \lambda_{(1,1)}^{(1)}] \oplus 3 [\iota \lambda_{(1,1)}^{(2)}] \oplus 2 [\iota \lambda_{(3,0)}^{(1)}] \oplus [\iota \lambda_{(2,2)}^{(1)}]) [\rho] & = & 5 [\iota \lambda_{(1,0)}^{(1)}] \oplus 2 [\iota \lambda_{(1,0)}^{(2)}] \oplus 7 [\iota \lambda_{(0,2)}^{(1)}] \nonumber \\
& & \; \oplus 6 [\iota \lambda_{(2,1)}^{(1)}] \oplus 3 [\iota \lambda_{(4,0)}^{(1)}], \label{eqn:D(12star)-rho_x_lambda(2,2)} \\
(2 [\iota \lambda_{(1,1)}^{(1)}] \oplus 2 [\iota \lambda_{(1,1)}^{(2)}] \oplus 2 [\iota \lambda_{(3,0)}^{(1)}] \oplus 2 [\iota \lambda_{(2,2)}^{(1)}]) [\rho] & = & 4 [\iota \lambda_{(1,0)}^{(1)}] \oplus 2 [\iota \lambda_{(1,0)}^{(2)}] \oplus 6 [\iota \lambda_{(0,2)}^{(1)}] \nonumber \\
& & \; \oplus 6 [\iota \lambda_{(2,1)}^{(1)}] \oplus 4 [\iota \lambda_{(4,0)}^{(1)}]. \label{eqn:D(12star)-rho_x_lambda(4,1)}
\end{eqnarray}
Then from (\ref{eqn:D(12star)-rho_x_lambda(1,1)})-(\ref{eqn:D(12star)-rho_x_lambda(4,1)}) we obtain the following sector products:
\begin{eqnarray*}
{} [\iota \lambda_{(1,1)}^{(1)}] [\rho] & = & [\iota \lambda_{(1,0)}^{(1)}] \oplus [\iota \lambda_{(1,0)}^{(2)}] \oplus [\iota \lambda_{(0,2)}^{(1)}], \\
{} [\iota \lambda_{(1,1)}^{(2)}] [\rho] & = & [\iota \lambda_{(1,0)}^{(1)}] \oplus [\iota \lambda_{(0,2)}^{(1)}] \oplus [\iota \lambda_{(2,1)}^{(1)}], \\
{} [\iota \lambda_{(3,0)}^{(1)}] [\rho] & = & [\iota \lambda_{(0,2)}^{(1)}] \oplus [\iota \lambda_{(2,1)}^{(1)}] \oplus [\iota \lambda_{(4,0)}^{(1)}], \\
{} [\iota \lambda_{(2,2)}^{(1)}] [\rho] & = & [\iota \lambda_{(2,1)}^{(1)}] \oplus [\iota \lambda_{(4,0)}^{(1)}].
\end{eqnarray*}

Next, from $[\iota \lambda_{(1,0)}] [\rho] = [\iota \lambda_{(0,1)}] \oplus [\iota \lambda_{(2,0)}]$ and (\ref{eqn:D(12star)-lambda(1,0)}) we have
\begin{equation} \label{eqn:D(12star)-rho_x_lambda(1,0)}
([\iota \lambda_{(1,0)}^{(1)}] \oplus [\iota \lambda_{(1,0)}^{(2)}]) [\rho] = 2 [\iota \lambda_{(0,1)}^{(1)}] \oplus 2 [\iota \lambda_{(0,1)}^{(2)}] \oplus [\iota \lambda_{(2,0)}^{(1)}].
\end{equation}
By considering $[\iota \lambda_{(0,2)}] [\rho] = [\iota \lambda_{(0,1)}] \oplus [\iota \lambda_{(1,2)}]$ and (\ref{eqn:D(12star)-lambda(0,2)}) we obtain $([\iota \lambda_{(1,0)}^{(1)}] \oplus [\iota \lambda_{(1,0)}^{(2)}] \oplus [\iota \lambda_{(0,2)}^{(1)}]) [\rho] = 3 [\iota \lambda_{(0,1)}^{(1)}] \oplus 2 [\iota \lambda_{(0,1)}^{(2)}] \oplus 2 [\iota \lambda_{(2,0)}^{(1)}] \oplus [\iota \lambda_{(1,2)}^{(1)}]$. Then from (\ref{eqn:D(12star)-rho_x_lambda(1,0)}) we see that
\begin{equation} \label{eqn:D(12star)-rho_x_lambda(0,2)(1)}
[\iota \lambda_{(0,2)}^{(1)}] [\rho] = [\iota \lambda_{(0,1)}^{(1)}] \oplus [\iota \lambda_{(2,0)}^{(1)}] \oplus [\iota \lambda_{(1,2)}^{(1)}].
\end{equation}
From $[\iota \lambda_{(2,1)}] [\rho]$, (\ref{eqn:D(12star)-lambda(2,1)}) and (\ref{eqn:D(12star)-rho_x_lambda(0,2)(1)}) we find
\begin{equation} \label{eqn:D(12star)-rho_x_lambda(2,1)}
(2 [\iota \lambda_{(1,0)}^{(1)}] \oplus [\iota \lambda_{(1,0)}^{(2)}] \oplus [\iota \lambda_{(2,1)}^{(1)}]) [\rho] = 3 [\iota \lambda_{(0,1)}^{(1)}] \oplus 3 [\iota \lambda_{(0,1)}^{(2)}] \oplus 3 [\iota \lambda_{(2,0)}^{(1)}] \oplus [\iota \lambda_{(1,2)}^{(1)}] \oplus [\iota \lambda_{(0,4)}^{(1)}].
\end{equation}
Similarly, by considering $[\iota \lambda_{(1,3)}] [\rho]$ and $[\iota \lambda_{(0,5)}] [\rho]$, and using (\ref{eqn:D(12star)-lambda(4,0)}), (\ref{eqn:D(12star)-lambda(1,3)}) and $[\iota \lambda_{(0,5)}] = [\iota \lambda_{(4,0)}]$, we have the following:
\begin{eqnarray}
(2 [\iota \lambda_{(1,0)}^{(1)}] \oplus [\iota \lambda_{(1,0)}^{(2)}] \oplus 2 [\iota \lambda_{(2,1)}^{(1)}] \oplus [\iota \lambda_{(4,0)}^{(1)}]) [\rho] & = & 3 [\iota \lambda_{(0,1)}^{(1)}] \oplus 3 [\iota \lambda_{(0,1)}^{(2)}] \oplus 4 [\iota \lambda_{(2,0)}^{(1)}] \nonumber \\
& & \; \oplus 3 [\iota \lambda_{(1,2)}^{(1)}] \oplus 2 [\iota \lambda_{(0,4)}^{(1)}], \label{eqn:D(12star)-rho_x_lambda(1,3)} \\
([\iota \lambda_{(1,0)}^{(1)}] \oplus [\iota \lambda_{(1,0)}^{(2)}] \oplus [\iota \lambda_{(2,1)}^{(1)}] \oplus [\iota \lambda_{(4,0)}^{(1)}]) [\rho] & = & 2 [\iota \lambda_{(0,1)}^{(1)}] \oplus 2 [\iota \lambda_{(0,1)}^{(2)}] \oplus 2 [\iota \lambda_{(2,0)}^{(1)}] \nonumber \\
& & \; \oplus 2 [\iota \lambda_{(1,2)}^{(1)}] \oplus 2 [\iota \lambda_{(0,4)}^{(1)}]. \label{eqn:D(12star)-rho_x_lambda(0,5)}
\end{eqnarray}
Then from (\ref{eqn:D(12star)-rho_x_lambda(1,0)}), (\ref{eqn:D(12star)-rho_x_lambda(2,1)})-(\ref{eqn:D(12star)-rho_x_lambda(0,5)}) we obtain the following sector products:
\begin{eqnarray*}
{} [\iota \lambda_{(1,0)}^{(1)}] [\rho] & = & [\iota \lambda_{(0,1)}^{(1)}] \oplus [\iota \lambda_{(0,1)}^{(2)}] \oplus [\iota \lambda_{(2,0)}^{(1)}], \\
{} [\iota \lambda_{(1,0)}^{(2)}] [\rho] & = & [\iota \lambda_{(0,1)}^{(1)}] \oplus [\iota \lambda_{(0,1)}^{(2)}], \\
{} [\iota \lambda_{(2,1)}^{(1)}] [\rho] & = & [\iota \lambda_{(2,0)}^{(1)}] \oplus [\iota \lambda_{(1,2)}^{(1)}] \oplus [\iota \lambda_{(0,4)}^{(1)}], \\
{} [\iota \lambda_{(4,0)}^{(1)}] [\rho] & = & [\iota \lambda_{(1,2)}^{(1)}] \oplus [\iota \lambda_{(0,4)}^{(1)}].
\end{eqnarray*}

Next, since $[\iota \lambda_{(0,1)}] [\rho] = [\iota \lambda_{(0,0)}] \oplus [\iota \lambda_{(1,1)}]$, from (\ref{eqn:D(12star)-lambda(0,1)}) we have
\begin{equation} \label{eqn:D(12star)-rho_x_lambda(0,1)}
([\iota \lambda_{(0,1)}^{(1)}] \oplus [\iota \lambda_{(0,1)}^{(2)}]) [\rho] = 2 [\iota \lambda_{(0,0)}] \oplus 2 [\iota \lambda_{(1,1)}^{(1)}] \oplus [\iota \lambda_{(1,1)}^{(2)}].
\end{equation}
By considering $[\iota \lambda_{(2,0)}] [\rho] = [\iota \lambda_{(1,1)}] \oplus [\iota \lambda_{(3,0)}]$ and (\ref{eqn:D(12star)-lambda(2,0)}) we obtain $([\iota \lambda_{(0,1)}^{(1)}] \oplus [\iota \lambda_{(0,1)}^{(2)}] \oplus [\iota \lambda_{(2,0)}^{(1)}]) [\rho] = 2 [\iota \lambda_{(0,0)}] \oplus 3 [\iota \lambda_{(1,1)}^{(1)}] \oplus 2 [\iota \lambda_{(1,1)}^{(2)}] \oplus [\iota \lambda_{(3,0)}^{(1)}]$. Then from (\ref{eqn:D(12star)-rho_x_lambda(0,1)}) we see that
\begin{equation} \label{eqn:D(12star)-rho_x_lambda(2,0)(1)}
[\iota \lambda_{(2,0)}^{(1)}] [\rho] = [\iota \lambda_{(1,1)}^{(1)}] \oplus [\iota \lambda_{(1,1)}^{(2)}] \oplus [\iota \lambda_{(3,0)}^{(1)}].
\end{equation}
From $[\iota \lambda_{(1,2)}] [\rho]$, (\ref{eqn:D(12star)-lambda(1,2)}) and (\ref{eqn:D(12star)-rho_x_lambda(2,0)(1)}) we obtain
\begin{equation} \label{eqn:D(12star)-rho_x_lambda(1,2)}
(2 [\iota \lambda_{(0,1)}^{(1)}] \oplus [\iota \lambda_{(0,1)}^{(2)}] \oplus [\iota \lambda_{(1,2)}^{(1)}]) [\rho] = 3 [\iota \lambda_{(0,0)}] \oplus 3 [\iota \lambda_{(1,1)}^{(1)}] \oplus 3 [\iota \lambda_{(1,1)}^{(2)}] \oplus [\iota \lambda_{(3,0)}^{(1)}] \oplus [\iota \lambda_{(2,2)}^{(1)}].
\end{equation}
Similarly, by considering $[\iota \lambda_{(3,1)}] [\rho]$ and $[\iota \lambda_{(0,4)}] [\rho]$, and using (\ref{eqn:D(12star)-lambda(3,1)}) and (\ref{eqn:D(12star)-lambda(0,4)}), we have the following:
\begin{eqnarray}
(2 [\iota \lambda_{(0,1)}^{(1)}] \oplus [\iota \lambda_{(0,1)}^{(2)}] \oplus 2 [\iota \lambda_{(1,2)}^{(1)}] \oplus [\iota \lambda_{(0,4)}^{(1)}]) [\rho] & = & 3 [\iota \lambda_{(0,0)}] \oplus 3 [\iota \lambda_{(1,1)}^{(1)}] \oplus 4 [\iota \lambda_{(1,1)}^{(2)}] \nonumber \\
& & \; \oplus 3 [\iota \lambda_{(3,0)}^{(1)}] \oplus 3 [\iota \lambda_{(2,2)}^{(1)}], \label{eqn:D(12star)-rho_x_lambda(3,1)} \\
([\iota \lambda_{(0,1)}^{(1)}] \oplus [\iota \lambda_{(0,1)}^{(2)}] \oplus [\iota \lambda_{(1,2)}^{(1)}] \oplus [\iota \lambda_{(0,4)}^{(1)}]) [\rho] & = & 2 [\iota \lambda_{(0,0)}] \oplus 2 [\iota \lambda_{(1,1)}^{(1)}] \oplus 2 [\iota \lambda_{(1,1)}^{(2)}] \nonumber \\
& & \; \oplus 2 [\iota \lambda_{(3,0)}^{(1)}] \oplus 2 [\iota \lambda_{(2,2)}^{(1)}]. \label{eqn:D(12star)-rho_x_lambda(0,4)}
\end{eqnarray}
Then from (\ref{eqn:D(12star)-rho_x_lambda(0,1)}), (\ref{eqn:D(12star)-rho_x_lambda(1,2)})-(\ref{eqn:D(12star)-rho_x_lambda(0,4)}) we obtain the following sector products:
\begin{eqnarray*}
{} [\iota \lambda_{(0,1)}^{(1)}] [\rho] & = & [\iota \lambda_{(0,0)}] \oplus [\iota \lambda_{(1,1)}^{(1)}] \oplus [\iota \lambda_{(1,1)}^{(2)}], \\
{} [\iota \lambda_{(0,1)}^{(2)}] [\rho] & = & [\iota \lambda_{(0,0)}] \oplus [\iota \lambda_{(1,1)}^{(1)}], \\
{} [\iota \lambda_{(1,2)}^{(1)}] [\rho] & = & [\iota \lambda_{(1,1)}^{(2)}] \oplus [\iota \lambda_{(3,0)}^{(1)}] \oplus [\iota \lambda_{(2,2)}^{(1)}], \\
{} [\iota \lambda_{(0,4)}^{(1)}] [\rho] & = & [\iota \lambda_{(3,0)}^{(1)}] \oplus [\iota \lambda_{(2,2)}^{(1)}].
\end{eqnarray*}

\begin{figure}[tb]
\begin{center}
  \includegraphics[width=130mm]{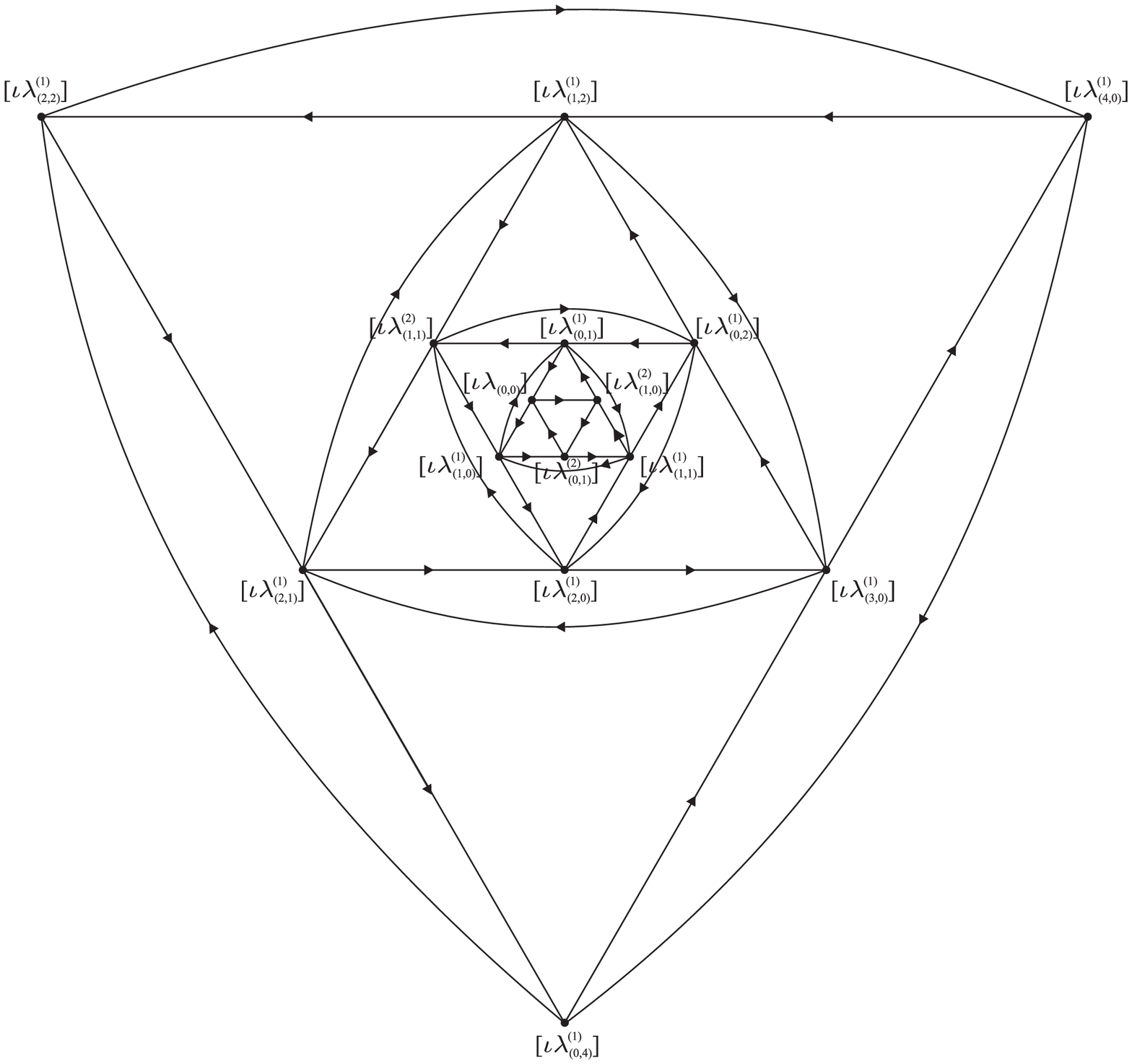}
 \caption{$M$-$N$ graph for the $\mathcal{D}^{(12)\ast}$ $SU(3)$-GHJ subfactor}\label{fig:fig-GHJ-13}
\end{center}
\end{figure}

We thus obtain the graph $\mathcal{D}^{(12)\ast}$ as the nimrep graph for the $M$-$N$ system, illustrated in Figure \ref{fig:fig-GHJ-13}, and the associated modular invariant is $Z_{\mathcal{D}^{(12)\ast}}$.

\subsection{The type I parent} \label{Sect:TypeIparent}

Thus we have constructed subfactors which realize all of the $SU(3)$ modular invariants, except for the $\mathcal{E}_4^{(12)}$ case, since the existence of this subfactor is not yet shown.
However, for the modular invariant associated to the graph $\mathcal{E}_4^{(12)}$, we have $Z_{\mathcal{E}_{MS}^{(12)\ast}} = Z_{\mathcal{E}_{MS}^{(12)}} C$, where $C$ is the modular invariant associated to the graph $\mathcal{A}^{(12)\ast}$. Since both $Z_{\mathcal{E}_{MS}^{(12)}}$, $C$ are shown to be realised by subfactors, the result of \cite[Theorem 3.6]{evans/pinto:2003} shows that the modular invariant $Z_{\mathcal{E}_{MS}^{(12)\ast}}$ is also realised by a subfactor.

The $M$-$N$ graph $\mathcal{G}$ of a subfactor $N \subset M$ is defined by the matrix $\Delta_{\rho}$ which gives the decomposition of the $M$-$N$ sectors with respect to multiplication by the fundamental representation $\rho$. Similarly, multiplication by the conjugate representation defines the matrix $\Delta_{\overline{\rho}} = \Delta_{\rho}^T$ which is the adjacency matrix of the conjugate graph $\widetilde{\mathcal{G}}$. Then since ${}_N \mathcal{X}_N$ is commutative, the matrices $\Delta_{\rho}$ and $\Delta_{\rho}^T$ commute, i.e. $\Delta_{\rho}$ is normal. This provides a proof that the adjacency matrices of the $\mathcal{ADE}$ graphs are all normal, since each of the $\mathcal{ADE}$ graphs appears as the $M$-$N$ graph for a subfactor $N \subset M$.

The zero-column of the modular invariant $Z$ associated with the subfactor $N \subset M$ determines $\langle \alpha^+_j, \alpha^+_{j'} \rangle$ since $\alpha$ preserves the sector product
\begin{eqnarray}
\langle \alpha^+_j, \alpha^+_{j'} \rangle & = & \langle \alpha^+_j \alpha^+_{j'}, \mathrm{id} \rangle = \sum_{j''} N_{j,j'}^{j''} \langle \alpha^+_{j''}, \mathrm{id} \rangle \nonumber \\
& = & \sum_{j''} N_{j,j'}^{j''} Z_{j'',0}, \label{eqn:alpha-graph-Z}
\end{eqnarray}
and similarly the zero-row determines $\langle \alpha^-_j, \alpha^-_{j'} \rangle$. Then for all modular invariants with the same zero-column, the sectors $[\alpha^{\pm}_1]$ satisfy the same equation (\ref{eqn:alpha-graph-Z}) and hence have the same nimrep graphs. Let $v$ be an isometry which intertwines the identity and the canonical endomorphism $\gamma = \iota \overline{\iota}$. Proposition 3.2 in \cite{bockenhauer/evans:2000} states that the following conditions are equivalent:
\begin{itemize}
\item[1.] $Z_{\lambda,0} = \langle \theta, \lambda \rangle$ for all $\lambda \in {}_N \mathcal{X}_N$.
\item[2.] $Z_{0,\lambda} = \langle \theta, \lambda \rangle$ for all $\lambda \in {}_N \mathcal{X}_N$.
\item[3.] Chiral locality holds: $\varepsilon^+(\theta,\theta)v^2 = v^2$.
\end{itemize}
The chiral locality condition, which can be expressed in terms of the single inclusion $N \subset M$ and the braiding, expresses local commutativity (locality) of the extended net, if $N \subset M$ arises from a net of subfactors \cite{longo/rehren:1995}. Chiral locality holds if and only if the dual canonical endomorphism is visible in the vacuum row, $[\theta] = \bigoplus_{\lambda} Z_{0,\lambda} [\lambda]$ (and hence in the vacuum column also).

We will call the inclusion $N \subset M$ \textbf{type I} if and only if one of the above equivalent conditions 1-3 hold. Otherwise we will call the inclusion \textbf{type II}.
Note that the inclusions obtained for the $\mathcal{E}_1^{(12)}$ and $\mathcal{E}_2^{(12)}$ graphs realize the same modular invariant $Z_{\mathcal{E}^{(12)}}$, but the inclusion for $\mathcal{E}_1^{(12)}$ is type I whilst the inclusion for $\mathcal{E}_2^{(12)}$ is type II. This shows that it is possible for a type I modular invariant to be realized by a type II inclusion, and suggests that care needs to be taken with the type I, II labelling of modular invariants. The nimrep graph of $[\alpha^{\pm}_1]$ for the identity modular invariant is the fusion graph of the original $N$-$N$ system, whilst the nimrep graph of $[\alpha^{\pm}_1]$ for the modular invariants associated to $\mathcal{D}^{(3k+3)}$ and $\mathcal{E}^{(8)}$ were computed in \cite{bockenhauer/evans:1999i}, and for $\mathcal{E}_1^{(12)}$ and $\mathcal{E}^{(24)}$ in \cite{bockenhauer/evans:1999ii}. In these cases we have $Z_{\lambda,0} = \langle \theta, \lambda \rangle$ for all $\lambda \in {}_N \mathcal{X}_N$, for $\theta$ given in (\ref{theta-A(n)})-(\ref{theta-E(24)}).
The principal graph of the inclusion $\alpha^{\pm}_1(N) \subset N$ is then the nimrep graph of $[\alpha^{\pm}_1]$. The other modular invariants all have the same zero-column as one of these modular invariants, and hence the nimrep graph of $[\alpha^{\pm}_1]$ for these modular invariants must be the graph given by the type I parent of $Z$, that is, the type I modular invariant which has the same first column as $Z$.
The results are summarized in Table \ref{table:summary_of_M-N_system}, where ``Type'' refers to the type of the inclusion $N \subset M$ given by the $SU(3)$-GHJ construction, where the distinguished vertex $\ast_{\mathcal{G}}$ is the vertex with lowest Perron-Frobenius weight.\footnote{Note, we have only showed the $\mathcal{A}^{\ast}$ and $\mathcal{D}^{\ast}$ case for $n=12$. We have not done any computations for the $D^{(n)}$ graphs, $n \not \equiv 0 \textrm{ mod } 3$.}

\begin{table}[hbt]
\begin{center}
\begin{tabular}{|c|c|c|c|c|} \hline
GHJ graph & Modular invariant & Type & $M$-$N$ graph & Type I parent \\
\hline\hline $\mathcal{A}^{(n)}$ & $Z_{\mathcal{A}^{(n)}}$ & I & $\mathcal{A}^{(n)}$ & $\mathcal{A}^{(n)}$ \\
\hline $\mathcal{A}^{(n)\ast}$ & $Z_{\mathcal{A}^{(n)\ast}} = C$ & II & $\mathcal{A}^{(n)\ast}$ & $\mathcal{A}^{(n)}$ \\
\hline $\mathcal{D}^{(3k)}$ & $Z_{\mathcal{D}^{(3k)}}$ & I & $\mathcal{D}^{(3k)}$ & $\mathcal{D}^{(3k)}$ \\
\hline $\mathcal{D}^{(n)}$ \quad ($n \not \equiv 0 \textrm{ mod } 3$) & $Z_{\mathcal{D}^{(n)}}$ & II & ? & $\mathcal{A}^{(n)}$ \\
\hline $\mathcal{D}^{(3k)\ast}$ & $Z_{\mathcal{D}^{(3k)\ast}} = Z_{\mathcal{D}^{(3k)}} C$ & II & $\mathcal{D}^{(3k)\ast}$ & $\mathcal{D}^{(3k)}$ \\
\hline $\mathcal{D}^{(n)\ast}$ \quad ($n \not \equiv 0 \textrm{ mod } 3$) & $Z_{\mathcal{D}^{(n)\ast}} = Z_{\mathcal{D}^{(n)}} C$ & II & $\mathcal{D}^{(n)\ast}$ & $\mathcal{A}^{(n)}$ \\
\hline $\mathcal{E}^{(8)}$ & $Z_{\mathcal{E}^{(8)}}$ & I & $\mathcal{E}^{(8)}$ & $\mathcal{E}^{(8)}$ \\
\hline $\mathcal{E}^{(8)\ast}$ & $Z_{\mathcal{E}^{(8)\ast}} = Z_{\mathcal{E}^{(8)}} C$ & II & $\mathcal{E}^{(8)\ast}$ & $\mathcal{E}^{(8)}$ \\
\hline $\mathcal{E}_1^{(12)}$ & $Z_{\mathcal{E}^{(12)}} = Z_{\mathcal{E}^{(12)}} C$ & I & $\mathcal{E}_1^{(12)}$ & $\mathcal{E}^{(12)}_1$ \\
\hline $\mathcal{E}_2^{(12)}$ & $Z_{\mathcal{E}^{(12)}} = Z_{\mathcal{E}^{(12)}} C$ & II & $\mathcal{E}_2^{(12)}$ & $\mathcal{E}^{(12)}_1$ \\
\hline $\mathcal{E}_3^{(12)}$ & - & - & - & - \\
\hline $\mathcal{E}_4^{(12)}$ & $Z_{\mathcal{E}_{MS}^{(12)\ast}} = Z_{\mathcal{E}_{MS}^{(12)}} C$ & II & $\mathcal{E}^{(12)}_4$ & $\mathcal{D}^{(12)}$ \\
\hline $\mathcal{E}_5^{(12)}$ & $Z_{\mathcal{E}_{MS}^{(12)}}$ & II & $\mathcal{E}_5^{(12)}$ & $\mathcal{D}^{(12)}$ \\
\hline $\mathcal{E}^{(24)}$ & $Z_{\mathcal{E}^{(24)}} = Z_{\mathcal{E}^{(24)}} C$ & I & $\mathcal{E}^{(24)}$ & $\mathcal{E}^{(24)}$ \\
\hline
\end{tabular}\\
\caption{The $SU(3)$ modular invariants realized by $SU(3)$-GHJ subfactors} \label{table:summary_of_M-N_system}
\end{center}
\end{table}

For $\mathcal{E}_4^{(12)}$, we do not show the existence of the Ocneanu cells, and hence do not have a GHJ subfactor here. However, we have shown that the $Z_{\mathcal{E}_{MS}^{(12)\ast}}$ modular invariant is realised as a braided subfactor. The corresponding nimrep is not computed here, but if (\ref{theta-E4(12)}) is a dual canonical endomorphism, then its nimrep graph is shown to be $\mathcal{E}_4^{(12)}$. This would be the case if $\mathcal{E}_4^{(12)}$ carries a cell system.

\paragraph{Acknowledgements}

This paper is based on work in \cite{pugh:2008}. The first author was partially supported by the EU-NCG network in Non-Commutative Geometry MRTN-CT-2006-031962, and the second author was supported by a scholarship from the School of Mathematics, Cardiff University.

\end{document}